\documentclass[12pt,a4paper]{article}
\newcommand{\R}{\mathbb{R}}

\def\build#1_#2^#3{\mathrel{
\mathop{\kern 0pt#1}\limits_{#2}^{#3}}}

\def\cq{$\hfill \square$}

\def\d{{\rm d}}

\def\eps{\varepsilon}

\def\ba{\begin{eqnarray*}}
\def\ea{\end{eqnarray*}}

\def\dim{{\rm dim\,}}

\usepackage{amsmath}
\usepackage{amsfonts}
\usepackage{amssymb}
\usepackage{amsthm}
\usepackage{graphics}
\begin{document}

\title{The genealogy of self-similar fragmentations with negative index 
as a continuum random tree}
\author{B\'en\'edicte Haas\thanks{%
LPMA, Universit\'e Pierre et Marie Curie, 175 rue du Chevaleret,
F-75013
Paris. \texttt{haas@ccr.jussieu.fr}}\ \ \& Gr\'egory Miermont\thanks{%
DMA, \'Ecole Normale Sup\'erieure, \& LPMA, Universit\'e Pierre et
Marie Curie, 45, rue d'Ulm, 75230 Paris Cedex 05.
\texttt{miermont@dma.ens.fr}}}
\date{}
\maketitle

\begin{abstract}
We encode a certain class of stochastic fragmentation processes,
namely self-similar fragmentation processes with a negative index
of self-similarity, into a metric family tree which belongs to the
family of Continuum Random Trees of Aldous. When the splitting
times of the fragmentation are dense near 0, the tree can in turn
be encoded into a continuous height function, just as the Brownian
Continuum Random Tree is encoded in a normalized Brownian
excursion. Under mild hypotheses, we then compute the Hausdorff
dimensions of these trees, and the maximal H\"older exponents of
the height functions.
\end{abstract}


\bigskip
\noindent{\bf A.M.S. Classification. }60G18, 60J25, 60G09.

\newtheorem{defn}{Definition} \newtheorem{prp}{Proposition} \newtheorem{lmm}{%
Lemma} \newtheorem{thm}{Theorem} \newtheorem{crl}{Corollary}

\newpage

\section{Introduction}

Self-similar fragmentation processes describe the evolution of an object
that falls apart, so that different fragments keep on collapsing
independently with a rate that depends on their sizes to a certain power,
called the \emph{index} of the self-similar fragmentation. A genealogy is
naturally associated to such fragmentation processes, by saying that the
common ancestor of two fragments is the block that included these fragments
for the last time, before a dislocation had definitely separated them. With
an appropriate coding of the fragments, one guesses that there should be a
natural way to define a genealogy tree, rooted at the initial fragment,
associated to any such fragmentation. It would be natural to put a metric on
this tree, e.g.\ by letting the distance from a fragment to the root of the
tree be the time at which the fragment disappears.

Conversely, it turns out that trees have played a key role in models
involving self-similar fragmentations, notably, Aldous and Pitman \cite
{jpda98sac} have introduced a way to log the so-called Brownian \emph{%
Continuum Random Tree} (CRT) \cite{aldouscrt93} that is related to the
standard additive coalescent. Bertoin \cite{bertsfrag02} has shown that a
fragmentation that is somehow dual to the Aldous-Pitman fragmentation can be
obtained as follows. Let $\mathcal{T }_{B}$ be the Brownian CRT, which is
considered as an ``infinite tree with edge-lengths'' (formal definitions are
given below). Let $\mathcal{T }_{t}^{1}, \mathcal{T }_{t}^{2},\ldots $ be
the distinct tree components of the forest obtained by removing all the
vertices of $\mathcal{T }$ that are at distance less than $t$ from the root,
and arranged by decreasing order of ``size''. Then the sequence $F_{B}(t)$
of these sizes defines as $t$ varies a self-similar fragmentation. A moment
of thought points out that the notion of genealogy defined above precisely
coincides with the tree we have fragmented in this way, since a split occurs
precisely at branchpoints of the tree. Fragmentations of CRT's that are
different from the Brownian one and that follow the same kind of
construction have been studied in \cite{mierfmoins}.

The goal of this paper is to show that any self-similar fragmentation
process with negative index can be obtained by a similar construction as
above, for a certain instance of CRT. We are interested in negative indices,
because in most interesting cases when the self-similarity index is
non-negative, all fragments have an ``infinite lifetime'', meaning that the
pieces of the fragmentation remain macroscopic at all times. In this case,
the family tree defined above will be unbounded and without endpoints, hence
looking completely different from the Brownian CRT. By contrast, as soon as
the self-similarity index is negative, a loss of mass occurs, that makes the
fragments disappear in finite time (see \cite{bertafrag02}). In this case,
the metric family tree will be a bounded object, and in fact, a CRT. To
state our results, we first give a rigorous definition of the involved
objects. Call
\begin{equation*}
S=\left\{ \mathbf{s}=(s_{1},s_{2},\ldots ):s_{1}\geq s_{2}\geq \ldots \geq
0;\sum_{i\geq 1}s_{i}\leq 1\right\} ,
\end{equation*}
and endow it with the topology of pointwise convergence.

\begin{defn}
A Markovian $S$-valued process $(F(t),t\geq 0)$ starting at
$(1,0,\ldots)$ is a ranked self-similar
fragmentation with index $\alpha \in \mathbb{R}$ if it is continuous in
probability and satisfies the following fragmentation property. For every $%
t,t^{\prime }\geq 0$, given $F(t)=(x_{1},x_{2},\ldots )$, $F(t+t^{\prime })$
has the same law as the decreasing rearrangement of the sequences $%
x_{1}F^{(1)}(x_{1}^{\alpha }t^{\prime }),x_{2}F^{(2)}(x_{2}^{\alpha
}t^{\prime }),\ldots $, where the $F^{(i)}$'s are independent copies of $F$.
\end{defn}

By a result of Bertoin \cite{bertsfrag02} and Berestycki \cite{berest02},
the laws of such fragmentation processes are characterized by a $3$-tuple $%
(\alpha ,c,\nu )$, where $\alpha $ is the index, $c\geq 0$ is an ``erosion''
constant, and $\nu $ is a $\sigma $-finite measure on $S$ that integrates $%
\mathbf{s}\mapsto 1-s_{1}$ such that $\nu (\{(1,0,0\ldots )\})=0$.
Informally, $c$ measures the rate at which fragments melt continuously (a
phenomenon we will not be much interested in here), while $\nu $ measures
instantaneous breaks of fragments: a piece with size $x$ breaks into
fragments with masses $x\mathbf{s}$ at rate $x^{\alpha }\nu (\mathrm{d}%
\mathbf{s})$. Notice that some mass can be lost within a sudden break: this
happens as soon as $\nu (\sum_{i}s_{i}<1)\neq 0$, but we will not be
interested in this phenomenon here either. The loss of mass phenomenon
stated above is completely different from erosion or sudden loss of mass: it
is due to the fact that small fragments tend to decay faster when $\alpha <0$%
.

On the other hand, let us define the notion of CRT. An $\mathbb{R}$-tree
(with the terminology of Dress and Terhalle \cite{drter}; it is called \emph{%
continuum tree set} in Aldous \cite{aldouscrt93}) is a complete metric space
$(T,d)$, whose elements are called \emph{vertices}, which satisfies the
following two properties:

\begin{itemize}
\item  For $v,w\in T$, there exists a unique geodesic $[[v,w]]$ going from $%
v $ to $w$, i.e.\ there exists a unique isomorphism $\varphi
_{v,w}:[0,d(v,w)]\rightarrow T$ with $\varphi _{v,w}(0)=v$ and $\varphi
_{v,w}(d(v,w))=w$, and its image is called $[[v,w]]$.

\item  For any $v,w\in T$, the only non-self-intersecting path going from $v$
to $w$ is $[[v,w]]$, i.e.\ for any continuous injective function $s\mapsto
v_{s}$ from $[0,1]$ to $T$ with $v_{0}=v$ and $v_{1}=w$, $\{v_{s}:s\in
\lbrack 0,1]\}=[[v,w]]$.
\end{itemize}

We will furthermore consider $\mathbb{R}$-trees that are \emph{rooted}, that
is, one vertex is distinguished as being the root, and we call it 
$\varnothing $. A \emph{leaf} is a vertex which does not belong to $
[[\varnothing ,w[[:=\varphi _{\varnothing ,w}([0,d(\varnothing ,w)))$ for
any vertex $w$. Call $\mathcal{L} (T)$ the set of leaves of $T$, and $
\mathcal{S}(T)=T\setminus \mathcal{L} (T)$ its skeleton. An $\mathbb{R}$
-tree is \emph{leaf-dense} if $T$ is the closure of $\mathcal{L} (T)$. We
also call \emph{height} of a vertex $v$ the quantity $\mathrm{ht}
(v)=d(\varnothing ,v)$. Last, for $T $ an $\mathbb{R} $-tree and $a>0$, we
let $a\otimes T$ be the $\mathbb{R}$-tree in which all distances are
multiplied by $a$.

\begin{defn}
A continuum tree is a pair $(T,\mu )$ where $T$ is an $\mathbb{R}$-tree and 
$\mu $ is a probability measure on $T$, called the mass measure, which is
non-atomic and satisfies $\mu (\mathcal{L}(T))=1$ and such that for every
non-leaf vertex $w$, $\mu \{v\in T:[[\varnothing ,v]]\cap \lbrack \lbrack
\varnothing ,w]]=[[\varnothing ,w]]\}>0$. The set of vertices just defined
is called the fringe subtree rooted at $w$. A CRT is a random variable $
\omega \mapsto (T(\omega ),\mu (\omega ))$ on a probability space $(\Omega ,
\mathcal{F},P)$ which values are continuum trees.
\end{defn}

Notice that the definition of a continuum tree implies that the 
$\mathbb{R}$-tree $T$ satisfies certain extra properties, for example, its set 
of leaves must be uncountable and have no isolated point. Also, the definition
of a CRT is a little inaccurate as we did not endow the space of $\R$-trees
with a $\sigma$-field. This problem is in fact circumvented by the fact that
CRTs are in fact entirely described by the sequence of their {\em marginals}, 
that is, of the subtrees spanned by the root and $k$ leaves chosen with law
$\mu$ given $\mu$, and these subtrees, which are interpreted as finite 
{\em trees with edge-lengths}, are random variables (see Sect.\ \ref{twel}).  
The reader should keep in mind that by the ``law'' of a CRT we mean the
sequence of these marginals. 

For $(T,\mu )$ a continuum tree, and for every $t\geq 0$, let 
$T_{1}(t),T_{2}(t),\ldots $ be the tree components of $\{v\in T:\mathrm{ht}
(v)>t\}$, ranked by decreasing order of $\mu $-mass. A continuum random tree
$(T,\mu )$ is said to be \textit{self-similar} with index $\alpha <0$ if for
every $t\geq 0$, conditionally on $(\mu (T_{i}(t)),i\geq 1)$, $
(T_{i}(t),i\geq 0)$ has the same law as $(\mu (T_{i}(t))^{-\alpha }\otimes
T^{(i)},i\geq 1)$ where the $T^{(i)}$'s are independent copies of $T$. 

Our first result is

\begin{thm}
\label{T1} Let $F$ be a ranked self-similar fragmentation process with
characteristic $3$-tuple $(\alpha ,c,\nu )$, with $\alpha <0$. Suppose also
that $F$ is not constant, that $c=0$ and $\nu (\sum_{i}s_{i}<1)=0$. Then
there exists an $\alpha $-self-similar CRT $(\mathcal{T}_{F},\mu _{F})$ such
that, writing $F^{\prime }(t)$ for the decreasing sequence of masses of
connected components of the open set $\{v\in \mathcal{T}_{F}:\mathrm{ht}%
(v)>t\}$, the process $(F^{\prime }(t),t\geq 0)$ has the same law as $F$.
The tree $\mathcal{T}_{F}$ is leaf-dense if and only if $\nu $ has infinite
total mass.
\end{thm}

The next statement is a kind of converse to this theorem.

\begin{prp}
\label{reciproque} Let $(\mathcal{T},\mu )$ be a self-similar CRT with index
$\alpha <0$. Then the process $F(t)=((\mu (\mathcal{T}_{i}(t),i\geq 1),t\geq
0)$ is a ranked self-similar fragmentation with index $\alpha $, it has no
erosion and its dislocation measure $\nu $ satisfies $\nu
(\sum_{i}s_{i}<1)=0 $. Moreover, $\mathcal{T}_{F}$ and $\mathcal{T}$ have
the same law.
\end{prp}

These results are proved in Sect.\ \ref{crtf}. There probably exists some
notion of continuum random tree extending the former which would include
fragmentations with erosion or with sudden loss of mass, but such
fragmentations usually are less interesting.

The next result, to be proved in Sect.\ \ref{hausdorff}, deals with the
Hausdorff dimension of the set of leaves of the CRT $\mathcal{T }_{F}$.

\begin{thm}
\label{T2} Let $F$ be a ranked self-similar fragmentation with
characteristics $(\alpha ,c,\nu )$ satisfying the hypotheses of Theorem \ref
{T1}. Writing $\dim_{\cal H}$ for Hausdorff dimension, one
has
\begin{equation}
\dim_{\cal H}\left(\mathcal{L}(\mathcal{T}_{F})\right) =\frac{%
1}{\left| \alpha \right| }\text{ a.s}.
\end{equation}
as soon as $\int_{S}\left( s_{1}^{-1}-1\right) \nu ({\mathrm{d}\mathbf{s}}%
)<\infty $
\end{thm}

Some comments about this formula. First, notice that under the extra
integrability assumption on $\nu $, the dimension of the whole tree is $%
\dim_{\cal H}(\mathcal{T}_{F})=(1/|\alpha |)\vee 1$ because
the skeleton $\mathcal{S}(\mathcal{T}_{F})$ has dimension $1$ as a countable
union of segments. The value $-1$ is therefore critical for $\alpha $, since
the above formula shows that the dimension of $\mathcal{T}_{F}$ as to be $1$
as soon as $\alpha \leq -1$. It was shown in a previous work by Bertoin \cite
{bertafrag02} that when $\alpha <-1$, for every fixed $t$ the number of
fragments at time $t$ is a.s.\ finite, so that $-1$ is indeed the threshold
under which fragments decay extremely fast. One should then picture the CRT $%
\mathcal{T}_{F}$ as a ``dead tree'' looking like a handful of thin sticks
connected to each other, while when $|\alpha |<1$ the tree looks more like a
dense ``bush''. Last, the integrability assumption in the theorem seems to
be reasonably mild; its heuristic meaning is that when a fragmentation
occurs, the largest resulting fragment is not too small. In particular, it
is always satisfied in the case of fragmentations for which $\nu
(s_{N+1}>0)=0$, since then $s_{1}>1/N$ for $\nu $-a.e.\ $\mathbf{s}$. Yet,
we point out that when $\int_{S}\left( s_{1}^{-1}-1\right) \nu (\mathrm{d}%
\mathbf{s})=\infty $, one anyway obtains the following bounds for the
Hausdorff dimension of $\mathcal{L}(\mathcal{T}_{F})$:
\begin{equation*}
\frac{\varrho }{\left| \alpha \right| }\leq \dim_{\cal H}
\left( \mathcal{L}(\mathcal{T}_{F})\right) \leq \frac{1}{\left| \alpha
\right| }\text{ a.s}.
\end{equation*}
where
\begin{equation}
\varrho :=\sup \left\{ p\leq 1:\int_{S}\left( s_{1}^{-p}-1\right) \nu ({%
\mathrm{d}\mathbf{s}})<\infty \right\} .  \label{219}
\end{equation}

It is worth noting that these results allow as a special case to compute the
Hausdorff dimension of the so-called \emph{stable trees} of Duquesne and Le
Gall \cite{duqleg02}, which were used to construct fragmentations in the
manner of Theorem \ref{T1} in \cite{mierfmoins}. The dimension of the stable
tree (as well as finer results of Hausdorff measures on more general L\'evy
trees) has been obtained independently in \cite{duqlegprep}. The stable tree
is a CRT whose law depends on parameter $\beta\in(1,2]$, and it satisfies
the required self-similarity property of Proposition \ref{reciproque} with
index $1/\beta-1$. We check that the associated dislocation measure
satisfies the integrability condition of Theorem \ref{T2} in Sect.\ \ref
{checkst}, so that

\begin{crl}
\label{stable} Fix $\beta \in (1,2]$. The $\beta $-stable tree has Hausdorff
dimension $\beta /(\beta -1)$.
\end{crl}

An interesting process associated to a given continuum tree $(T,\mu)$ is the
so-called \emph{cumulative height profile} $\bar{W}_{T}(h)=\mu\{v\in T:%
\mathrm{ht}(v)\leq h\}$, which is non-decreasing and bounded by $1$ on 
$\mathbb{R}_+$.  It may happen that the Stieltjes measure
$\mathrm{d}\bar{W}_{T}(h)$ 
is absolutely continuous with respect to Lebesgue measure, in which
case its density $(W_{T}(h),h\geq 0)$ is called the \emph{height profile},
or \emph{width process} of the tree. In our setting, for any fragmentation $%
F $ satisfying the hypotheses of Theorem \ref{T1}, the cumulative height
profile has the following interpretation: one has $(\bar{W}_{\mathcal{T}%
_{F}}(h),h\geq 0)$ has the same law as $(M_F(h),h\geq 0)$, where $%
M_F(h)=1-\sum_{i\geq 1}F_i(h)$ is the total mass lost by the fragmentation
at time $h$. Detailed conditions for existence (or non-existence) of the
width profile $\mathrm{d} M_F(h)/\mathrm{d} h$ have been given in \cite
{haas03}. It was also proved there that under some mild assumptions $\mathrm{%
dim\,} _{\mathcal{H}}\left( \mathrm{d} M_{F}\right) \geq 1\wedge A/\left|
\alpha \right| $ a.s., where $A$ is a $\nu $-dependent parameter introduced
in $\left( \ref{278}\right) $, Section 3 below. The upper bound we obtain
for $\dim_{\cal H}\left( \mathcal{L}(\mathcal{T}_{F})\right)
$ allows us to complete this result:

\begin{crl}
\label{coromass} Let $F$ be a ranked self-similar fragmentation with same
hypotheses\ as in Theorem \ref{T1}. Then $\dim_{\cal H}\left(
\mathrm{d}M_{F}\right) \leq 1\wedge 1/\left| \alpha \right| $ a.s.
\end{crl}

Notice that this result re-implies the fact from \cite{haas03} that the
height profile does not exist as soon as $|\alpha |\geq 1$.

The last motivation of this paper (Sect.\ \ref{heightprocess}) is about
relations between CRTs and their so-called encoding \emph{height processes}.
The fragmentation $F_B$ of \cite{bertsfrag02}, as well as the fragmentations
from \cite{mierfmoins}, were defined out of certain random functions $%
(H_s,0\leq s\leq 1)$. Let us describe briefly the construction of $F_B$. Let
$B^{\mathrm{exc}}$ be the standard Brownian excursion with duration $1$, and
consider the open set $\{s\in [0,1]:2B^{\mathrm{exc}}_s>t\}$. Write $F(t)$
for the decreasing sequence of the lengths of its interval components. Then $%
F$ has the same law as the fragmentation $F_B^{\prime}$ defined out of the
Brownian CRT in the same way as in Theorem \ref{T1}. This is immediate from
the description of Le Gall \cite{legall93} and Aldous \cite{aldouscrt93} of
the Brownian tree as being encoded in the Brownian excursion. To be concise,
define a pseudo-metric on $[0,1]$ by letting $\overline{d}(s,s^{\prime})=
2B^{\mathrm{exc}}_s+2B^{\mathrm{exc}}_{s^{\prime}}-4\inf_{u\in
[s,s^{\prime}]}B^{\mathrm{exc}}_u$, with the convention that $%
[s,s^{\prime}]=[s^{\prime},s]$ if $s^{\prime}<s$. We can define a true
metric space by taking the quotient with respect to the equivalence relation
$s\equiv s^{\prime}\iff \overline{d}(s,s^{\prime})=0$. Call $( \mathcal{T }%
_B,d)$ this metric space. Write $\mu_B$ for the measure induced on $\mathcal{%
T }_B$ by Lebesgue measure on $[0,1]$. Then $( \mathcal{T }_B,\mu_B) $ is
the Brownian CRT, and the equality in law of the fragmentations $F_B$ and $%
F_B^{\prime}$ follows immediately from the definition of the mass measure.
Our next result generalizes this construction.

\begin{thm}
\label{T3} Let $F$ be a ranked self-similar fragmentation with same
hypotheses as in Theorem \ref{T1}, and suppose $\nu $ has infinite total
mass. Then there exists a continuous random function $(H_{F}(s),0\leq s\leq
1)$, called the height function, such that $H_{F}(0)=H_{F}(1)$, $H_{F}(s)>0$
for every $s\in (0,1)$, and such that $F$ has the same law as the
fragmentation $F^{\prime }$ defined by: $F^{\prime }(t)$ is the decreasing
rearrangement of the lengths of the interval components of the open set $%
I_{F}(t)=\{s\in (0,1):H_{F}(s)>t\}$.
\end{thm}

An interesting point in this construction is also that it shows that a large
class of self-similar fragmentation with negative index has a natural \emph{%
interval representation}, given by $(I_F(t),t\geq 0)$.

In parallel to the computation of the Hausdorff dimension of the CRTs built
above, we are able to estimate H\"older coefficients for the height
processes of these CRTs. Our result is

\begin{thm}
\label{T4} Suppose $\nu (S)=\infty $, and set
\begin{eqnarray*}
\vartheta _{\mathrm{low}}:= &&\sup \left\{ b>0:\lim_{x\downarrow 0}x^{b}\nu
(s_{1}<1-x)=\infty \right\} , \\
\vartheta _{\mathrm{up}}:= &&\inf \left\{ b>0:\lim_{x\downarrow 0}x^{b}\nu
(s_{1}<1-x)=0\right\} .
\end{eqnarray*}
Then the height process $H_{F}$ is a.s.\ H\"{o}lder-continuous of order $%
\gamma $ for every $\gamma <\vartheta _{\mathrm{low}}\wedge |\alpha |$, and,
provided that $\int_{S}(s_{1}^{-1}-1)\nu (\mathrm{d}\mathbf{s})<\infty ,$
a.s.\ not H\"{o}lder-continuous of order $\gamma $ for every $\gamma
>\vartheta _{\mathrm{up}}\wedge |\alpha |$.
\end{thm}

Again we point out that one actually obtains an upper bound for the maximal
H\"{o}lder coefficient even when $\int_{S}(s_{1}^{-1}-1)\nu (\mathrm{d}%
\mathbf{s})=\infty :$ with $\varrho $ defined by $\left( \ref{219}\right) ,$
a.s. $H_{F}$ cannot be H\"{o}lder-continuous of order $\gamma $ for any $%
\gamma >\vartheta _{\mathrm{up}}\wedge |\alpha |/\varrho .$

Note that $\vartheta _{\mathrm{low}},\vartheta _{\mathrm{up}}$ depend only
on the characteristics of the fragmentation process, and more precisely, on
the behavior of $\nu $ when $s_{1}$ is close to $1$. By contrast, our
Hausdorff dimension result for the tree depended on a hypothesis on the
behavior of $\nu $ when $s_{1}$ is near $0$. Remark also that $\vartheta _{%
\mathrm{up}}$ may be strictly smaller than $1$. Therefore, the Hausdorff
dimension of $\mathcal{T}_{F}$ is in general not equal to the inverse of the
maximal H\"{o}lder coefficient of the height process, as one could have
expected. However, this turns out to be true in the case of the stable tree,
as will be checked in Section \ref{heightstable}:

\begin{crl}
\label{C2} The height process of the stable tree with index $\beta \in (1,2]$
is a.s.\ H\"{o}lder-continuous of any order $\gamma <1-1/\beta $, but a.s.\
not of order $\gamma >1-1/\beta $.
\end{crl}

When $\beta =2$, this just states that the Brownian excursion is
H\"{o}lder-continuous of any order $<1/2$, a result that is well-known for
Brownian motion and which readily transfers to the normalized Brownian
excursion (e.g.\ by rescaling the first excursion of Brownian motion whose
duration is greater than $1$). The general result had been obtained in \cite
{duqleg02} by completely different methods.

Last, we mention that most of our results extend to a more general class of
fragmentations in which a fragment with mass $x$ splits to give fragments
with masses $x$\textbf{s}$,$ $\mathbf{s}\in S,$ at rate $\varsigma (x)\nu ($d%
$\mathbf{s})$ for some non-negative continuous function $\varsigma $ on $%
\left( 0,1\right] $ (see \cite{haas02} for a rigorous definition). The
proofs of the above theorems easily adapt to give the following results:
when $\liminf_{x\rightarrow 0}x^{-b}\varsigma (x)>0$ for some $b<0,$ the
fragmentation can be encoded as above into a CRT and, provided that $\nu $
is infinite, into a height function. The set of leaves of the CRT then has a
Hausdorff dimension \ smaller than $1/\left| b\right| $ and the height
function is $\gamma $-H\"{o}lder continuous for every $\gamma <\vartheta _{%
\text{low}}\wedge \left| b\right| .$ If moreover $\limsup_{x\rightarrow
0}x^{-a}\varsigma (x)<\infty $ for some $a<0$ and $\int_{S}\left(
s_{1}^{-1}-1\right) \nu (\mathrm{d}\mathbf{s})<\infty ,$ the Hausdorff
dimension is larger than $1/\left| a\right| $ and the height function cannot
have a H\"{o}lder coefficient $\gamma >\vartheta _{\text{sup}}\wedge \left|
a\right| $.

\section{The CRT $\mathcal{T }_F$}

\label{crtf}

Building the CRT $\mathcal{T }_F$ associated to a ranked fragmentation $F$
will be done by determining its ``marginals'', i.e.\ the subtrees spanned by
a finite but arbitrary number of randomly chosen leaves. To this purpose, it
will be useful to use partitions-valued fragmentations, which we first
define, as well as a certain family of trees with edge-lengths.

\subsection{Exchangeable partitions and partitions-valued self-similar
fragmentations}

\label{pvssf}

Let $\mathcal{P}_{\infty }$ be the set of (unordered) partitions of $\mathbb{%
N}=\{1,2,\ldots \}$ and $[n]=\{1,2,\ldots ,n\}$. We adopt the following
ordering convention: for $\pi \in \mathcal{P}_{\infty }$, we let $(\pi
_{1},\pi _{2},\ldots )$ be the blocks of $\pi $, so that $\pi _{i}$ is the
block containing $i$ provided that $i$ is the smallest integer of the block
and $\pi _{i}=\varnothing $ otherwise. We let $\mathbb{O}=\{\{1\},\{2\},%
\ldots \}$ be the partition of $\mathbb{N}$ into singletons. If $B\subset
\mathbb{N}$ and $\pi \in \mathcal{P}_{\infty }$ we let $\pi \cap B$ (or $\pi
|_{B}$) be the restriction of $\pi $ to $B$, i.e.\ the partition of $B$
whose collection of blocks is $\{\pi _{i}\cap B,i\geq 1\}$. If $\pi \in
\mathcal{P}_{\infty }$ and $B\in \pi $ is a block of $\pi $, we let
\begin{equation*}
|B|=\lim_{n\rightarrow \infty }\frac{\#(B\cap \lbrack n])}{n}
\end{equation*}
be the asymptotic frequency of the block $B$, whenever it exists. A random
variable $\pi $ with values in $\mathcal{P}_{\infty }$ is called \emph{%
exchangeable} if its law is invariant by the natural action of permutations
of $\mathbb{N}$ on $\mathcal{P}_{\infty }$. By a theorem of Kingman \cite
{kingman78,aldous85}, all the blocks of such random partitions admit
asymptotic frequencies a.s. For $\pi $ whose blocks have asymptotic
frequencies, we let $|\pi |\in S$ be the decreasing sequence of these
frequencies. Kingman's theorem more precisely says that the law of any
exchangeable random partition $\pi $ is a (random) ``paintbox process'', a
term we now explain. Take $\mathbf{s}\in S$ (the paintbox) and consider a
sequence $U_{1},U_{2},\ldots $ of i.i.d.\ variables in $\mathbb{N}\cup \{0\}$
(the colors) with $P(U_{1}=j)=s_{j}$ for $j\geq 1$ and $P(U_{1}=0)=1-%
\sum_{k}s_{k}$. Define a partition $\pi $ on $\mathbb{N}$ by saying that $%
i\neq j$ are in the same block if and only if $U_{i}=U_{j}\neq 0$ (i.e.\ $i$
and $j$ have the same color, where $0$ is considered as colorless). Call $%
\rho _{\mathbf{s}}(\mathrm{d}\pi )$ its law, the \emph{$\mathbf{s}$-paintbox}
law. Kingman's theorem says that the law of any random partition is a mixing
of paintboxes, i.e.\ it has the form $\int_{\mathbf{s}\in S}m(\mathrm{d}{%
\mathbf{s}})\rho _{\mathbf{s}}(\mathrm{d}\pi )$ for some probability measure
$m$ on $S$. A useful consequence is that the block of an exchangeable
partition $\pi $ containing $1$, or some prescribed integer $i$, is a \emph{%
size-biased pick} from the blocks of $\pi $, i.e.\ the probability it equals
a non-singleton block $\pi _{j}$ conditionally on $(|\pi _{j}|,j\geq 1)$
equals $|\pi _{j}|$. Similarly,

\begin{lmm}
\label{usefull} Let $\pi $ be an exchangeable random partition which is
a.s.\ different from the trivial partition $\mathbb{O}$, and $B$ an infinite
subset of $\mathbb{N}$. For any $i\in \mathbb{N}$, let
\begin{equation*}
\widetilde{i}=\inf \{j\geq i:j\in B\mbox{ and }\{j\}\notin \pi \},
\end{equation*}
then $\widetilde{i}<\infty $ a.s.\ and the block $\widetilde{\pi}$ of $\pi $
containing $\widetilde{i}$ is a size-biased pick among the non-singleton
blocks of $\pi $, i.e.\ if we denote these by $\pi _{1}^{\prime },\pi
_{2}^{\prime },\ldots $,
\begin{equation*}
P(\widetilde{\pi}=\pi _{k}^{\prime }|(|\pi _{j}^{\prime }|,j\geq 1))=|\pi
_{k}^{\prime }|/\sum_{j}|\pi _{j}^{\prime }|.
\end{equation*}
\end{lmm}

For any sequence of partitions $(\pi ^{(i)},i\geq 1),$ define $%
\pi=\bigcap_{i\geq 1}\pi ^{(i)}$ by
\begin{equation*}
k\overset{\pi}{\sim }j\iff k\overset{\pi ^{(i)}}{\sim }j\quad \forall i\geq
1.
\end{equation*}

\begin{lmm}
\label{Lemmafreqasym} Let $(\pi^{(i)},i\geq 1)$ be a sequence of independent
exchangeable partitions and set $\pi :=\bigcap_{i\geq 1}\pi ^{(i)}.$ Then,
a.s.\ for every $j\in \mathbb{N},$%
\begin{equation*}
\left| \pi _{j}\right| =\prod_{i\geq 1}\left| \pi _{k(i,j)}^{(i)}\right| ,
\end{equation*}
where $(k(i,j),j\geq 1)$ is defined so that $\pi_{j}=\bigcap_{i\geq 1}\pi
_{k(i,j)}^{(i)}.$
\end{lmm}

{\noindent \textbf{Proof. }} First notice that $k(i,j)\leq j$ for all $i\geq
1$ a.s. This is clear when $\pi _{j}\neq \varnothing $, since $j\in \pi _{j}$
and then $j\in \pi _{k(i,j)}^{(i)}.$ When $\pi _{j}=\varnothing ,$ $j\in \pi
_{m}$ for some $m<j $ and then $m$ and $j$ belong to the same block of $\pi
^{(i)}$ for all $i\geq 1.$ Thus $k(i,j)\leq m<j.$ Using then the paintbox
construction of exchangeable partitions explained above and the independence
of the $\pi ^{(i)}$'s$,\ $we see that the r.v. $\prod_{i\geq 1}\mathbf{1}%
_{\{m\in \pi _{k(i,j)}^{(i)}\}},$ $m\geq j+1,$ are iid conditionally on $%
(|\pi _{k(i,j)}^{(i)}|,i\geq 1)$ with a mean equal to $\prod_{i\geq 1}|\pi
_{k(i,j)}^{(i)}|.$ The law of large numbers therefore gives
\begin{equation*}
\prod_{i\geq 1}\left| \pi _{k(i,j)}^{(i)}\right| =\lim_{n\rightarrow \infty }%
\frac{1}{n}\sum_{j+1\leq m\leq n}\prod_{i\geq 1}\mathbf{1}_{\left\{ m\in \pi
_{k(i,j)}^{(i)}\right\} }\text{ a.s}.
\end{equation*}
On the other hand, the random variables $\prod_{i\geq 1}\mathbf{1}_{\{ m\in
\pi _{k(i,j)}^{(i)}\} }=\mathbf{1}_{\left\{ m\in \pi _{j}\right\} }, $ $%
m\geq j+1,$ are i.i.d.\ conditionally on $\left| \pi _{j}\right| $ with mean
$\left| \pi _{j}\right| $ and then the limit above converges a.s. to $\left|
\pi _{j}\right| ,$ again by the law of large numbers. $\hfill \square$

We now turn our attention to partitions-valued fragmentations.

\begin{defn}
Let $(\Pi (t),t\geq 0)$ be a Markovian $\mathcal{P}_{\infty }$-valued
process with $\Pi (0)=\{\mathbb{N},\varnothing ,\varnothing ,\ldots \}$ that
is continuous in probability and exchangeable as a process (meaning that the
law of $\Pi $ is invariant by the action of permutations). Call it a
partition-valued self-similar fragmentation with index $\alpha \in \mathbb{R}
$ if moreover $\Pi (t)$ admits asymptotic frequencies for all $t$, a.s., if
the process $(|\Pi (t)|,t\geq 0)$ is continuous in probability, and if the
following fragmentation property is satisfied. For $t,t^{\prime }\geq 0$,
given $\Pi (t)=(\pi _{1},\pi _{2},\ldots )$, the sequence $\Pi (t+t^{\prime
})$ has the same law as the partition with blocks $\pi _{1}\cap \Pi
^{(1)}(|\pi _{1}|^{\alpha }t^{\prime }),\pi _{2}\cap \Pi ^{(2)}(|\pi
_{2}|^{\alpha }t^{\prime }),\ldots $, where $(\Pi ^{(i)},i\geq 1)$ are
independent copies of $\Pi $.
\end{defn}

Bertoin \cite{bertsfrag02} has shown that any such fragmentation is also
characterized by the same $3$-tuple $(\alpha ,c,\nu )$ as above, meaning
that the laws of partition-valued and ranked self-similar fragmentations are
in a one-to-one correspondence. In fact, for every $(\alpha ,c,\nu )$, one
can construct a version of the partition-valued fragmentation $\Pi $ with
parameters $(\alpha ,c,\nu )$, and then $(|\Pi (t)|,t\geq 0)$ is the ranked
fragmentation with parameters $(\alpha ,c,\nu )$. Let us build this version
now. It is done following \cite{berthfrag01,bertsfrag02} by a Poissonian
construction. Recall the notation $\rho _{\mathbf{s}}(\mathrm{d}\pi )$, and
define $\kappa _{\nu }(\mathrm{d}\pi )=\int_{S}\nu (\mathrm{d}{\mathbf{s}}%
)\rho _{\mathbf{s}}(\mathrm{d}\pi )$. Let $\#$ be the counting measure on $%
\mathbb{N}$ and let $(\Delta _{t},k_{t})$ be a $\mathcal{P}_{\infty }\times
\mathbb{N}$-valued Poisson point process with intensity $\kappa _{\nu
}\otimes \#$. We may construct a process $(\Pi ^{0}(t),t\geq 0)$ by letting $%
\Pi ^{0}(0)$ be the trivial partition $(\mathbb{N},\varnothing ,\varnothing
,\ldots )$, and saying that $\Pi ^{0} $ jumps only at times $t$ when an atom
$(\Delta _{t},k_{t})$ occurs. When this is the case, $\Pi ^{0}$ jumps from
the state $\Pi ^{0}(t-)$ to the following partition $\Pi ^{0}(t)$: replace
the block $\Pi _{k_{t}}^{0}(t-)$ by $\Pi _{k_{t}}^{0}(t-)\cap \Delta _{t}$,
and leave the other blocks unchanged. Such a construction can be made
rigorous by considering restrictions of partitions to the first $n$ integers
and by a consistency argument. Then $\Pi ^{0}$ has the law of the
fragmentation with parameters $(0,0,\nu )$.

Out of this ``homogeneous'' fragmentation, we construct the $(\alpha ,0,\nu
) $-fragmentation by introducing a time-change. Call $\lambda _{i}(t)$ the
asymptotic frequency of the block of $\Pi ^{0}(t)$ that contains $i$, and
write
\begin{equation}
T_{i}(t)=\inf \left\{ u\geq 0:\int_{0}^{u}\lambda _{i}(r)^{-\alpha }\mathrm{d%
} r>t\right\} .  \label{1}
\end{equation}
Last, for every $t\geq 0$ we let $\Pi (t)$ be the random partition such that
$i,j$ are in the same block of $\Pi (t)$ if and only if they are in the same
block of $\Pi ^{0}(T_{i}(t))$, or equivalently of $\Pi ^{0}(T_{j}(t))$. Then
$(\Pi (t),t\geq 0)$ is the wanted version. Let $\left( \mathcal{G}(t),t\geq
0\right) $ be the natural filtration generated by $\Pi $ completed up to $P$%
-null sets. According to \cite{bertsfrag02}, the fragmentation property
holds actually for $\mathcal{G}$-stopping times and we shall refer to it as
the \textit{strong fragmentation property.} In the homogeneous case, we will
rather call $\mathcal{G}^{0}$ the natural filtration.

When $\alpha<0$, the loss of mass in the ranked fragmentations shows up at
the level of partitions by the fact that a positive fraction of the blocks
of $\Pi(t)$ are singletons for some $t>0$. This last property of
self-similar fragmentations with negative index allows to build a collection
of trees with edge-lengths.

\subsection{Trees with edge-lengths}\label{twel}

A tree is a finite connected graph with no cycles. It is \emph{rooted} when
a particular vertex (the root) is distinguished from the others, in this
case the edges are by convention oriented, pointing from the root, and we
define the out-degree of a vertex $v$ as being the number of edges that
point outward $v$. A \emph{leaf} in a rooted tree is a vertex with
out-degree $0$. For $k\geq 1$, let $\mathbf{T}_k$ be the set of rooted trees
with exactly $k$ labeled leaves (the names of the labels may change
according to what we see fit), the other vertices (except the root) begin
unlabeled , and such that the root is the only vertex that has out-degree $1$%
. If $\mathbf{t}\in \mathbf{T}_k$, we let $E(\mathbf{t})$ be the set of its
edges.

A \emph{tree with edge-lengths} is a pair $\vartheta =(\mathbf{t},\mathbf{e}%
) $ for $\mathbf{t}\in \bigcup_{k\geq 1}\mathbf{T}_{k}$ and $\mathbf{e}%
=(e_{i},i\in E(t))\in (\mathbb{R}_{+}\setminus \{0\})^{E(\mathbf{t})}$. Call
$\mathbf{t}$ the \emph{skeleton} of $\vartheta $. Such a tree is naturally
equipped with a distance $d(v,w)$ on the set of its vertices, by adding the
lengths of edges that appear in the unique path connecting $v$ and $w$ in
the skeleton (which we still denote by $[[v,w]]$). The height of a vertex is
its distance to the root. We let $\mathbb{T}_{k}$ be the set of trees with
edge-lengths whose skeleton is in $\mathbf{T}_{k}$. For $\vartheta \in
\mathbb{T}_{k}$, let $e_{\mathrm{root}}$ be the length of the unique edge
connected to the root, and for $e<e_{\mathrm{root}}$ write $\vartheta -e $
for the tree with edge-lengths that has same skeleton and same edge-lengths
as $\vartheta $, but for the edge pointing outward the root which is
assigned length $e_{\mathrm{root}}-e$.

We also define an operation $\mathtt{MERGE}$ as follows. Let $n\geq2$ and
take $\vartheta_1,\vartheta_2,\ldots,\vartheta_n$ respectively in $\mathbb{T}%
_{k_1},\mathbb{T}_{k_2},\ldots,\mathbb{T}_{k_n}$, with leaves $(L^1_i,1\leq
i\leq k_1),(L^2_i,1\leq i\leq k_2),\ldots,(L^n_i,1\leq i\leq k_n)$
respectively. Let also $e>0$. The tree with edge-lengths $\mathtt{MERGE}%
((\vartheta_1,\ldots,\vartheta_n);e)\in \mathbb{T}_{\sum_i k_i}$ is defined
by merging together the roots of $\vartheta_1,\ldots,\vartheta_n$ into a
single vertex $\bullet$, and by drawing a new edge $\mathrm{root}\to\bullet$
with length $e$.

Last, for $\vartheta\in \mathbb{T}_k$ and $i$ vertices $v_1,\ldots,v_i$,
define the subtree spanned by the root and $v_1,\ldots,v_i$ as follows. For
every $p\neq q$, let $b(v_p,v_q)$ be the branchpoint of $v_p$ and $v_q$,
that is, the highest point in the tree that belongs to $[[\mathrm{root}%
,v_p]]\cap [[\mathrm{root},v_q]]$. The spanned tree is the tree with
edge-lengths whose vertices are the root, the vertices $v_1,\ldots,v_i$ and
the branchpoints $b(v_p,v_q)$, $1\leq p\neq q\leq i$, and whose edge-lengths
are given by the respective distances between this subset of vertices of the
original tree.

\subsection{Building the CRT}

Now for $B\subset \mathbb{N}$ finite, define $\mathcal{R}(B)$, a random
variable with values in $\mathbb{T}_{\#B}$, whose leaf-labels are of the
form $L_i$ for $i\in \mathbb{N}$ , as follows. Let $D_i=\inf\{t\geq
0:\{i\}\in\Pi(t)\}$ be the first time when $\{i\}$ ``disappears'', i.e.\ is
isolated in a singleton of $\Pi(t)$. For $B$ a finite subset of $\mathbb{N}$
with at least two elements, let $D_B=\inf\{t\geq 0:\#(B\cap \Pi(t))\neq 1\}$
be the first time when the restriction of $\Pi(t)$ to $B$ is non-trivial,
i.e.\ has more than one block. By convention, $D_{\{i\}}=D_i$. For every $%
i\geq 1$, define $\mathcal{R}(\{i\})$ as a single edge $\mathrm{root}\to L_i$%
, and assign this edge the length $D_i$. For $B$ with $\#B\geq 2$, let $%
B_1,\ldots,B_i$ be the non-empty blocks of $B\cap\Pi(D_B)$, arranged in
increasing order of least element, and define a tree $\mathcal{R}(B)$
recursively by
\begin{equation*}
\mathcal{R}(B)=\mathtt{MERGE}((\mathcal{R}(B_1)-D_B,\ldots,\mathcal{R}%
(B_i)-D_B);D_B).
\end{equation*}
Last, define $\mathcal{R}(k)=\mathcal{R}([k])$. Notice that by definition of
the distance, the distance between $L_i$ and $L_j$ in $\mathcal{R}(k)$ for
any $k\geq i\vee j$ equals $D_i+D_j-2D_{\{i,j\}}$.

We now state the key lemma that allows to describe the CRT out of the family
$(\mathcal{R}(k),k\geq 1)$ which is the candidate for the marginals of $%
\mathcal{T }_{F}$. By Aldous \cite{aldouscrt93}, it suffices to check two
properties, called \emph{consistency} and \emph{leaf-tightness}. Notice that
in \cite{aldouscrt93}, only binary trees (in which branchpoint have
out-degree $2$) are considered, but as noticed therein, this translates to
our setting with minor changes.

\begin{lmm}
(i) The family $(\mathcal{R}(k),k\geq 1)$ is consistent in the sense that
for every $k$ and $j\leq k$, $\mathcal{R}(j)$ has the same law as the
subtree of $\mathcal{R}(k)$ spanned by the root and $j$ distinct leaves $%
L^k_1,\ldots,L^k_j$ taken uniformly at random from the leaves $%
L_1,\ldots,L_k $ of $\mathcal{R}(k)$, independently of $\mathcal{R}(k)$.

(ii) The family $(\mathcal{R}(k),k\geq 1)$ is leaf-tight, that is, with the
above notations,
\begin{equation*}
\min_{2\leq j\leq k}d(L^k_1,L^k_j)\build\to_{}^{p}0.
\end{equation*}
\end{lmm}

{\noindent \textbf{Proof. }} The consistency property is an immediate
consequence of the fact that the process $\Pi $ is exchangeable. Taking $j$
leaves uniformly out of the $k$ ones of $\mathcal{R}(k)$ is just the same as
if we had chosen exactly the leaves $L_{1},L_{2},\ldots ,L_{j}$, which give
rise to the tree $\mathcal{R}(j)$, and this is (i).

For (ii), first notice that we may suppose by exchangeability that $%
L_1^k=L_1 $. The only point is then to show that the minimal distance of
this leaf to the leaves $L_{2},\ldots ,L_{k}$ tends to $0 $ in probability
as $k\rightarrow \infty $. Fix $\eta >0$ and for $\varepsilon>0$ write $%
t_{\varepsilon}^{1}=\inf \{t\geq 0:|\Pi _{1}(t)|<\varepsilon\}$, where $\Pi
_{1}(t)$ is the block of $\Pi (t)$ containing $1$. Then $t_{\varepsilon}^{1}$
is a stopping time with respect to the natural filtration $( \mathcal{F }%
_{t},t\geq 0)$ associated to $\Pi $ and $t_{\varepsilon}^{1}\uparrow D_{1}$
as $\varepsilon\downarrow 0$. By the strong Markov property and
exchangeability, one has that if $K(\varepsilon)=\inf \{k>1:k\in \Pi
_{1}(t_{\varepsilon}^{1})\}$, then $P(D_{1}+D_{K(\varepsilon)}-2t_{%
\varepsilon}^{1}<\eta )=E[P_{\Pi
(t_{\varepsilon}^{1})}(D_{1}+D_{K(\varepsilon)}<\eta )]$ where $P_{\pi }$ is
the law of the fragmentation $\Pi $ started at $\pi $ (the law of $\Pi $
under $P_{\pi }$ is the same as that of the family of partitions $\left(
\left\{ \text{blocks of }\pi _{1}\cap \Pi ^{(1)}(|\pi _{1}|^{\alpha }t),\pi
_{2}\cap \Pi ^{(2)}(|\pi _{2}|^{\alpha }t),\ldots \right\} ,t\geq 0\right) $
where the $\Pi ^{(i)}$'s$,i\geq 1,$ are independent copies of $\Pi $ under $%
P_{\{\mathbb{N},\varnothing ,\varnothing ,\ldots \}}$). By the self-similar
fragmentation property and exchangeability this is greater than $%
P(D_{1}+D_{2}<\varepsilon^{\alpha }\eta )$, which in turn is greater than $%
P(2\tau <\varepsilon^{\alpha }\eta )$ where $\tau $ is the first time where $%
\Pi (t)$ becomes the partition into singletons, which by \cite{bertafrag02}
is finite a.s. This last probability thus goes to $1$ as $%
\varepsilon\downarrow 0$. Taking $\varepsilon=\varepsilon(n)\downarrow 0$
quickly enough as $n\rightarrow \infty $ and applying the Borel-Cantelli
lemma, we a.s.\ obtain a sequence $K(\varepsilon(n))$ such that $%
d(L_{1},L_{K(n)})\leq D_{1}+D_{K(\varepsilon(n))}-2t_{\varepsilon(n)}<\eta $%
. Hence the result. $\hfill \square$

For a rooted $\mathbb{R}$-tree $T$ and $k$ vertices $v_{1},\ldots ,v_{k}$,
we define exactly as for marked trees the subtree spanned by the root and $%
v_{1},\ldots ,v_{k}$, as an element of $\mathbb{T}_{k}$. A consequence of
\cite[Theorem 3]{aldouscrt93} is then:

\begin{lmm}
\label{C1hm} There exists a CRT $(\mathcal{T}_{\Pi },\mu _{\Pi })$ such that
if $Z_{1},\ldots ,Z_{k}$ is a sample of $k$ leaves picked independently
according to $\mu _{\Pi }$ conditionally on $\mu _{\Pi }$, the subtree of $%
\mathcal{T}_{\Pi }$ spanned by the root and $Z_{1},\ldots ,Z_{k}$ has the
same law as $\mathcal{R}(k)$.
\end{lmm}

In the sequel, sequences like $(Z_1,Z_2,\ldots)$ will be called \emph{%
exchangeable sequences with directing measure $\mu_{\Pi}$}.

\noindent \textbf{Proof of Theorem \ref{T1}.}\emph{\ } We have to check that
the tree $\mathcal{T}_{\Pi }$ of the preceding lemma gives rise to a
fragmentation process with the same law as $F=|\Pi |$. By construction, we
have that for every $t\geq 0$ the partition $\Pi (t)$ is such that $i$ and $%
j $ are in the same block of $\Pi (t)$ if and only if $L_{i}$ and $L_{j}$
are in the same connected component of $\{v\in \mathcal{T}_{\Pi }:\mathrm{ht}%
(v)>t\}$. Hence, the law of large numbers implies that if $F^{\prime }(t)$
is the decreasing sequence of the $\mu $-masses of these connected
components, then $F^{\prime }(t)=F(t)$ a.s.\ for every $t$. Hence, $%
F^{\prime }$ is a version of $F$, so we can set $\mathcal{T}_{F}=\mathcal{T}%
_{\Pi }$. That $\mathcal{T}_{F}$ is $\alpha $-self-similar is an immediate
consequence of the fragmentation and self-similar properties of $F$.

We now turn to the last statement of Theorem \ref{T1}. With the notation of
Lemma \ref{C1hm} we will show that the path $[[\varnothing ,Z_{1}]]$ is
almost-surely in the closure of the set of leaves of $\mathcal{T }_{F}$ if
and only if $\nu (S)=\infty $. Then it must hold by exchangeability that so
do the paths $[[\varnothing ,Z_{i}]]$ for every $i\geq 1$, and this is
sufficient because the definition of the CRTs imply that $\mathcal{S}(
\mathcal{T }_{F})=\bigcup_{i\geq 1}[[\varnothing ,Z_{i}[[$, see \cite[Lemma
6]{aldouscrt93} (the fact that $\mathcal{T }_{F}$ is a.s.\ compact will be
proved below). To this end, it suffices to show that for any $a\in (0,1)$,
the point $aZ_{1}$ of $[[\varnothing ,Z_{1}]]$ that is at a proportion $a$
from $\varnothing $ (the point $\varphi _{\varnothing ,Z_{1}}(ad(\varnothing
,Z_{1}))$ with the above notations) can be approached closely by leaves,
that is, for $\eta >0$ there exists $j>1$ such that $d(aZ_{1},Z_{j})<\eta $.
It thus suffices to check that for any $\delta >0$
\begin{equation}  \label{sansordre}
P(\exists 2\leq j\leq k:|D_{\{1,j\}}-aD_{1}|<\delta \mbox{ and
}D_{j}-D_{\{1,j\}}<\delta )\build\rightarrow _{k\rightarrow \infty }^{{}}1,
\end{equation}
with the above notations derived from $\Pi $ (this is a slight variation of
\cite[(iii) a).\ Theorem 15]{aldouscrt93}).

Suppose that $\nu (S)=\infty $. Then for every rational $r>0$ such that $%
|\Pi _{1}(r)|\neq 0$ and for every $\delta >0$, the block containing $1$
undergoes a fragmentation in the time-interval $\left( r,r+\delta /2\right) $%
. This is obvious from the Poisson construction of the self-similar
fragmentation $\Pi $ given above, because $\nu $ is an infinite measure so
there is an infinite number of atoms of $(\Delta _{t},k_{t})$ with $k_{t}=1$
in any time-interval with positive length. It is then easy that there exists
an infinite number of elements of $\Pi _{1}(r)$ that are isolated in
singletons of $\Pi (r+\delta )$, e.g.\ because of Lemma \ref{compact} below
which asserts that only a finite number of the blocks of $\Pi (r+\delta /2)$
``survive'' at time $r+\delta $, i.e.\ is not completely reduced to
singletons. Thus, an infinite number of elements of $\Pi _{1}(r)$ correspond
to leaves of some $\mathcal{R}(k)$ for $k$ large enough. By taking $r$ close
to $aD_{1}$ we thus have the result.

On the other hand, if $\nu(S)<\infty$, it follows from the Poisson
construction that the state $(1,0,\ldots)$ is a holding state, so the first
fragmentation occurs at a positive time, so the root cannot be approached by
leaves. $\hfill \square$

\noindent\textbf{Remark. } We have seen that we may actually build
simultaneously the trees $(\mathcal{R}(k), k\geq 1)$ on the same probability
space as a measurable functional of the process $(\Pi(t),t\geq 0)$. This
yields, by redoing the ``special construction'' of Aldous \cite{aldouscrt93}%
, a \emph{stick-breaking construction} of the tree $\mathcal{T }_F$, by now
considering the trees $\mathcal{R}(k)$ as $\mathbb{R}$-trees obtained as
finite unions of segments rather than trees with edge-lengths (one can check
that it is possible to switch between the two notions). The mass measure is
then defined as the limit of the empirical measure on the leaves $%
L_1,\ldots,L_n$. The special CRT thus constructed is a subset of $\ell^1$ in
\cite{aldouscrt93}, but we consider it as universal, i.e.\ up to
isomorphism. The tree $\mathcal{R}(k+1)$ is then obtained from $\mathcal{R}%
(k)$ by branching a new segment with length $D_{k+1}-\max_{B\subset[k%
],B\neq\varnothing}D_{B\cup\{k\}}$, and $\mathcal{T }_F$ can be
reinterpreted as the completion of the metric space $\bigcup_{k\geq 1}%
\mathcal{R}(k)$. On the other hand, call $L_1,L_2,\ldots$ as before the
leaves of $\bigcup_{k\geq 1}\mathcal{R}(k)$, $L_k$ being the leaf
corresponding to the $k$-th branch. One of the subtleties of the special
construction of \cite{aldouscrt93} is that $L_1,L_2,\ldots$ is not itself an
exchangeable sample with the mass measure as directing law. However,
considering such a sample $Z_1,Z_2,\ldots$, we may construct a random
partition $\Pi^{\prime}(t)$ for every $t$ by letting $i\sim^{\Pi^{%
\prime}(t)}j$ if and only if $Z_i$ and $Z_j $ are in the same connected
component of the forest $\{v\in \mathcal{T }_F:\mathrm{ht}(v)>t\}$. Then
easily $\Pi^{\prime}(t)$ is again a partition-valued self-similar
fragmentation, and in fact $|\Pi^{\prime}(t)|=F(t)$ a.s.\ for every $t$ so $%
\Pi^{\prime}$ has same law as $\Pi$ ($\Pi^{\prime}$ can be interpreted as a
``relabeling'' of the blocks of $\Pi$). As a conclusion, up to this
relabeling, we may and will assimilate $\mathcal{T }_F$ as the completion of
the increasing union of the trees $\mathcal{R}(k)$, while $L_1,L_2,\ldots$
will be considered as an exchangeable sequence with directing law $\mu_F$.

\noindent \textbf{Proof of Proposition \ref{reciproque}.}\emph{\ } The fact
that the process $F$ defined out of a CRT $(\mathcal{T},\mu )$ with the
stated properties is a $S$-valued self-similar fragmentation with index $%
\alpha $ is straightforward and left to the reader. The treatment of the
erosion and sudden loss of mass is a little more subtle. Let $%
Z_{1},Z_{2},\ldots $ be an exchangeable sample directed by the measure $\mu $%
, and for every $t\geq 0$ define a random partition $\Pi (t)$ by saying that
$i$ and $j$ are in the same block of $\Pi (t)$ if $Z_{i}$ and $Z_{j}$ fall
in the same tree component of $\{v\in \mathcal{T}:\mathrm{ht}(v)>t\}$. By
the arguments above, $\Pi $ defines a self-similar partition-valued
fragmentation such that $|\Pi (t)|=F(t)$ a.s.\ for every $t$. Notice that if
we show that the erosion coefficient $c=0$ and that no sudden loss of mass
occur, it will immediately follow that $\mathcal{T}$ has the same law as $%
\mathcal{T}_{F}$.

Now suppose that $\nu(\sum_i s_i<1)\neq 0$. Then (e.g.\ by the Poisson
construction of fragmentations described above) there exists a.s.\ two
distinct integers $i$ and $j$ and a time $D$ such that $i$ and $j$ are in
the same block of $\Pi(D-)$ but $\{i\}\in \Pi(D)$ and $\{j\}\in \Pi(D)$.
This implies that $Z_i=Z_j$, so $\mu$ has a.s.\ an atom and $( \mathcal{T }%
,\mu)$ cannot be a CRT. On the other hand, suppose that the erosion
coefficient $c>0 $. Again from the Poisson construction, we see that there
a.s.\ exists a time $D$ such that $\{1\}\notin \Pi(D-)$ but $\{1\}\in \Pi(D)$%
, and nevertheless $\Pi(D)\cap\Pi_1(D-)$ is not the trivial partition $%
\mathbb{O}$. Taking $j$ in a non-trivial block of this last partition and
denoting its death time by $D^{\prime}$, we obtain that the distance from $%
Z_1$ to $Z_j$ is $D^{\prime}-D$, while the height of $Z_1$ is $D$ and that
of $Z_j$ is $D^{\prime}$. This implies that $Z_1$ is a.s.\ not in the set of
leaves of $\mathcal{T }$, again contradicting the definition of a CRT. $%
\hfill \square$

\section{Hausdorff dimension of $\mathcal{T}_{F}$}

\label{hausdorff}

Let $(M,d)$ be a compact metric space. For $\mathcal{E}\subseteq M$, the
Hausdorff dimension of $\mathcal{E}$ is the real number
\begin{equation}
\dim_{\cal H}(\mathcal{E}):=\inf \left\{ \gamma >0:m_{\gamma
}(\mathcal{E})=0\right\}=\sup \left\{ \gamma >0:m_{\gamma }(\mathcal{E}%
)=\infty \right\} ,  \label{100}
\end{equation}
where
\begin{equation}
m_{\gamma }(\mathcal{E}):=\sup_{\varepsilon >0}\inf \sum\limits_{i}\Delta
(E_{i})^{\gamma },  \label{101}
\end{equation}
the infimum being taken over all collections $(E_{i},i\geq 1)$ of subsets of
$\mathcal{E}$ with diameter $\Delta (E_{i})\leq \varepsilon $, whose union
covers $\mathcal{E}$. This dimension is meant to measure the ``fractal
size'' of the considered set. For background on this subject, we mention
\cite{Falconer} (in the case $M=\mathbb{R}^{n}$, but the generalization to
general metric spaces of the results we will need is straightforward). 

The goal of this Section is to prove Theorem 2 and more generally that
\begin{equation*}
\frac{\varrho }{\left| \alpha \right| }\leq \dim_{\cal H}
\left( \mathcal{L}(\mathcal{T}_{F})\right) \leq \frac{1}{\left| \alpha
\right| }\text{ a.s}.
\end{equation*}
where $\varrho $ is the $\nu $-dependent parameter defined by $\left( \ref
{219}\right) $. 
The proof is divided in the two usual upper and lower bound parts. In
Section 3.1, we first prove that $\mathcal{T}_{F}$ is indeed compact and
that $\dim_{\cal H}\left( \mathcal{L}(\mathcal{T}_{F})\right)
\leq 1/\left| \alpha \right| $ a.s., which is true without any extra
integrability assumption on $\nu .$ We then show that this upper bound
yields $\dim_{\cal H}\left( \mathrm{d}M_{F}\right) \leq
1\wedge 1/\left| \alpha \right| $ a.s. (Corollary \ref{coromass}), the
Hausdorff dimension of $\mathrm{d}M_{F}$ being defined as
\begin{equation*}
\dim_{\cal H}\left( \mathrm{d}M_{F}\right) :=\inf \left\{
\dim_{\cal H}(E):\mathrm{d}M_{F}(E)=1\right\} .
\end{equation*}
Sections 3.2 to 3.4 are devoted to the lower bound $\dim_{\cal H}
\left( \mathcal{L}(\mathcal{T}_{F})\right) \geq \varrho /\left| \alpha
\right| $ a.s.\ This is obtained by using appropriate subtrees of $\mathcal{T%
}_{F}$ (we will see that the most naive way to apply Frostman's energy
method with the mass measure $\mu _{F}$ fails in general). That Theorem 2
applies to stable trees is proved in Sect.\ \ref{checkst}.

\bigskip

\subsection{Upper bound}


We begin by stating the expected

\begin{lmm}
\label{compact} The tree $\mathcal{T }_F$ is a.s. compact.
\end{lmm}

{\noindent \textbf{Proof. }} For $t\geq 0$ and $\varepsilon>0$, denote by $%
N_{t}^{\varepsilon}$ the number of blocks of $\Pi (t)$ not reduced to
singletons that are not entirely reduced to dust at time $t+\varepsilon $.
We first prove that $N_{t}^{\varepsilon}$ is a.s.\ finite. Let $(\Pi
_{i}(t),i\geq 1)$ be the blocks of $\Pi (t),$ and $(\left| \Pi
_{i}(t)\right| ,i\geq 1)$, their respective asymptotic frequencies. For
integers $i$ such that $\left| \Pi _{i}(t)\right| >0$, that is $\Pi
_{i}(t)\neq \varnothing $ and $\Pi _{i}(t)$ is not reduced to a singleton,
let $\tau _{i}:=\inf \left\{ s>t:\Pi _{i}(t)\cap \Pi (s)=\mathbb{O}\right\} $
be the first time at which the block $\Pi _{i}(t)$ is entirely reduced to
dust. Applying the fragmentation property at time $t,$ we may write $\tau
_{i}$ as $\tau _{i}=t+\left| \Pi _{i}(t)\right| ^{\left| \alpha \right| }%
\widetilde{\tau }_{i}$ where $\widetilde{\tau }_{i}$ is a r.v.\ independent
of $\mathcal{G}\left( t\right) $ that has same distribution as $\tau =\inf
\{t\geq 0:\Pi (t)=\mathbb{O}\}$, the first time at which the fragmentation
is entirely reduced to dust. Now, fix $\varepsilon >0.$ The number of blocks
of $\Pi (t)$ that are not entirely reduced to dust at time $t+\varepsilon ,$
which could be \textit{a priori} infinite, is then given by
\begin{equation*}
N_{t}^{\varepsilon }=\sum_{i:\left| \Pi _{i}(t)\right| >0}\mathbf{1}%
_{\left\{ \left| \Pi _{i}(t)\right| ^{\left| \alpha \right| }\widetilde{\tau
}_{i}>\varepsilon \right\} }.
\end{equation*}
>From Proposition 15 in \cite{haas02}, we know that there exist two
constants $C_{1},C_{2}$ such that $P(\tau >t)\leq C_{1}e^{-C_{2}t}$ for all $%
t\geq 0.$ Consequently, for all $\delta >0,$%
\begin{eqnarray}
E\left[ N_{t}^{\varepsilon }\mid \mathcal{G}\left( t\right) \right] &\leq
&C_{1}\sum_{i:\left| \Pi _{i}(t)\right| >0}e^{-C_{2}\varepsilon \left| \Pi
_{i}(t)\right| ^{\alpha }}  \label{102} \\
&\leq &C(\delta )\varepsilon ^{-\delta }\sum_{i}\left| \Pi _{i}(t)\right|
^{\left| \alpha \right| \delta },  \notag
\end{eqnarray}
where $C(\delta )=\sup_{x\in \mathbb{R}^{+}}\left\{ C_{1}x^{\delta
}e^{-C_{2}x}\right\} <\infty .$ Since $\sum_{i}\left| \Pi _{i}(t)\right|
\leq 1$ a.s, this shows by taking $\delta =1/|\alpha |$ that $%
N_{t}^{\varepsilon }<\infty $ a.s. 

Let us now construct a covering of $\mathrm{supp\,}\left( \mu \right) $ with
balls of radius $5\varepsilon $. Recall that we may suppose that the tree $%
\mathcal{T}_{F}$ is constructed together with an exchangeable leaf sample $%
(L_{1},L_{2},\ldots )$ directed by $\mu _{F}$. For each $l\in \mathbb{N\cup }%
\left\{ 0\right\} $\textbf{,} we introduce the set
\begin{equation*}
B_{l}^{\varepsilon }=\left\{ k\in \mathbb{N}:\left\{ k\right\} \notin \Pi
(l\varepsilon ),\left\{ k\right\} \in \Pi (\left( l+1\right) \varepsilon
)\right\} ,
\end{equation*}
some of which may be empty when $\nu (S)<\infty ,$ since the tree is not
leaf-dense. For $l\geq 1,$ the number of blocks of the partition $%
B_{l}^{\varepsilon }\cap \Pi (\left( l-1\right) \varepsilon )$ of $%
B_{l}^{\varepsilon }$ is less than or equal to $N_{(l-1)\varepsilon
}^{\varepsilon }$ and so is a.s. finite. 
Since the fragmentation is entirely reduced to dust at time $\tau <\infty $
a.s., $N_{l\varepsilon }^{\varepsilon }$ is equal to zero for $l\geq \tau
/\varepsilon $ and then, defining
\begin{equation*}
N_{\varepsilon }:=\sum_{l=0}^{\left[ \tau /\varepsilon \right]
}N_{l\varepsilon }^{\varepsilon }
\end{equation*}
we have $N_{\varepsilon }<\infty $ a.s.\ ($\left[ \tau /\varepsilon \right] $
denotes here the largest integer smaller than $\tau /\varepsilon ).$ Now,
consider a finite random sequence of pairwise distinct integers $\sigma
(1),...,\sigma (N_{\varepsilon })$ such that 
for each $1\leq l\leq \left[ \tau /\varepsilon \right] $ and each non-empty
block of $B_{l}^{\varepsilon }\cap \Pi (\left( l-1\right) \varepsilon ),$
there is a $\sigma (i),1\leq i\leq N_{\varepsilon },$ in this block. Then
each leaf $L_{j}$ belongs then to a ball of center $L_{\sigma (i)},$ for an
integer $1\leq i\leq N_{\varepsilon },$ and of radius $4\varepsilon $.
Indeed, fix $j\geq 1.$ It is clear that the sequence $\left(
B_{l}^{\varepsilon }\right) _{l\in \mathbb{N}\cup \left\{ 0\right\} }$ forms
a partition of $\mathbb{N}$. Thus, there exists a unique block $%
B_{l}^{\varepsilon }$ containing $j$ and in this block we consider the
integer $\sigma (i)$ that belongs to the same block as $j$ in the partition $%
B_{l}^{\varepsilon }\cap \Pi (((l-1)\vee 0)\varepsilon )$.
By definition (see Section 2.3), the distance between the leaves $L_{j}$ and
$L_{\sigma (i)}$ is $d(L_{j},L_{\sigma (i)})=D_{j}+D_{\sigma
(i)}-2D_{\left\{ j,\sigma (i)\right\} }.$ By construction, $j$ and $\sigma
(i)$ belong to the same block of $\Pi (\left( \left( l-1\right) \vee
0\right) \varepsilon )$ and both die before $\left( l+1\right) \varepsilon .$
In other words, max($D_{j},D_{\sigma (i)})\leq \left( l+1\right) \varepsilon
$ and $D_{\left\{ j,\sigma (i)\right\} }\geq \left( \left( l-1\right) \vee
0\right) \varepsilon ,$ which implies that $d(L_{j},L_{\sigma (i)})\leq
4\varepsilon .$ Therefore, we have covered the set of leaves $\left\{
L_{j},j\geq 1\right\} $ by at most $N_{\varepsilon }$ balls of radius $%
4\varepsilon .$ Since the sequence $\left( L_{j}\right) _{j\geq 1}$ is dense
in $\mathrm{supp\,}\left( \mu \right) ,$ this induces by taking balls with
radius $5\varepsilon $ instead of $4\varepsilon $ a covering of $\mathrm{%
supp\,}\left( \mu \right) $ by $N_{\varepsilon }$ balls of radius $%
5\varepsilon .$ This holds for all $\varepsilon >0$ so $\mathrm{supp\,}%
\left( \mu \right) $ is a.s.\ compact. The compactness of $\mathcal{T}_{F}$
follows. $\hfill \square $

Let us now prove the upper bound for $\dim_{\cal H}(\mathcal{L%
}( \mathcal{T }_F))$. The difficulty for finding a ``good'' covering of the
set $\mathcal{L}( \mathcal{T }_F)$ is that as soon as $\nu$ is infinite,
this set is dense in $\mathcal{T }_F$, and thus one cannot hope to find its
dimension by the plain box-counting method, because the skeleton $\mathcal{S}%
( \mathcal{T }_F)$ has a.s.\ Hausdorff dimension $1$ as a countable union of
segments. However, we stress that the covering with balls of radius $%
5\varepsilon$ of the previous lemma is a good covering of the \emph{whole
tree}, because the box-counting method leads to the right bound $\mathrm{%
dim\,}_{\mathcal{H}}( \mathcal{T }_F)\leq (1/|\alpha|)\vee1$, and this is
sufficient when $|\alpha|<1$. When $|\alpha|\geq 1$ though, we may lose the
details of the structure of $\mathcal{L}( \mathcal{T }_F)$. We will thus try
to find a sharp ``cutset'' for the tree, motivated by the computation of the
dimension of leaves of discrete infinite trees.

\noindent\textbf{Proof of Theorem \ref{T2}: upper bound. } For every $i\in%
\mathbb{N}$ and $t\geq 0$ let $\Pi_{(i)}(t)$ be the block of $\Pi(t)$
containing $i$ and for $\varepsilon>0$ let
\begin{equation*}
t_i^{\varepsilon}=\inf\{t\geq 0:|\Pi_{(i)}(t)|<\varepsilon\}.
\end{equation*}
Define a partition $\Pi^{\varepsilon}$ by $i\sim^{\Pi^{\varepsilon}}j$ if
and only if $\Pi_{(i)}(t^{\varepsilon}_i)=\Pi_{(j)}(t^{\varepsilon}_j)$. One
easily checks that this random partition is exchangeable, moreover it has
a.s.\ no singleton. Indeed, notice that for any $i$, $\Pi_{(i)}(t^{%
\varepsilon}_i)$ is the block of $\Pi(t^{\varepsilon}_i)$ that contains $i$,
and this block cannot be a singleton because the process $%
(|\Pi_{(i)}(t)|,t\geq 0)$ reaches $0$ continuously. Therefore, $%
\Pi^{(\varepsilon)}$ admits asymptotic frequencies a.s., and these
frequencies sum to $1$. Then let
\begin{equation*}
\tau_i^{\varepsilon}=\sup_{j\in\Pi_{(i)}(t^{\varepsilon}_i)} \inf\{t\geq
t^{\varepsilon}_i:|\Pi_{(j)}(t)|=0\}-t^{\varepsilon}_i
\end{equation*}
be the time after $t^{\varepsilon}_i$ when the fragment containing $i$
vanishes entirely (notice that $\tau_i^{\varepsilon}=\tau_j^{\varepsilon}$
whenever $i\sim^{\Pi^\eps}j$). We also let $b^{\varepsilon}_i$ be the unique
vertex of $[[\varnothing,L_i]]$ at distance $t^{\varepsilon}_i$ from the
root, notice that again $b^{\varepsilon}_i=b^{\varepsilon}_j$ whenever $%
i\sim^{\Pi^{\varepsilon}}j$.

We claim that
\begin{equation*}
\mathcal{L}( \mathcal{T }_F)\subseteq \bigcup_{i\in \mathbb{N}}\overline{B}%
(b^{\varepsilon}_i,\tau^{\varepsilon}_i),
\end{equation*}
where $\overline{B}(v,r)$ is the closed ball centered at $v$ with radius $r$
in $\mathcal{T }_F$. Indeed, for $L\in\mathcal{L}( \mathcal{T }_F)$, let $b_L
$ be the vertex of $[[\varnothing,L]]$ with minimal height such that $\mu_F(
\mathcal{T }_{b_L})<\varepsilon$, where $\mathcal{T }_{b_L}$ is the fringe
subtree of $\mathcal{T }_F$ rooted at $b_L$. Since $b_L\in\mathcal{S}(
\mathcal{T }_F)$, $\mu_F( \mathcal{T }_{b_L})>0$ and there exists infinitely
many $i$'s with $L_i\in \mathcal{T }_{b_L}$. But then, it is immediate that
for any such $i$, $t^{\varepsilon}_i=\mathrm{ht}(b_L)=\mathrm{ht}%
(b^{\varepsilon}_i)$. Since $(L_i,i\geq 1)$ is dense in $\mathcal{L}(
\mathcal{T }_F)$, and since for every $j$ with $L_j\in \mathcal{T }%
_{b^{\varepsilon}_i}$ one has $d(b^{\varepsilon}_i,L_j)\leq
\tau^{\varepsilon}_i$ by definition, it follows that $L\in \overline{B}%
(b^{\varepsilon}_i,\tau^{\varepsilon}_i)$. Therefore, $(\overline{B}%
(b^{\varepsilon}_i,\tau^{\varepsilon}_i),i\geq 1)$ is a covering of $%
\mathcal{L}( \mathcal{T }_F)$.

The next claim is that this covering is fine as $\varepsilon\downarrow 0$,
namely
\begin{equation*}
\sup_{i\in\mathbb{N}}\tau^{\varepsilon}_i \build\to_{\eps\downarrow0}^{}0
\qquad \mbox { a.s.}
\end{equation*}
Indeed, if it were not the case, we would find $\eta>0$ and $i_n,n\geq 0$,
such that $\tau^{1/2^{n}}_{i_n}\geq\eta$ and $d(b^{1/2^{n}}_{i_n},L_{i_n})%
\geq \eta /2$ for every $n$. Since $\mathcal{T }_F$ is compact, we may
suppose up to extraction that $L_{i_n}\to v$ for some $v\in \mathcal{T }_F$.
Now, since $\mu_F( \mathcal{T }_{b^{1/2^{n}}_{i_n}})\leq 2^{-n}$, it follows
that we may find a vertex $b\in[[\varnothing,v]]$ at distance at least $%
\eta/4$ from $v$, such that $\mu_F( \mathcal{T }_b)=0$, and this does not
happen a.s.

To conclude, by the self-similarity property applied at the $(\mathcal{G}%
(t),t\geq0)$-stopping time $t^{\varepsilon}_i$, $\tau^{\varepsilon}_i$ has
the same law as $|\Pi_{(i)}(t^{\varepsilon}_i)|^{|\alpha|}\tau$, where $\tau$
has same law as $\inf\{t\geq 0:|\Pi(t)|=(0,0,\ldots)\}$ and is taken
independent of $|\Pi_{(i)}(t^{\varepsilon}_i)|$. Therefore,
\begin{equation}  \label{recouv}
E\left[\sum_{i:\Pi_{(i)}(t^{\varepsilon}_i)=\Pi_i(t^{\varepsilon}_i)}
(\tau^{\varepsilon}_i)^{1/|\alpha|} \right]=E[\tau^{1/|\alpha|}]E\left[%
\sum_{i\geq 1}|\Pi_i(t^{\varepsilon}_i)|\right] =E[\tau^{1/|\alpha|}]<\infty
\end{equation}
(we have just chosen one $i$ representing each class of $\Pi^{\varepsilon}$
above). The fact that $E[\tau^{1/|\alpha|}]$ is finite comes from the fact
that $\tau$ has exponential moments. Because our covering is a fine covering
as $\varepsilon\downarrow0$, it finally follows that (with the above
notations)
\begin{equation*}
m_{1/|\alpha|}(\mathcal{L}( \mathcal{T }_F))\leq
\liminf_{\varepsilon\downarrow0}\sum_{i:\Pi_{(i)}(t^{\varepsilon}_i)=%
\Pi_i(t^{\varepsilon}_i)} (\tau^{\varepsilon}_i)^{1/|\alpha|} \qquad %
\mbox{a.s.},
\end{equation*}
which is a.s.\ finite by (\ref{recouv}) and Fatou's Lemma. $\hfill \square$

\noindent {\textbf{Proof of Corollary \ref{coromass}.}} By Theorem \ref{T1},
the measure d$M_{F}$ has same law as d$\overline{W}_{\mathcal{T}_{F}},$ the
Stieltjes measure associated to the cumulative height profile $\overline{W}_{
\mathcal{T}_{F}}(t)=\mu _{F}\left\{ v\in \mathcal{T}_{F}:\text{ht}(v)\leq
t\right\},t\geq 0.$ To bound from above the Hausdorff dimension of 
$\d\overline{W}_{\mathcal{T}_{F}},$ note that
\begin{equation*}
\text{d}\overline{W}_{\mathcal{T}_{F}}(\text{ht}\left( \mathcal{L}\left(
\mathcal{T}_{F}\right) \right) )=\int_{\mathcal{T}_{F}}\mathbf{1}_{\left\{
\text{ht}(v)\in \text{ht}\left( \mathcal{L}\left( \mathcal{T}_{F}\right)
\right) \right\} }\mu _{F}(\text{d}v)=1
\end{equation*}
since $\mu _{F}(\mathcal{L}\left( \mathcal{T}_{F}\right) )=1$. By definition
of $\dim_{\cal H}\left( \text{d}\overline{W}_{\mathcal{T}%
_{F}}\right) , $ it is thus sufficient to show that $\dim_{\cal H}
\left( \text{ht}(\mathcal{L}\left( \mathcal{T}_{F}\right)
)\right) \leq 1/\left| \alpha \right| $ a.s$.$ To do so, remark that ht is
Lipschitz and that this property easily leads to
\begin{equation*}
\dim_{\cal H}\left( \text{ht}(\mathcal{L}\left( \mathcal{T}%
_{F}\right) )\right) \leq \dim_{\cal H}\left( \mathcal{L}%
\left( \mathcal{T}_{F}\right) \right) .
\end{equation*}
The conclusion hence follows from the majoration $\dim_{\cal H}
\left( \mathcal{L}\left( \mathcal{T}_{F}\right) )\right) \leq 1/\left|
\alpha \right| $ proved above. $\hfill \square$

\subsection{A first lower bound}

Recall that Frostman's energy method to prove that $\dim_{\cal H}
(\mathcal{E})\geq \gamma $ where $\mathcal{E}$ is a subset of a metric
space $(M,d)$ is to find a nonzero positive measure $\eta (\mathrm{d}x)$ on $%
\mathcal{E}$ such that $\int_{\mathcal{E}}\int_{\mathcal{E}}\frac{\eta (%
\mathrm{d}x)\eta (\mathrm{d}y)}{d(x,y)^{\gamma }}<\infty $. A naive approach
for finding a lower bound of the Hausdorff dimension of $\mathcal{T}_{F}$ is
thus to apply this method by taking $\eta =\mu _{F}$ and $\mathcal{E}=%
\mathcal{L}(\mathcal{T}_{F})$. The result states as follows.

\begin{lmm}
\label{firstm} For any fragmentation process $F$ satisfying the hypotheses
of Theorem \ref{T1}, one has
\begin{equation*}
\dim_{\cal H}(\mathcal{L}(\mathcal{T}_{F}))\geq \frac{A}{%
|\alpha |}\wedge \left( 1+\frac{\underline{p}}{|\alpha |}\right) ,
\end{equation*}
where
\begin{equation}
\underline{p}:=-\inf \left\{ q:\int_{S}\left( 1-\sum_{i\geq
1}s_{i}^{q+1}\right) \nu (\mathrm{d}\mathbf{s})>-\infty \right\} \in \lbrack
0,1],  \label{210}
\end{equation}
and
\begin{equation}
A:=\sup \left\{ a\leq 1:\int_{S}\sum_{1\leq i<j}s_{i}^{1-a}s_{j}\nu (\mathrm{%
d}\mathbf{s})<\infty \right\} \in \lbrack 0,1].  \label{278}
\end{equation}
\end{lmm}

{\noindent \textbf{Proof. }} By Lemma \ref{C1hm} (recall that $\left(
\mathcal{T}_{\Pi },\mu _{\Pi }\right) =\left( \mathcal{T}_{F},\mu
_{F}\right) $ by Theorem 1) we have
\begin{equation*}
\int_{\mathcal{T}_{F}}\int_{\mathcal{T}_{F}}\frac{\mu _{F}(\mathrm{d}x)\mu
_{F}(\mathrm{d}y)}{d(x,y)^{\gamma }}\overset{a.s.}{=}E\left[ \frac{1}{%
d(L_{1},L_{2})^{\gamma }}\left| \mathcal{T}_{F},\mu _{F}\right. \right]
\end{equation*}
so that
\begin{equation*}
E\left[ \int_{\mathcal{T}_{F}}\int_{\mathcal{T}_{F}}\frac{\mu _{F}(\mathrm{d}%
x)\mu _{F}(\mathrm{d}y)}{d(x,y)^{\gamma }}\right] =E\left[ \frac{1}{%
d(L_{1},L_{2})^{\gamma }}\right]
\end{equation*}
and by definition, $d(L_{1},L_{2})=D_{1}+D_{2}-2D_{\left\{ 1,2\right\} }$.
Applying the strong fragmentation property at the stopping time $D_{\left\{
1,2\right\} }$, we can rewrite $D_{1}$ and $D_{2}$ as
\begin{equation*}
\begin{array}{cc}
D_{1}=D_{\left\{ 1,2\right\} }+\lambda _{1}^{\left| \alpha \right|
}(D_{\left\{ 1,2\right\} })\widetilde{D}_{1} & \text{ \ }D_{2}=D_{\left\{
1,2\right\} }+\lambda _{2}^{\left| \alpha \right| }(D_{\left\{ 1,2\right\} })%
\widetilde{D}_{2}
\end{array}
\end{equation*}
where $\lambda _{1}(D_{\left\{ 1,2\right\} })$ (resp. $\lambda
_{2}(D_{\left\{ 1,2\right\} }))$ is the asymptotic frequency of the block
containing 1 (resp. 2) at time $D_{\left\{ 1,2\right\} }$ and $\widetilde{D}%
_{1}$ and $\widetilde{D}_{2}$ are independent with the same law as $D_{1}$
and independent of $\mathcal{G}\left( D_{\left\{ 1,2\right\} }\right) .$
Therefore,
\begin{equation*}
d(L_{1},L_{2})=\lambda _{1}^{\left| \alpha \right| }(D_{\left\{ 1,2\right\}
})\widetilde{D}_{1}+\lambda _{2}^{\left| \alpha \right| }(D_{\left\{
1,2\right\} })\widetilde{D}_{2},
\end{equation*}
and
\begin{equation}
E\left[ \frac{1}{d(L_{1},L_{2})^{\gamma }}\right] \leq 2E\left[ \lambda
_{1}^{\alpha \gamma }(D_{\left\{ 1,2\right\} })\mathbf{;}\lambda
_{1}(D_{\left\{ 1,2\right\} })\geq \lambda _{2}(D_{\left\{ 1,2\right\} })%
\right] E\left[ D_{1}^{-\gamma }\right] .  \label{7}
\end{equation}
By \cite[Lemma 2]{haas03} the first expectation in the right-hand side of
inequality $\left( \ref{7}\right) $ is finite as soon as $\left| \alpha
\right| \gamma <A$, while by \cite[Sect.\ 4.2.1]{haas02} the second
expectation is finite as soon as $\gamma <1+\underline{p}/\left| \alpha
\right| $. 
That $\dim_{\cal H}(\mathcal{L}(\mathcal{T}_{F}))\overset{a.s.%
}{\geq }\left( \left( A/\left| \alpha \right| \right) \wedge \left( 1+%
\underline{p}/\left| \alpha \right| \right) \right) $ follows. $\hfill
\square $

Let us now make a comment about this bound. For dislocation measures such
that $\nu (s_{N+1}>0)=0$ for some $N\geq 1,$ the constant $A$ equals $1$
since for all $a>-1,$%
\begin{equation*}
\int_{S}\sum\limits_{i<j}^{{}}s_{i}^{1+a}s_{j}\nu (\mathrm{d}\mathbf{s})\leq
\int_{S}(N-1)\sum_{2\leq j\leq N}s_{j}\nu (\mathrm{d}\mathbf{s})\leq
(N-1)\int_{S}\left( 1-s_{1}^{{}}\right) \nu (\mathrm{d}\mathbf{s})<\infty .
\end{equation*}
In such cases, if moreover $\underline{p}=1$, the ``naive'' lower
bound of Lemma \ref{firstm} is thus equal to $1/\left| \alpha
\right| $. A typical setting in which this holds is when $\nu
(S)<\infty $ and $\nu (s_{N+1}>0)=0$ and therefore, for such
dislocation measures the ``naive'' lower bound is also the best
possible.

\subsection{A subtree of $\mathcal{T }_F$ and a reduced fragmentation}

In the general case, in order to improve this lower bound, we will thus try
to transform the problem on $F$ into a problem on an auxiliary fragmentation
that satisfies the hypotheses above. The idea is as follows: fix an integer $%
N$ and $0<\varepsilon <1.$ Consider the subtree $\mathcal{T }%
_{F}^{N,\varepsilon}\subset \mathcal{T }_{F}$ constructed from $\mathcal{T}%
_{F}$ by keeping, at each branchpoint, the $N$ largest fringe subtrees
rooted at this branchpoint (that is the subtrees with the largest masses)
and discarding the others in order to yield a tree in which branchpoints
have out-degree at most $N$. Also, we remove the accumulation of
fragmentation times by discarding all the fringe subtrees rooted at the
branchpoints but the largest one, as soon as the proportion of its mass
compared to the others is larger than $1-\varepsilon$. Then there exists a
probability $\mu _{F}^{N,\varepsilon}$ such that $( \mathcal{T }%
_{F}^{N,\varepsilon},\mu _{F}^{N,\varepsilon})$ is a CRT, to which we will
apply the energy method.

Let us make the definition precise. Define $\mathcal{L}^{N,\varepsilon}%
\subset \mathcal{L}( \mathcal{T }_{F})$ to be the set of leaves $L$ such
that for every branchpoint $b\in \lbrack \lbrack \varnothing ,L]],$ $L\in
\mathcal{F}_{b}^{N,\varepsilon }$ with $\mathcal{F}_{b}^{N,\varepsilon }$
defined by
\begin{equation}
\left\{
\begin{array}{ll}
\mathcal{F}_{b}^{N,\varepsilon }= \mathcal{T }_{b}^{1}\cup \ldots \cup
\mathcal{T }_{b}^{N} & \mbox{if }\mu _{F}( \mathcal{T }_{b}^{1})/\mu
_{F}\left( \bigcup_{i\geq 1} \mathcal{T }_{b}^{i}\right) \leq 1-\varepsilon
\\
\mathcal{F}_{b}^{N,\varepsilon }= \mathcal{T }_{b}^{1} & \mbox{if }\mu _{F}(
\mathcal{T }_{b}^{1})/\mu _{F}\left( \bigcup_{i\geq 1} \mathcal{T }%
_{b}^{i}\right) >1-\varepsilon
\end{array}
\right. ,  \label{condtfneps}
\end{equation}
where $\mathcal{T }_{b}^{1}, \mathcal{T }_{b}^{2}\ldots $ are the connected
components of the fringe subtree of $\mathcal{T }_{F}$ rooted at $b$, from
whom $b$ has been removed (the connected components of $\{v\in \mathcal{T }%
_{F}:\mathrm{ht}(v)>b\}$) and ranked in decreasing order of $\mu _{F}$-mass.
Then let $\mathcal{T }_{F}^{N,\varepsilon}\subset \mathcal{T }_{F}$ be the
subtree of $\mathcal{T }_{F}$ spanned by the root and the leaves of $%
\mathcal{L}^{N,\varepsilon}$, i.e.
\begin{equation*}
\mathcal{T }_{F}^{N,\varepsilon}=\{v\in \mathcal{T }_{F}:\exists L\in
\mathcal{L}^{N,\varepsilon},v\in \lbrack \lbrack \varnothing ,L]]\}.
\end{equation*}
The set $\mathcal{T }_{F}^{N,\varepsilon}\subset \mathcal{T }_{F}$ is
plainly connected and closed in $\mathcal{T }_{F}$, thus an $\mathbb{R}$%
-tree.

Now let us try to give a sense to ``taking at random a leaf in $\mathcal{T }%
_{F}^{N,\varepsilon}$''. In the case of $\mathcal{T }_{F}$, it was easy
because, from the partition-valued fragmentation $\Pi $, it sufficed to look
at the fragment containing $1$ (or some prescribed integer). Here, it is not
difficult (as we will see later) that the corresponding leaf $L_{1}$ a.s.\
never belongs to $\mathcal{T }_{F}^{N,\varepsilon}$ when the dislocation
measure $\nu $ charges the set $\left\{ s_{1}>1-\varepsilon \right\} \cup
\left\{ s_{N+1}>0\right\} $. Therefore, we will have to use several random
leaves of $\mathcal{T }_{F}$. For any leaf $L\in \mathcal{L}( \mathcal{T }%
_{F})\setminus \mathcal{L}( \mathcal{T }_{F}^{N,\varepsilon})$ let $b(L)$ be
the highest vertex $v$ of $[[\varnothing ,L]]$ such that 
$v\in \mathcal{T}_F^{N,\varepsilon}$. Call it the branchpoint of $L$ and 
$\mathcal{T }_{F}^{N,\varepsilon}$.

Now take at random a leaf $Z_{1}$ of $\mathcal{T }_{F}$ with law $\mu _{F}$
conditionally on $\mu _{F}$, and define recursively a sequence $(Z_{n},n\geq
1)$ with values in $\mathcal{T }_{F}$ as follows. Let $Z_{n+1}$ be
independent of $Z_{1},\ldots ,Z_{n}$ conditionally on $( \mathcal{T }%
_{F},\mu _{F},b(Z_{n}))$, and take it with conditional law
\begin{equation*}
P(Z_{n+1}\in \cdot | \mathcal{T }_{F},\mu _{F},b(Z_{n}))=\mu _{F}(\cdot \cap
\mathcal{F}_{b(Z_{n})}^{N,\varepsilon})/\mu _{F}(\mathcal{F}%
_{b(Z_{n})}^{N,\varepsilon}).
\end{equation*}

\begin{lmm}
\label{sampneps} Almost surely, the sequence $(Z_{n},n\geq 1)$ converges to
a random leaf $Z^{N,\varepsilon }\in \mathcal{L}(\mathcal{T}%
_{F}^{N,\varepsilon })$. If $\mu _{F}^{N,\varepsilon }$ denotes the
conditional law of $Z^{N,\varepsilon }$ given $(\mathcal{T}_{F},\mu _{F})$,
then $(\mathcal{T}_{F}^{N,\varepsilon },\mu _{F}^{N,\varepsilon })$ is a
CRT, provided $\varepsilon $ is small enough.
\end{lmm}

To prove this and for later use we first reconnect this discussion to
partition-valued fragmentations. Recall from Sect.\ \ref{pvssf} the
construction of the homogeneous fragmentation $\Pi ^{0}$ with
characteristics $(0,0,\nu )$ out of a $\mathcal{P}_{\infty }\times \mathbb{N}
$-valued Poisson point process $((\Delta _{t},k_{t}),t\geq 0)$ with
intensity $\kappa _{\nu }\otimes \#$. For any partition $\pi \in \mathcal{P}%
_{\infty }$ that admits asymptotic frequencies whose ranked sequence is $%
\mathbf{s}$, write $\pi _{i}^{\downarrow }$ for the block of $\pi $ with
asymptotic frequency $s_{i}$ (with some convention for ties, e.g.\ taking
the order of least element). We define a function $\mathtt{GRIND}%
^{N,\varepsilon}:\mathcal{P}_{\infty }\rightarrow \mathcal{P}_{\infty }$
that reduces the smallest blocks of the partition to singletons as follows.
If $\pi $ does not admit asymptotic frequencies, let $\mathtt{GRIND}%
^{N,\varepsilon}(\pi )=\pi $, else let
\begin{equation*}
\mathtt{GRIND}^{N,\varepsilon }(\pi )=\left\{
\begin{array}{lc}
\left( \pi _{1}^{\downarrow },...,\pi _{N}^{\downarrow },\text{singletons}%
\right) & \text{ if }s_{1}\leq 1-\varepsilon \\
\left( \pi _{1}^{\downarrow },\text{singletons}\right) & \text{ if }%
s_{1}>1-\varepsilon .
\end{array}
\right.
\end{equation*}
Now for each $t\geq 0$ write $\Delta _{t}^{N,\varepsilon }=\mathtt{GRIND}%
^{N,\varepsilon}(\Delta _{t})$, so $((\Delta
_{t}^{N,\varepsilon},k_{t}),t\geq 0)$ is a $\mathcal{P}_{\infty }\times
\mathbb{N}$-valued Poisson point process with intensity measure $\kappa
_{\nu ^{N,\varepsilon }}\otimes \#,$ where $\nu ^{N,\varepsilon }$ is the
image of $\nu $ by the function
\begin{equation*}
\mathbf{s}\in S\mapsto \left\{
\begin{array}{l}
\left( s_{1},...,s_{N},0,...\right) \text{ if }s_{1}\leq 1-\varepsilon \\
\left( s_{1},0,...\right) \text{ if }s_{1}>1-\varepsilon .
\end{array}
\right.
\end{equation*}
From this Poisson point process we construct first a version $\Pi
^{0,N,\varepsilon }$ of the $\left( 0,0,\nu ^{N,\varepsilon }\right) $
fragmentation, as explained in Section 2.1. For every time $t\geq 0,$ the
partition $\Pi ^{0,N,\varepsilon }(t)$ is finer than $\Pi ^{0}(t)$ and the
blocks of $\Pi ^{0,N,\varepsilon }(t)$ non-reduced to singleton are blocks
of $\Pi ^{0}(t).$ Next, using the times-change $\left( \ref{1}\right) ,$ we
construct from $\Pi ^{0,N,\varepsilon }$ a version of the $\left( \alpha
,0,\nu ^{N,\varepsilon }\right) $ fragmentation, that we denote by $\Pi
^{N,\varepsilon }.$

Note that for dislocation measures $\nu $ such that $\nu ^{N,\varepsilon
}\left( \sum s_{i}<1\right) =0$, Theorem \ref{T2} is already proved, by
the previous subsection. For the rest of this subsection and next
subsection, we shall thus focus on dislocation measures $\nu $ such that $%
\nu ^{N,\varepsilon }\left( \sum s_{i}<1\right) >0.$ In that case, in $\Pi
^{0,N,\varepsilon}$ (unlike for $\Pi ^{0}$) each integer $i$ is eventually
isolated in a singleton a.s.\ within a sudden break and this is why a $\mu
_{F}$-sampled leaf on $\mathcal{T }_{F}$ cannot be in $\mathcal{T }%
_{F}^{N,\varepsilon}$, in other words, $\mu _{F}$ and $\mu
_{F}^{N,\varepsilon}$ are a.s.\ singular. Recall that we may build $\mathcal{%
T }_{F}$ together with an exchangeable $\mu _{F}$-sample of leaves $%
L_{1},L_{2},\ldots $ on the same probability space as $\Pi $ (or $\Pi ^{0}$%
). We are going to use a subfamily of $(L_{1},L_{2},\ldots )$ to build a
sequence with the same law as $(Z_{n},n\geq 1)$ built above. Let $i_{1}=1$
and
\begin{equation*}
i_{n+1}=\inf \{i>i_{n}:L_{i_{n+1}}\in \mathcal{F}_{b(L_{i_{n}})}^{N,%
\varepsilon}\}.
\end{equation*}
It is easy that $(L_{i_{n}},n\geq 1)$ has the same law as $(Z_{n},n\geq 1)$.
From this, we build a decreasing family of blocks $B^{0,N,\varepsilon
}(t)\in \Pi ^{0}(t)$, $t\geq 0$, by letting $B^{0,N,\varepsilon}(t)$ be the
unique block of $\Pi ^{0}(t)$ that contains all but a finite number of
elements of $\{i_{1},i_{2},\ldots \}$.

Here is a useful alternative description of $B^{0,N,\varepsilon}(t)$. Let $%
D_{i}^{0,N,\varepsilon }$ be the death time of $i$ for the fragmentation $%
\Pi ^{0,N,\varepsilon}$ that is
\begin{equation*}
D_{i}^{0,N,\varepsilon}=\inf \{t\geq 0:\left\{ i\right\} \in \Pi
^{0,N,\varepsilon }(t)\}.
\end{equation*}
By exchangeability the $D_{i}^{0,N,\varepsilon}$'s are identically
distributed and $D_{1}^{0,N,\varepsilon}=\inf \{t\geq 0:k_{t}=1\mbox{ and }%
\{1\}\in \Delta _{t}^{N,\varepsilon}\}$ so it has an exponential law with
parameter $\int_{S}(1-\sum_{i}s_{i})\nu ^{N,\varepsilon }($d${\mathbf{s}}).$
Then notice that $B^{0,N,\varepsilon}(t)$ is the block admitting $i_{n}$ as
least element when $D_{i_{n}}^{0,N,\varepsilon }\leq
t<D_{i_{n+1}}^{0,N,\varepsilon }$. Indeed, by construction we have
\begin{equation*}
i_{n+1}=\inf \{i\in
B^{0,N,\varepsilon}(D_{i_{n}}^{0,N,\varepsilon}-):\{i\}\notin \Pi
^{0,N,\varepsilon}(D_{i_{n}}^{0,N,\varepsilon})\}.
\end{equation*}
Moreover, the asymptotic frequency $\lambda _{1}^{0,N,\varepsilon}(t)$ of $%
B^{0,N,\varepsilon}(t)$ exists for every $t$ and equals the $\mu _{F}$-mass
of the tree component of $\{v\in \mathcal{T }_{F}:\mathrm{ht}(v)>t\}$
containing $L_{i_{n}}$ for $D_{i_{n}}^{0,N,\varepsilon }\leq
t<D_{i_{n+1}}^{0,N,\varepsilon }$.


Notice that at time $D^{0,N,\varepsilon}_{i_n}$, either one non-singleton
block coming from $B^{0,N,\varepsilon}(D^{0,N,\varepsilon}_{i_n}-)$, or up
to $N$ non-singleton blocks may appear; by Lemma \ref{usefull}, $%
B^{0,N,\varepsilon}(D^{0,N,\varepsilon}_{i_n})$ is then obtained by taking
at random one of these blocks with probability proportional to its size.

\noindent \textbf{Proof of Lemma \ref{sampneps}.}\emph{\ } For $t\geq 0$ let
$\lambda ^{0,N,\varepsilon }(t)=|B^{0,N,\varepsilon }(t)|$ and
\begin{equation}
T^{0,N,\varepsilon }(t):=\inf \left\{ u\geq 0:\int_{0}^{u}\left( \lambda
^{0,N,\varepsilon }(r)\right) ^{-\alpha }\mathrm{d}r>t\right\}  \label{2}
\end{equation}
and write $B^{N,\varepsilon }(t):=B^{0,N,\varepsilon }(T^{0,N,\varepsilon
}(t)),$ for $T^{0,N,\varepsilon }(t)<\infty $ and $B^{N,\varepsilon
}(t)=\varnothing $ otherwise, so for all $t\geq 0,$ $B^{N,\varepsilon
}(t)\in \Pi ^{N,\varepsilon }(t)$. Let also $D_{i_{n}}^{N,\varepsilon
}:=T^{0,N,\varepsilon }(D_{i_{n}}^{0,N,\varepsilon })$ be the death time of $%
i_{n}$ in the fragmentation $\Pi ^{N,\varepsilon }$. It is easy that $%
b_{n}=b(L_{i_{n}})$ is the branchpoint of the paths $[[\varnothing
,L_{i_{n}}]]$ and $[[\varnothing ,L_{i_{n+1}}]]$, so the path $[[\varnothing
,b_{n}]]$ has length $D_{i_{n}}^{N,\varepsilon }$. The ``edges'' $%
[[b_{n},b_{n+1}]]$, $n\in \mathbb{N}$, have respective lengths $%
D_{i_{n+1}}^{N,\varepsilon }-D_{i_{n}}^{N,\varepsilon },$ $n\in \mathbb{N}.$
Since the sequence of death times $(D_{i_{n}}^{N,\varepsilon },n\geq 1)$ is
increasing and bounded by $\tau $ (the first time at which $\Pi $ is
entirely reduced to singletons), the sequence $(b_{n},n\geq 1)$ is Cauchy,
so it converges by completeness of $\mathcal{T}_{F}$. Now it is easy that $%
D_{i_{n}}^{0,N,\varepsilon }\rightarrow \infty $ as $n\rightarrow \infty $
a.s., so $\lambda ^{0,N,\varepsilon }(t)\rightarrow 0$ as $t\rightarrow
\infty $ a.s.\ (see also the next lemma). Therefore, it is easy by the
fragmentation property that $d(L_{i_{n}},b_{n})\rightarrow 0$ a.s.\ so $%
L_{i_{n}}$ is also Cauchy, with the same limit, and that the limit has to be
a leaf which we denote $L^{N,\varepsilon }$ (of course it has same
distribution as the $Z^{N,\varepsilon }$ of the lemma's statement). The fact
that $L^{N,\varepsilon }\in \mathcal{T}_{F}^{N,\varepsilon }$ a.s.\ is
obtained by checking (\ref{condtfneps}), which is true since it is verified
for each branchpoint $b\in \lbrack \lbrack \varnothing ,b_{n}]]$ for every $%
n\geq 1$ by construction.

We now sketch the proof that $( \mathcal{T }_{F}^{N,\varepsilon},\mu
_{F}^{N,\varepsilon})$ is indeed a CRT, leaving details to the reader. We
need to show non-atomicity of $\mu _{F}^{N,\varepsilon}$, but it is clear
that when performing the recursive construction of $Z^{N,\varepsilon}$ twice
with independent variables, $(Z_{n},n\geq 1)$ and $(Z_{n}^{\prime },n\geq 1)$
say, there exists a.s.\ some $n$ such that $Z_{n}$ and $Z_{n}^{\prime }$ end
up in two different fringe subtrees rooted at some of the branchpoints $%
b_{n} $, provided that $\varepsilon$ is small enough so that $\nu
(1-s_{1}\geq \varepsilon)\neq 0$ (see also below the explicit construction
of two independently $\mu _{F}^{N,\varepsilon}$-sampled leaves). On the
other hand, all of the subtrees of $\mathcal{T }_{F}$ rooted at the
branchpoints of $\mathcal{T }_{F}^{N,\varepsilon}$ have positive $\mu _{F}$%
-mass, so they will end up being visited by the intermediate leaves used to
construct a $\mu _{F}^{N,\varepsilon}$-i.i.d.\ sample, so the condition $\mu
_{F}^{N,\varepsilon}(\{v\in \mathcal{T }_{F}^{N,\varepsilon}:[[\varnothing
,v]]\cap \lbrack \lbrack \varnothing ,w]]=[[\varnothing ,w]]\})>0$ for every
$w\in \mathcal{S}( \mathcal{T }_{F}^{N,\varepsilon})$ is satisfied. $\hfill
\square$

It will also be useful to sample two leaves $(L_{1}^{N,%
\varepsilon},L_{2}^{N,\varepsilon})$ that are independent with same
distribution $\mu _{F}^{N,\varepsilon}$ conditionally on $\mu
_{F}^{N,\varepsilon}$ out of the exchangeable family $L_{1},L_{2},\ldots$. A
natural way to do this is to use the family $(L_1,L_3,L_5,\ldots)$ to sample
the first leaf in the same way as above, and to use the family $%
(L_2,L_4,\ldots)$ to sample the other one. That is, let $%
j_{1}^{1}=1,j_{1}^{2}=2$ and define recursively $(j_{n}^{1},j_{n}^{2},n\geq
1)$ by letting
\begin{equation*}
\left\{
\begin{array}{ll}
j_{n+1}^{1}=\inf \{j\in 2\mathbb{N}+1,j>j_{n}^{1}:L_{j}\in \mathcal{F}_{b(
L_{j_n^1})}^{N,\varepsilon}\} &  \\
j_{n+1}^{2}=\inf \{j\in 2\mathbb{N},j>j_{n+1}^{1}:L_{j}\in \mathcal{F}_{b(
L_{j_n^2})}^{N,\varepsilon}\} &
\end{array}
\right. .
\end{equation*}
It is easy to check that $(L_{j_{n}^{1}},n\geq 1)$ and $(L_{j_{n}^{2}},n\geq
1)$ are two independent sequences distributed as $(Z_{1},Z_{2},\ldots )$ of
Lemma \ref{sampneps}. Therefore, these sequences a.s.\ converge to limits $%
L_{1}^{N,\varepsilon},L_{2}^{N,\varepsilon}$, and these are independent with
law $\mu _{F}^{N,\varepsilon}$ conditionally on $\mu _{F}^{N,\varepsilon}$.
We let $\mathcal{D}_{k}=\mathrm{ht}(L_{k}^{N,\varepsilon})$, $k=1,2$.

Similarly as above, for every $t\geq 0$ we let $B_{k}^{0,N,\varepsilon}(t)$,
$k=1,2$ (resp.\ $B_{k}^{N,\varepsilon}(t)$) be the block of $\Pi ^{0}(t)$
(resp.\ $\Pi (t)$) that contains all but the first few elements of $%
\{j_{1}^{k},j_{2}^{k},\ldots \}$, and we call $\lambda
_{k}^{0,N,\varepsilon}(t)$ (resp.\ $\lambda _{k}^{N,\varepsilon}(t)$) its
asymptotic frequency. Last, let $\mathcal{D}_{\{1,2\}}^{0}=\inf \{t\geq
0:B_{1}^{0,N,\varepsilon}(t)\cap B_{2}^{0,N,\varepsilon}(t)=\varnothing \}$
(and define similarly $\mathcal{D}_{\{1,2\}})$. Notice that for $t<\mathcal{D%
}_{\{1,2\}}^{0}$, we have $B_{1}^{0,N,\varepsilon}(t)=B_{2}^{0,N,%
\varepsilon}(t)$, and by construction the two least elements of the blocks $%
(2\mathbb{N}+1)\cap B^{0,N,\varepsilon}_1(t)$ and $(2\mathbb{N})\cap
B^{0,N,\varepsilon}_1(t) $ are of the form $j_{n}^{1},j_{m}^{2}$ for some $%
n,m$. On the other hand, for $t\geq \mathcal{D}_{\{1,2\}}^{0}$, we have $%
B^{0,N,\varepsilon}_1(t)\cap B^{0,N\varepsilon}_2(t)=\varnothing$, and again
the least elements of $(2\mathbb{N}+1)\cap B^{0,N,\varepsilon}_1(t)$ and $(2%
\mathbb{N})\cap B^{0,N\varepsilon}_2(t)$ are of the the form $%
j_{n}^{1},j_{m}^{2}$ for some $n,m$. In any case, we let $%
j^{1}(t)=j_{n}^{1},j^2(t)=j_m^2$ for these $n,m$.

\subsection{Lower bound}

Since $\mu _{F}^{N,\varepsilon }$ is a measure on $\mathcal{L}(\mathcal{T}%
_{F})$, we want to show that for every $a<\varrho ,$ the integral $\int_{%
\mathcal{T}_{F}^{N,\varepsilon }}\int_{\mathcal{T}_{F}^{N,\varepsilon }}%
\frac{\mu _{F}^{N,\varepsilon }(\mathrm{d}x)\mu _{F}^{N,\varepsilon }(%
\mathrm{d}y)}{d(x,y)^{a/\left| \alpha \right| }}$ is a.s. finite for
suitable $N$ and $\varepsilon $. So consider $a<\varrho $, and note that
\begin{equation*}
E\left[ \int_{\mathcal{T}_{F}^{N,\varepsilon }}\int_{\mathcal{T}%
_{F}^{N,\varepsilon }}\frac{\mu _{F}^{N,\varepsilon }(\mathrm{d}x)\mu
_{F}^{N,\varepsilon }(\mathrm{d}y)}{d(x,y)^{a/\left| \alpha \right| }}\right]
=E\left[ \frac{1}{d(L_{1}^{N,\varepsilon },L_{2}^{N,\varepsilon })^{a/\left|
\alpha \right| }}\right] ,
\end{equation*}
where $d(L_{1}^{N,\varepsilon },L_{2}^{N,\varepsilon })=\mathcal{D}_{1}+%
\mathcal{D}_{2}-2\mathcal{D}_{\{1,2\}},$ with notations above.
The fragmentation property at the stopping time $\mathcal{D}_{\{1,2\}}$ lead
to
\begin{equation*}
\mathcal{D}_{k}=\mathcal{D}_{\{1,2\}}+\lambda _{k}^{N,\varepsilon }(\mathcal{%
D}_{\{1,2\}})^{|\alpha |}\widetilde{\mathcal{D}}_{k},\text{ }k=1,2,
\end{equation*}
where $\widetilde{\mathcal{D}}_{1},\widetilde{\mathcal{D}}_{2}$ are
independent with the same distribution as $\mathcal{D},$ the height of the
leaf $L^{N,\varepsilon }$ constructed above, and independent of $\mathcal{G}(%
\mathcal{D}_{\{1,2\}}).$ Therefore, the distance $d(L_{1}^{N,\varepsilon
},L_{2}^{N,\varepsilon })$ can be rewritten as
\begin{equation*}
d(L_{1}^{N,\varepsilon },L_{2}^{N,\varepsilon })=\left( \lambda
_{1}^{N,\varepsilon }(\mathcal{D}_{\{1,2\}})\right) ^{\left| \alpha \right| }%
\widetilde{\mathcal{D}}_{1}+\left( \lambda _{2}^{N,\varepsilon }(\mathcal{D}%
_{\{1,2\}})\right) ^{\left| \alpha \right| }\widetilde{\mathcal{D}}_{2}
\end{equation*}
and
\begin{equation*}
E\left[ d(L_{1}^{N,\varepsilon },L_{2}^{N,\varepsilon })^{-a/\left| \alpha
\right| }\right] \leq 2E\left[ \left( \lambda _{1}^{N,\varepsilon }(\mathcal{%
D}_{\{1,2\}})\right) ^{-a};\lambda _{1}^{N,\varepsilon }(\mathcal{D}%
_{\{1,2\}})\geq \lambda _{2}^{N,\varepsilon }(\mathcal{D}_{\{1,2\}})\right] E%
\left[ \mathcal{D}^{-a/\left| \alpha \right| }\right] .
\end{equation*}
Therefore, that $\dim_{\cal H}(\mathcal{L}(\mathcal{T}%
_{F}))\geq 1\vee \left( \varrho /\left| \alpha \right| \right) $ is directly
implied by the following Lemmas \ref{Corosub1} and \ref{Corosub2}.

\begin{lmm}
\label{Corosub1} 
The quantity $E[\mathcal{D}^{-\gamma }]$ is finite for every $0\leq \gamma
\leq \varrho /\left| \alpha \right| .$
\end{lmm}

The proof uses the following technical lemma. Recall that $%
\lambda^{N,\varepsilon}(t)=|B^{N,\varepsilon}(t)|$.

\begin{lmm}
One can write $\lambda ^{N,\varepsilon }=\exp \left( -\xi _{\rho (\cdot
)}\right) ,$ where $\xi $ (tacitly depending on $N,\varepsilon $) is a
subordinator with Laplace exponent
\begin{equation}
\Phi _{\xi }(q)=\int_{S}\bigg(\left( 1-s_{1}^{q}\right) \mathbf{1}_{\left\{
s_{1}>1-\varepsilon \right\} }+\sum_{i=1}^{N}\left( 1-s_{i}^{q}\right) \frac{%
s_{i}\mathbf{1}_{\left\{ s_{1}\leq 1-\varepsilon \right\} }}{s_{1}+...+s_{N}}%
\bigg)\nu (\mathrm{d}\mathbf{s}),\text{ }q\geq 0,  \label{3}
\end{equation}
and $\rho $ is the time-change
\begin{equation*}
\rho (t)=\inf \left\{ u\geq 0:\int_{0}^{u}\exp (\alpha \xi _{r})\mathrm{d}%
r>t\right\} ,\text{ }t\geq 0.
\end{equation*}
\end{lmm}

{\noindent \textbf{Proof. }} Recall the construction of the process $%
B^{0,N,\varepsilon }$ from $\Pi ^{0}$, which itself was constructed from a
Poisson process $(\Delta _{t},k_{t},t\geq 0)$. From the definition of $%
B^{0,N,\varepsilon}(t)$, we have
\begin{equation*}
B^{0,N,\varepsilon}(t)=\bigcap_{0\leq s\leq t}\bar{\Delta}%
_{s}^{N,\varepsilon},
\end{equation*}
where the sets $\bar{\Delta}_{s}^{N,\varepsilon}$ are defined as follows.
For each $s\geq 0$, let $i(s)$ be the least element of the block $%
B^{0,N,\varepsilon}(s-)$ (so that $B^{0,N,\varepsilon}(s-)=\Pi
_{i(s)}^{0}(s-)$), so $(i(s),s\geq 0)$ is an $( \mathcal{F }(s-),s\geq 0)$%
-adapted jump-hold process, and the process $\{\Delta _{s}:k_{s}=i(s),s\geq
0\}$ is a Poisson point process with intensity $\kappa _{\nu }$. Then for
each $s$ such that $k_{s}=i(s)$, $\bar{\Delta}_{s}^{N,\varepsilon}$ consists
in a certain block of $\Delta _{s}$, and precisely,
$\bar{\Delta}_{s}^{N,\varepsilon}$ is the block of $\Delta _{s}$ containing
\begin{equation*}
\inf \left\{ i\in B^{0,N,\varepsilon}(s-):\{i\}\notin \Delta
_{s}^{N,\varepsilon}\right\} ,
\end{equation*}
the least element of $B^{0,N,\varepsilon}(s-)$ which is not isolated in a
singleton of $\Delta _{s}^{N,\varepsilon}$ (such an integer must be of the
form $i_{n}$ for some $n$ by definition). Now $B^{0,N,\varepsilon}(s-)$ is $%
\mathcal{F }(s-)$-measurable, hence independent of $\Delta _{s}$. By Lemma
\ref{usefull}, $\bar{\Delta}_{s}^{N,\varepsilon}$ is thus a size-biased pick
among the non-void blocks of $\Delta _{s}^{N,\varepsilon}$, and by
definition of the function $\mathtt{GRIND}^{N,\varepsilon}$, the process $(|%
\bar{\Delta}_{s}^{N,\varepsilon}|,s\geq 0)$ is a $[0,1]$-valued Poisson
point process with intensity $\omega (s)$ characterized by
\begin{equation*}
\int_{\lbrack 0,1]}f(s)\omega (\mathrm{d} s)=\int_{S}\left( \mathbf{1}%
_{\{s_{1}>1-\varepsilon\}}f(s_{1})+\mathbf{1}_{\{s_{1}\leq
1-\varepsilon\}}\sum_{i=1}^{N}f(s_{i})\frac{s_{i}}{s_{1}+\ldots +s_{N}}%
\right) \nu (\text{d}{\mathbf{s}}),
\end{equation*}
for every positive measurable function $f$. Then 
$|B^{0,N,\varepsilon}(t)|=\prod_{0\leq
s\leq t}|\bar{\Delta}_{s}^{N,\varepsilon}|$ a.s.\ for every $t\geq 0$. To
see this, denote for every $k\geq 1$ by $\Delta _{s_{1}}^{N,\varepsilon
,k},\Delta _{s_{2}}^{N,\varepsilon ,k},$... the atoms $\Delta
_{s}^{N,\varepsilon },$ $s\leq t,$ such that $|\Delta _{s}^{N,\varepsilon
}|_{1}\in \lbrack 1-k^{-1},1-(k+1)^{-1}).$ Complete this a.s. finite
sequence of partitions by partitions \textbf{1} and call $\Gamma ^{(k)}$
their intersection, i.e. $\Gamma ^{(k)}:=\bigcap_{i\geq 1}(\Delta
_{s_{i}}^{N,\varepsilon ,k})$. By Lemma \ref{Lemmafreqasym}$,$ $|\Gamma
_{n_{k}}^{(k)}|\overset{a.s.}{=}\prod_{i\geq 1}|\overline{\Delta }%
_{s_{i}}^{N,\varepsilon ,k}|,$ where $n_{k} $ is the index of the block $%
\bigcap_{i\geq 1}\overline{\Delta }_{s_{i}}^{N,\varepsilon ,k}$ in the
partition $\Gamma ^{(k)}.$ These partitions $\Gamma ^{(k)},$ $k\geq 1,$ are
exchangeable and clearly independent. Applying again Lemma \ref
{Lemmafreqasym} gives $|\bigcap_{k\geq 1}\Gamma _{n_{k}}^{(k)}|\overset{a.s.%
}{=}\prod_{k\geq 1}\prod_{i\geq 1}|\overline{\Delta }_{s_{i}}^{N,\varepsilon
,k}|,$ which is exactly the equality mentioned above. The exponential
formula for Poisson processes then shows that $(\xi _{t},t\geq 0)=(-\log
(\lambda ^{0,N,\varepsilon}(t)),t\geq 0)$ is a subordinator with Laplace
exponent $\Phi _{\xi }$. The result is now obtained by noticing that (\ref{1}%
) rewrites $\lambda ^{N,\varepsilon}(t)=\lambda ^{0,N,\varepsilon}(\rho (t))$
in our setting. $\hfill \square$



\noindent \textbf{Proof of Lemma \ref{Corosub1}.}\emph{\ } By the previous
lemma, $\mathcal{D}=\inf \{t\geq 0:\lambda ^{N,\varepsilon }(t)=0\}$, which
equals $\int_{0}^{\infty }\exp (\alpha \xi _{t})$d$t$ by definition of $\rho
$. According to Theorem 25.17 in \cite{Sato}, if for some positive $\gamma $
the quantity
\begin{equation*}
\Phi _{\xi }(-\gamma ):=\int_{S}\left( \left( 1-s_{1}^{-\gamma }\right)
\mathbf{1}_{\left\{ s_{1}>1-\varepsilon \right\} }+\sum_{i=1}^{N}\left(
1-s_{i}^{-\gamma }\right) \frac{s_{i}\mathbf{1}_{\left\{ s_{i}>0\right\} }%
\mathbf{1}_{\left\{ s_{1}\leq 1-\varepsilon \right\} }}{s_{1}+...+s_{N}}%
\right) \nu (\mathrm{d}\mathbf{s})
\end{equation*}
is finite, then $E[\exp (\gamma \xi _{t})]<\infty $ for all $t\geq 0$ and it
equals $\exp (-t\Phi _{\xi }(-\gamma )).$ Notice that $\Phi _{\xi }(-\gamma
)>-\infty $ for $\gamma <\varrho \leq 1$. Indeed for such $\gamma $'s, $%
\int_{S}\left( s_{1}^{-\gamma }-1\right) \mathbf{1}_{\left\{
s_{1}>1-\varepsilon \right\} }\nu (\mathrm{d}\mathbf{s})<\infty $ by
definition and
\begin{equation*}
\int_{S}\left( \sum_{i=1}^{N}\left( s_{i}^{1-\gamma }-s_{i}\right) \frac{%
\mathbf{1}_{\left\{ s_{1}\leq 1-\varepsilon \right\} }}{s_{1}+...+s_{N}}%
\right) \nu (\mathrm{d}\mathbf{s})\leq N\int_{S}\frac{s_{1}^{1-\gamma }%
\mathbf{1}_{\left\{ s_{1}\leq 1-\varepsilon \right\} }}{s_{1}}\nu (\mathrm{d}%
\mathbf{s}),
\end{equation*}
which is finite by definition of $\varrho $ and since $\nu $ integrates $%
\left( 1-s_{1}\right) $. This implies in particular that $\xi _{t}$ has
finite expectation for every $t$, and it follows by \cite{CPY} that $E[%
\mathcal{D}^{-1}]<\infty $. Then, following the proof of Proposition 2 in
\cite{BY} and using again that $\Phi _{\xi }(-\gamma )>-\infty $ for $\gamma
<\varrho ,$%
\begin{equation*}
E\left[ \left( \int_{0}^{\infty }\exp (\alpha \xi _{t})\mathrm{d}t\right)
^{-k-1}\right] =\frac{-\Phi _{\xi }(-\left| \alpha \right| k)}{k}E\left[
\left( \int_{0}^{\infty }\exp (\alpha \xi _{t})\mathrm{d}t\right) ^{-k}%
\right]
\end{equation*}
for every integer $k<\varrho /\left| \alpha \right| $. Hence, using
induction, $E[(\int_{0}^{\infty }\exp (\alpha \xi _{t}))^{-k-1}]$ is finite
for $k=[\varrho /|\alpha |]$ if $\varrho /|\alpha |\notin \mathbb{N}$ and
for $k=\varrho /|\alpha |-1$ else. In both cases, we see that $E[\mathcal{D}%
^{-\gamma }]<\infty $ for every $\gamma \leq \varrho /|\alpha |$. $\hfill
\square $

\begin{lmm}
\label{Corosub2} For any $a<\varrho $, there exists $N,\varepsilon $ such
that
\begin{equation*}
E\left[ \left( \lambda _{1}^{N,\varepsilon }(\mathcal{D}_{\{1,2\}})\right)
^{-a};\lambda _{1}^{N,\varepsilon }(\mathcal{D}_{\{1,2\}})\geq \lambda
_{2}^{N,\varepsilon }(\mathcal{D}_{\{1,2\}})\right] <\infty .
\end{equation*}
\end{lmm}

The ingredient for proving Lemma \ref{Corosub2} is the following lemma,
which uses the notations around the construction of the leaves $%
(L_{1}^{N,\varepsilon},L_{2}^{N,\varepsilon})$.

\begin{lmm}
\label{Coro2interm} With the convention $\log (0)=-\infty $, the process
\begin{equation*}
\sigma (t)=-\log \left| B_{1}^{0,N,\varepsilon }(t)\cap
B_{2}^{0,N,\varepsilon }(t)\right| \quad ,\quad t\geq 0
\end{equation*}
is a killed subordinator (its death time is $\mathcal{D}_{\{1,2\}}^{0}$)
with Laplace exponent
\begin{equation}
\Phi _{\sigma }(q)=\mathtt{k}^{N,\varepsilon }+\int_{S}\bigg(\left(
1-s_{1}^{q}\right) \mathbf{1}_{\left\{ s_{1}>1-\varepsilon \right\}
}+\sum_{i=1}^{N}\left( 1-s_{i}^{q}\right) \frac{s_{i}^{2}\mathbf{1}_{\left\{
s_{1}\leq 1-\varepsilon \right\} }}{\left( s_{1}+...+s_{N}\right) ^{2}}\bigg)%
\nu (\mathrm{d}\mathbf{s}),\text{ }q\geq 0\text{,}  \label{4}
\end{equation}
where the killing rate $\mathtt{k}^{N,\varepsilon }:=\int_{S}\sum_{i\neq
j}s_{i}s_{j}\frac{\mathbf{1}_{\left\{ s_{1}\leq 1-\varepsilon \right\} }}{%
\left( s_{1}+...+s_{N}\right) ^{2}}\nu (\mathrm{d}\mathbf{s})\in \left(
0,\infty \right) .$ Moreover, the pair
\begin{equation*}
(l_{1}^{N,\varepsilon },l_{2}^{N,\varepsilon })=\exp (\sigma (\mathcal{D}%
_{\{1,2\}}^{0}-))(\lambda _{1}^{0,N,\varepsilon }(\mathcal{D}%
_{\{1,2\}}^{0}),\lambda _{2}^{0,N,\varepsilon }(\mathcal{D}_{\{1,2\}}^{0}))
\end{equation*}
is independent of $\sigma (\mathcal{D}_{\{1,2\}}^{0}-)$ with law
characterized by
\begin{equation*}
E\left[ f\left( l_{1}^{N,\varepsilon },l_{2}^{N,\varepsilon }\right) \right]
=\frac{1}{\mathtt{k}^{N,\varepsilon }}\int_{S}\sum_{1\leq i\neq j\leq
N}f(s_{i},s_{j})\frac{s_{i}s_{j}\mathbf{1}_{\left\{ s_{1}\leq 1-\varepsilon
\right\} }\mathbf{1}_{\left\{ s_{i}>0\right\} }\mathbf{1}_{\left\{
s_{j}>0\right\} }}{\left( s_{1}+...+s_{N}\right) ^{2}}\nu (\mathrm{d}\mathbf{%
s})
\end{equation*}
for any positive measurable function $f$.
\end{lmm}

{\noindent \textbf{Proof. }}
We again use the Poisson construction of $\Pi ^{0}$ out of $(\Delta
_{t},k_{t},t\geq 0)$ and follow closely the proof of the intermediate lemma
used in the proof of Lemma \ref{Corosub1}. For every $t\geq 0$ we have
\begin{equation*}
B_{k}^{0,N,\varepsilon}(t)=\bigcap_{0\leq s\leq t}\bar{\Delta}_{s}^{k}\quad
,\quad k=1,2,
\end{equation*}
where $\bar{\Delta}_{s}^{k}$ is defined as follows. Let $J^{k}(s),k=1,2$ be
the integers such that $B_{k}^{0,N,\varepsilon}(s-)=\Pi _{J^{k}(s)}^{0}(s-)$%
, so $\{\Delta _{s}:k_{s}=J^{k}(s),s\geq 0\}$, $k=1,2$ are two Poisson
processes with same intensity $\kappa _{\nu }$, which are equal for $s$ in
the interval $[0,\mathcal{D}_{\{1,2\}}^{0})$. Then for $s$ with $%
k_{s}=J^{k}(s)$, let $\bar{\Delta}_{s}^{k}$ be the block of $\Delta _{s}$
containing $j^{k}(s)$. If $B_{1}^{0,N,\varepsilon}(s-)=B_{2}^{0,N,%
\varepsilon}(s-)$ notice that $j^{1}(s),j^{2}(s)$ are the two least integers
of $(2\mathbb{N}+1)\cap B_{1}^{0,N,\varepsilon}(s-)$ and $(2\mathbb{N})\cap
B_2^{0,N,\varepsilon}(s-)$ respectively that are not isolated as singletons
of $\Delta _{s}^{N,\varepsilon}$, so $\bar{\Delta}_{s}^{1}=\bar{\Delta}%
_{s}^{2}$ if these two integers fall in the same block of $\Delta
_{s}^{N,\varepsilon}$. Hence by a variation of Lemma \ref{usefull}, $(|\bar{%
\Delta}_{s}^{1}\cap \bar{\Delta}_{s}^{2}|,s\geq 0)$ is a Poisson process
whose intensity is the image measure of $\kappa _{\nu ^{N,\varepsilon}}(\pi
\mathbf{1}_{\{1\sim 2\}})$ by the map $\pi \mapsto |\pi |$, and killed at an
independent exponential time (namely $\mathcal{D}_{\{1,2\}}^{0})$ with
parameter $\kappa _{\nu ^{N,\varepsilon}}(1\nsim 2)$ (here $1\sim 2$ means
that $1$ and $2$ are in the same block of $\pi $). This implies (\ref{4}).

The time $\mathcal{D}_{\{1,2\}}^{0}$ is the first time when the two
considered integers fall into two distinct blocks of $\Delta
_{s}^{N,\varepsilon}$. It is then easy by the Poissonian construction and
the paintbox representation to check that these blocks have asymptotic
frequencies $(l_{1}^{N,\varepsilon},l_{2}^{N,\varepsilon})$ which are
independent of $\sigma (\mathcal{D}_{\{1,2\}}^{0}-)$, and have the claimed
law. $\hfill \square$

\noindent \textbf{Proof of Lemma \ref{Corosub2}.}\emph{\ } First notice,
from the fact that self-similar fragmentations are time-changed homogeneous
fragmentations, that
\begin{equation*}
(\lambda _{1}^{N,\varepsilon }(\mathcal{D}_{\{1,2\}}),\lambda
_{2}^{N,\varepsilon }(\mathcal{D}_{\{1,2\}}))\build=_{{}}^{d}(\lambda
_{1}^{0,N,\varepsilon }(\mathcal{D}_{\{1,2\}}^{0}),\lambda
_{2}^{0,N,\varepsilon }(\mathcal{D}_{\{1,2\}}^{0})).
\end{equation*}
Thus, with the notations of the intermediate lemma,
\begin{eqnarray*}
\lefteqn{E\left[ \left( \lambda _{1}^{N,\varepsilon }(\mathcal{D}%
_{\{1,2\}})\right) ^{-a};\lambda _{1}^{N,\varepsilon }(\mathcal{D}%
_{\{1,2\}})\geq \lambda _{2}^{N,\varepsilon }(\mathcal{D}_{\{1,2\}})\right] }
\\
&=&E\left[ \text{exp}(a\sigma (\mathcal{D}_{\{1,2\}}^{0}-)\right] E\left[
\left( l_{1}^{N,\varepsilon }\right) ^{-a};l_{1}^{N,\varepsilon }\geq
l_{2}^{N,\varepsilon }\right] .
\end{eqnarray*}
First, define for every $a>0$ $\ \Phi _{\sigma }(-a)$ by replacing $q$ by $%
-a $ in $\left( \ref{4}\right) $ and then remark that $\Phi _{\sigma
}(-a)>-\infty $ when $a<\varrho .$ Indeed, $\int_{S}\left(
s_{1}^{-a}-1\right) \mathbf{1}_{\left\{ s_{1}>1-\varepsilon \right\} }\nu (%
\mathrm{d}{\mathbf{s}})$ is then finite and, since $\sum_{1\leq i\leq
N}s_{i}^{2-a}\leq \left( \sum_{1\leq i\leq N}s_{i}\right) ^{2-a}$ $\ \
(2-a\geq 1),$%
\begin{equation*}
\sum_{1\leq i\leq N}\left( s_{i}^{2-a}-s_{i}^{2}\right) \frac{\mathbf{1}%
_{\left\{ s_{1}\leq 1-\varepsilon \right\} }}{\left( s_{1}+...+s_{N}\right)
^{2}}\leq \frac{\mathbf{1}_{\left\{ s_{1}\leq 1-\varepsilon \right\} }}{%
s_{1}^{a}}
\end{equation*}
which, by assumption, is integrable with respect to $\nu $. Then, consider a
subordinator $\widetilde{\sigma }$ with Laplace transform $\Phi _{\sigma }-%
\mathtt{k}^{N,\varepsilon }$ and independent of $\mathcal{D}_{\{1,2\}}^{0}, $
such that $\sigma =\widetilde{\sigma }$ on $(0,\mathcal{D}_{\{1,2\}}^{0}).$
As in the proof of Lemma \ref{Corosub1}, we use Theorem 25.17 of \cite{Sato}%
, which gives $E\left[ \exp (a\widetilde{\sigma }(t))\right] =\exp \left(
-t\left( \Phi _{\sigma }(-a)-\mathtt{k}^{N,\varepsilon }\right) \right) $
for all $t\geq 0.$ Hence, by independence of $\widetilde{\sigma }$ and $%
\mathcal{D}_{\{1,2\}}^{0}$,
\begin{eqnarray*}
E\left[ \exp (a\sigma(\mathcal{D}_{\{1,2\}}^{0}-)\right]&=& E\left[\exp(a%
\widetilde{\sigma}(\mathcal{D}^0_{\{1,2\}}))\right] \\
& =&\mathtt{k}^{N,\varepsilon } \int_{0}^{\infty }\exp (-t\mathtt{k}%
^{N,\varepsilon })\exp \left( -t(\Phi _{\sigma }(-a)-\mathtt{k}%
^{N,\varepsilon })\right) \mathrm{d}t,
\end{eqnarray*}
which is finite if and only if $\Phi _{\sigma }(-a)>0.$ Recall that $\Phi
_{\sigma }(-a)$ is equal to
\begin{equation}
\int_{S}\left( 1-s_{1}^{-a}\right) \mathbf{1}_{\left\{ s_{1}>1-\varepsilon
\right\} }\nu (\mathrm{d}{\mathbf{s}})+\int_{S}\left( \sum_{1\leq i\neq
j\leq N}\!\!\!s_{i}s_{j}+\sum_{1\leq i\leq N}\!\!\left(
s_{i}^{2}-s_{i}^{2-a}\right) \right) \frac{\mathbf{1}_{\left\{ s_{1}\leq
1-\varepsilon \right\} }}{\left( s_{1}+...+s_{N}\right) ^{2}}\nu (\mathrm{d}{%
\mathbf{s}}).  \label{6}
\end{equation}
Since
\begin{equation*}
\sum_{1\leq i\neq j\leq N}s_{i}s_{j}+\sum_{1\leq i\leq
N}(s_{i}^{2}-s_{i}^{2-a})=(\sum_{1\leq i\leq N}s_{i})^{2}-\sum_{1\leq i\leq
N}s_{i}^{2-a},
\end{equation*}
the integrand in the second term converges to $\left(
1-\sum_{i}s_{i}^{2-a}\right) \mathbf{1}_{\left\{ s_{1}\leq 1-\varepsilon
\right\} }$ as $N\rightarrow \infty $ and is dominated by $\left(
1+s_{1}^{-a}\right) \mathbf{1}_{\left\{ s_{1}\leq 1-\varepsilon \right\} }.$
So, by dominated convergence, the second term of $\left( \ref{6}\right) $
converges to $\int_{S}(1-\sum_{i}s_{i}^{2-a})\mathbf{1}_{\{s_{1}\leq
1-\varepsilon \}}\nu ($d${\mathbf{s}})$ as $N\rightarrow \infty $. This last
integral converges to a strictly positive quantity as $\varepsilon\downarrow
0$, and since $\int_{S}\left( 1-s_{1}^{-a}\right) \mathbf{1}_{\left\{
s_{1}>1-\varepsilon \right\} }\nu (\mathrm{d}{\mathbf{s}})\rightarrow 0$ as $%
\varepsilon \rightarrow 0,$ $\Phi _{\sigma }(-a)$ is strictly positive for $%
N $ and $1/\varepsilon $ large enough. Hence $E[\exp (a\sigma (\mathcal{D}%
_{\{1,2\}}^{0}-))]<\infty $ for $N$ and $1/\varepsilon $ large enough.

On the other hand, Lemma \ref{Coro2interm} implies that the finiteness of $%
E[(l_{1}^{N,\varepsilon })^{-a}\mathbf{1}_{\{l_{1}^{N,\varepsilon }\geq
l_{2}^{N,\varepsilon }\}}]$ is equivalent to that of $\int_{S}\sum_{1\leq
i\neq j\leq N}s_{i}^{1-a}s_{j}\frac{\mathbf{1}_{\{s_{1}\leq 1-\varepsilon \}}%
}{(s_{1}+...+s_{N})^{2}}\nu ($d${\mathbf{s}})$. But this integral is finite
for every integers $N$ and every $0<\varepsilon <1,$ since $\sum_{1\leq
i\neq j\leq N}s_{i}^{1-a}s_{j}\leq N^{2}s_{1}^{2-a}$ and $\nu $ integrates $%
s_{1}^{-a}\mathbf{1}_{\left\{ s_{1}\leq 1-\varepsilon \right\} }.$ Hence the
result. $\hfill \square $

\subsection{Dimension of the stable tree}

\label{checkst}

This section is devoted to the proof of Corollary \ref{stable}. Recall from
\cite{mierfmoins} that the fragmentation $F_{-}$ associated to the $\beta $%
-stable tree has index $1/\beta -1$ (where $\beta \in (1,2]$). In the case $%
\beta =2$, the tree is the Brownian CRT and the fragmentation is binary (it
is the fragmentation $F_{B}$ of the Introduction), so that the integrability
assumption of Theorem 2 is satisfied and then the dimension is $2$. So
suppose $\beta <2$. The main result of \cite{mierfmoins} is that the
dislocation measure $\nu _{-}(\mathrm{d}\mathbf{s})$ of $F_{-}$ has the form
\begin{equation*}
\nu _{-}(\mathrm{d}\mathbf{s})=C(\beta )E\left[ T_{1};\frac{\Delta T_{[0,1]}%
}{T_{1}}\in \mathrm{d}\mathbf{s}\right]
\end{equation*}
for some constant $C(\beta )$, where $(T_{x},x\geq 0)$ is a stable
subordinator with index $1/\beta $ and $\Delta T_{[0,1]}=(\Delta _{1},\Delta
_{2},\ldots )$ is the decreasing rearrangement of the sequence of jumps of $%
T $ accomplished within the time-interval $[0,1]$ (so that $\sum_{i}\Delta
_{i}=T_{1}$). By Theorem \ref{T2}, to prove Corollary \ref{stable} it thus
suffices to check that $E[T_{1}(T_{1}/\Delta _{1}-1)]$ is finite. The
problem is that computations involving jumps of subordinators are often
quite involved; they are sometimes eased by using size-biased picked jumps,
whose laws are more tractable. However, one can check that if $\Delta _{\ast
}$ is a size-biased picked jump among $(\Delta _{1},\Delta _{2},\ldots )$,
the quantity $E[T_{1}(T_{1}/\Delta _{\ast }-1)]$ is infinite, therefore we
really have to study the joint law of $(T_{1},\Delta _{1})$. This has been
done in Perman \cite{perm93}, but we will re-explain all the details we need
here.

Recall that the process $(T_{x},x\geq 0)$ can be put in the L\'{e}vy-It\^{o}
form $T_{x}=\sum_{0\leq y\leq x}\Delta (y)$, where $(\Delta (y),y\geq 0)$ is
a Poisson point process with intensity $cu^{-1-1/\beta }\mathrm{d}u$ (the
L\'{e}vy measure of $T$) for some constant $c>0$. Therefore, the law of the
largest jump of $T$ before time $1$ is characterized by
\begin{equation*}
P(\Delta _{1}<v)=P\left( \sup_{0\leq y\leq 1}\Delta (y)<v\right) =\exp
\left( -c\beta v^{-1/\beta }\right) \quad \quad v>0,
\end{equation*}
and by the restriction property of Poisson processes, conditionally on $%
\Delta _{1}=v$, one can write $T_{1}=v+T_{1}^{(v)}$, where $%
(T_{x}^{(v)},x\geq 0)$ is a subordinator with L\'{e}vy measure $%
cu^{-1-1/\beta }\mathbf{1}_{\{0\leq u\leq v\}}$d$u$. The Laplace transform
of $T_{x}^{(v)}$ is given by the L\'{e}vy-Khintchine formula
\begin{equation*}
E[\exp (-\lambda T_{x}^{(v)})]=\exp \left( -x\int_{0}^{v}\frac{%
c(1-e^{-\lambda u})}{u^{1+1/\beta }}\mathrm{d}u\right) \quad \lambda ,x\geq
0,
\end{equation*}
in particular, $T_{1}^{(v)}$ admits moments of all order (by differentiating
in $\lambda $) and $v^{-1}T_{1}^{(v)}$ has the same law as $T_{v^{-1/\beta
}}^{(1)}$ (by changing variables). We then obtain
\begin{eqnarray*}
E[T_{1}(T_{1}/\Delta _{1}-1)] &=&E\left[ \Delta _{1}\left( 1+\frac{%
T_{1}^{(\Delta _{1})}}{\Delta _{1}}\right) \frac{T_{1}^{(\Delta _{1})}}{%
\Delta _{1}}\right] \\
&=&K_{1}\int_{\mathbb{R}_{+}}\mathrm{d}v\,v^{-1/\beta }e^{-\beta
cv^{-1/\beta }}E\left[ \left( 1+\frac{T_{1}^{(v)}}{v}\right) \frac{%
T_{1}^{(v)}}{v}\right] \\
&=&K_{1}\int_{\mathbb{R}_{+}}\mathrm{d}v\,v^{-1/\beta }e^{-\beta
cv^{-1/\beta }}E\left[ \left( 1+T_{v^{-1/\beta }}^{(1)}\right)
T_{v^{-1/\beta }}^{(1)}\right]
\end{eqnarray*}
where $K_{1}=K(\beta )>0$. Since $T_{1}^{(1)}$ has a moment of orders $1$
and $2$, the expectation in the integrand is dominated by some $%
K_{2}v^{-1/\beta }+K_{3}v^{-2/\beta }$. It is then easy that the integrand
is integrable both near $0$ and $\infty $ since $\beta <2.$ Hence $%
\int_{S}\left( s_{1}^{-1}-1\right) \nu _{-}(\mathrm{d}\mathbf{s})<\infty .$

\section{The height function}

\label{heightprocess}

We now turn to the proof of the results related to the height function,
starting with Theorem \ref{T3}. The height function we are going to build
will in fact satisfy more than stated there: we will show that under the
hypotheses of Theorem \ref{T3}, there exists a process $H_F$ that encodes $%
\mathcal{T }_F$ in the sense given in the introduction, that is, $\mathcal{T
}_F$ is isometric to the quotient $((0,1),\overline{d})/\equiv$, where $%
\overline{d}(u,v)=H_F(u)+H_F(v)-2\inf_{s\in[u,v]}H_F(s)$ and $u\equiv v\iff
\overline{d}(u,v)=0$. Once we have proved this, the result is obvious since $%
I_F(t)/\equiv$ is the set of vertices of $\mathcal{T }_F$ that are above
level $t$.

\subsection{Construction of the height function}

Recall from \cite{aldouscrt93} that to encode a CRT, defined as a projective
limit of consistent random $\mathbb{R}-$trees $(\mathcal{R}(k),k\geq 1)$, in
a continuous height process, one first needs to enrich the structure of the $%
\mathbb{R}$-trees with consistent \emph{orders} on each set of children of
some node. The sons of a given node of $\mathcal{R}(k)$ are thus labelled as
first, second, etc... This induces a \emph{planar} representation of the
tree. This representation also induces a total order on the vertices of $%
\mathcal{R}(k)$, which we call $\preceq _{k}$, by the rule $v\preceq w$ if
either $v$ is an ancestor of $w$, or the branchpoint $b(v,w)$ of $v$ and $w$
is such that the edge leading toward $v$ is earlier than the edge leading
toward $w$ (for the ordering on children of $b(v,w)$). In turn, the
knowledge of $\mathcal{R}(k),\preceq _{k}$, or even of $\mathcal{R}(k)$ and
the restriction of $\preceq _{k}$ to the leaves $L_{1},\ldots ,L_{k}$ of $%
\mathcal{R}(k)$, allows to recover the planar structure of $\mathcal{R}(k)$.
The family of planar trees $(\mathcal{R}(k),\preceq _{k},k\geq 1)$ is said
to be \emph{consistent} if furthermore for every $1\leq j<k$ the planar tree
$(\mathcal{R}(j),\preceq _{j})$ has the same law as the planar subtree of $(%
\mathcal{R}(k),\preceq _{k})$ spanned by $j$ leaves $L_{1}^{1},\ldots
,L_{j}^{k}$ taken independently uniformly at random among the leaves of $%
\mathcal{R}(k)$.

We build such a consistent family out of the consistent family of unordered
trees $(\mathcal{R}(k),k\geq 1)$ as follows. Starting from the tree $%
\mathcal{R}(1)$, which we endow with the trivial order on its only leaf, we
build recursively the total order on $\mathcal{R}(k+1)$ from the order $%
\preceq _{k}$ on $\mathcal{R}(k)$, so that the restriction of $\preceq
_{k+1} $ to the leaves $L_{1},\ldots ,L_{k}$ of $\mathcal{R}(k)$ equals $%
\preceq _{k}$. Given $\mathcal{R}(k+1),\preceq _{k},$ let $b(L_{k+1})$ be
the father of $L_{k+1}$. We distinguish two cases:

\begin{enumerate}
\item  if $b(L_{k+1})$ is a vertex of $\mathcal{R}(k)$, which has $r$
children ${c_{1},c_{2},\ldots ,c_{r}}$ in $\mathcal{R}(k)$, choose $J$
uniformly in $\{1,2,\ldots ,r+1\}$ and let $c_{J-1}\preceq
_{k+1}L_{k+1}\preceq _{k+1}c_{J}$, that is, turn $L_{k+1}$ into the $j$-th
son of $b(L_{k+1})$ in $\mathcal{R}(k+1)$ with probability $1/(r+1)$ (here $%
c_{0}$ (resp.\ $c_{r+1}$) is the predecessor (resp.\ successor) of $c_{1}$
(resp.\ $c_{r}$) for $\preceq _{k}$; we simply ignore them if they do not
exist)

\item  else, $b(L_{k+1})$ must have a unique son $s$ besides $L_{k+1}$. Let $%
s^{\prime }$ be the predecessor of $s$ for $\preceq _{k}$ and $s^{\prime
\prime }$ its successor (if any), and we let $s^{\prime }\preceq
_{k+1}L_{k+1}\preceq _{k+1}s$ with probability $1/2$ and $s\preceq
_{k+1}L_{k+1}\preceq _{k+1}s^{\prime \prime }$ with probability $1/2$.
\end{enumerate}

It is easy to see that this procedure uniquely determines the law of the
total order $\preceq_{k+1}$ on $\mathcal{R}(k+1)$ given $\mathcal{R}%
(k+1),\preceq_k$, and hence the law of $(\mathcal{R}(k),\preceq_k,k\geq 1)$
(the important thing being that the order is total).

\begin{lmm}
\label{orderrk} The family of planar trees $(\mathcal{R}(k),\preceq
_{k},k\geq 1)$ is consistent. Moreover, given $\mathcal{R}(k)$, the law of $%
\preceq _{k}$ can be obtained as follows: for each vertex $v$ of $\mathcal{R}%
(k)$, endow the (possibly empty) set $\{c_{1}(v),\ldots ,c_{i}(v)\}$ of
children of $v$ in uniform random order, this independently over different
vertices.
\end{lmm}

{\noindent \textbf{Proof. }} The second statement is obvious by induction.
The first statement follows, since we already know that the family of
unordered trees $(\mathcal{R}(k),k\geq 1)$ is consistent. $\hfill \square$

As a consequence, there exists a.s.\ a unique total order $\preceq$ on the
set of leaves $\{L_1,L_2\ldots\}$ such that the restriction $%
\preceq_{|[k]}=\preceq_k$. One can check that this order extends to a total
order on the set $\mathcal{L}( \mathcal{T }_F)$ : if $L,L^{\prime}$ are
distinct leaves, we say that $L\preceq L^{\prime}$ if and only if there
exist two sequences $L_{\phi(k)}\preceq L_{\varphi(k)},k\geq 1$, the first
one decreasing and converging to $L$ and the second increasing and
converging to $L^{\prime}$. In turn, this extends to a total order (which we
still call $\preceq$) on the whole tree $\mathcal{T }_F$. Theorem \ref{T3}
is now a direct application of \cite[Theorem 15 (iii)]{aldouscrt93}, the
only thing to check being the conditions a) and b) therein (since we already
know that $\mathcal{T }_F$ is compact). Precisely, condition (iii) a)
rewritten to fit our setting spells:
\begin{equation*}
\lim_{k\to\infty} P(\exists 2\leq j\leq k:|D_{\{1,j\}}-aD_1|\leq \delta\, %
\mbox{ and }\,D_j-D_{\{1,j\}}<\delta \, \mbox{ and }\, L_j\preceq L_1)=1.
\end{equation*}
This is thus a slight modification of (\ref{sansordre}), and the proof goes
similarly, the difference being that we need to keep track of the order on
the leaves. Precisely, consider again some rational $r<aD_1$ close to $aD_1$%
, so that $|\Pi_1(r)|\neq0$. The proof of (\ref{sansordre}) shows that
within the time-interval $[r,r+\delta]$, infinitely many integers of $%
\Pi_1(r)$ have been isolated into singletons. Now, by definition of $\preceq$%
, the probability that any of these integers $j$ satisfies $L_j\preceq_j L_1$
is $1/2$. Therefore, infinitely many integers of $\Pi_1(r)$ give birth to a
leaf $L_j$ that satisfy the required conditions, a.s. The proof of
\cite[Condition (iii) b)]{aldouscrt93} is exactly similar, hence proving
Theorem \ref{T3}.

It is worth recalling the detailed construction of the process $H_{F}$,
which is taken from the proof of \cite[Theorem 15]{aldouscrt93} with a
slight modification (we use the leaves $L_{i}$ rather than a new sample $%
Z_{i},i\geq 1$, but one checks that the proof remains valid). Given the
continuum ordered tree $(\mathcal{T}_{F},\mu _{F},\preceq ,(L_{i},i\geq 1))$%
,
\begin{equation*}
U_{i}=\lim_{n\rightarrow \infty }\frac{\#\{j\leq n:L_{j}\preceq L_{i}\}}{n},
\end{equation*}
a limit that exists a.s. Then the family $(U_{i},i\geq 1)$ is distributed as
a sequence of independent sequence of uniformly distributed random variables
on $(0,1)$, and since $\preceq $ is a total order, one has $U_{i}\leq U_{j}$
if and only if $L_{i}\preceq L_{j}$. Next, define $H_{F}(U_{i})$ to be the
height of $L_{i}$ in $\mathcal{T}_{F}$, and extend it by continuity on $%
[0,1] $ (which is a.s.\ possible according to \cite[Theorem 15]{aldouscrt93}%
) to obtain $H_{F}$. In fact, one can define $\tilde{H}_{F}(U_{i})=L_{i}$
and extend it by continuity on $\mathcal{T}_{F}$, in which case $\tilde{H}%
_{F}$ is an isometry between $\mathcal{T}_{F}$ and $((0,1),\overline{d}%
)/\equiv $ that maps (the equivalence class of) $U_{i}$ to $L_{i}$ for $%
i\geq 1$, and which preserves order.

Writing $I_F(t)=\{s\in(0,1):H_F(s)>t\}$, and $|I_F(t)|$ for the decreasing
sequence of the lengths of the interval components of $I_F(t)$, we know from
the above that $(|I_F(t)|,t\geq 0)$ has the same law as $F$. More precisely,

\begin{lmm}
The processes $(|I_{F}(t)|,t\geq 0)$ and $(F(t),t\geq 0)$ are equal.
\end{lmm}

{\noindent \textbf{Proof. }} Let $\Pi ^{\prime }(t)$ be the partition of $%
\mathbb{N}$ such that $i\sim ^{\Pi ^{\prime }(t)}j$ is and only if $U_{i}$
and $U_{j}$ fall in the same interval component of $I_{F}(t).$ The isometry $%
\tilde{H}_{F}$ allows to assimilate $L_{i}$ to $U_{i}$, then the interval
component of $I_{F}(t)$ containing $U_{i}$ corresponds to the tree component
of $\{v\in \mathcal{T}_{F}:\mathrm{ht}(v)>t\}$ containing $L_{i}$, therefore
$U_{j}$ falls in this interval if and only if $i\sim ^{\Pi (t)}j$, and $\Pi
^{\prime }(t)=\Pi (t)$.
By the law of large numbers and the fact that $(U_{j},j\geq 1)$ is
distributed as a uniform i.i.d.\ sample, it follows that the length of the
interval equals the asymptotic frequency of the block of $\Pi (t)$
containing $i$, a.s. for every $t$.
One inverts the assertions ``a.s.'' and ``for every $t$'' by a simple
monotony argument, showing that if $(U_{i},i\geq 1)$ is a uniform i.i.d.\
sample, then a.s.\ for every sub-interval $(a,b)$ of $(0,1)$, the asymptotic
frequency $\lim_{n\rightarrow \infty }n^{-1}\#\{i\leq n:U_{i}\in (a,b)\}=b-a$
(use distribution functions). $\hfill \square $

We will also need the following result, which is slightly more accurate than
just saying, as in the introduction, that $(I_{F}(t),t\geq 0)$ is an
``interval representation'' of $F$:

\begin{lmm}
\label{isfrag} The process $(I_{F}(t),t\geq 0)$ is a self-similar interval
fragmentation, meaning that it is nested ($I_{F}(t^{\prime })\subseteq
I_{F}(t)$ for every $0\leq t\leq t^{\prime }$), continuous in probability,
and for every $t,t^{\prime }\geq 0$, given $I_{F}(t)=\bigcup_{i\geq 1}I_{i}$
where $I_{i}$ are pairwise disjoint intervals, $I_{F}(t+t^{\prime })$ has
the same law as $\bigcup_{i\geq 1}g_{i}(I_{F}^{(i)}(t^{\prime
}|I_{i}|^{\alpha }))$, where the $I_{F}^{(i)},i\geq 1$ are independent
copies of $I_{F}$, and $g_{i}$ is the orientation-preserving affine function
that maps $(0,1)$ to $I_{i}$.
\end{lmm}

Here, the ``continuity in probability'' is with respect to the Hausdorff
metric $D$ on compact subsets of $[0,1]$, and it just means that $%
P(D(I_{F}^{c}(t_{n}),I_{F}^{c}(t))>\varepsilon )\rightarrow 0$ as $%
n\rightarrow \infty $ for any sequence $t_{n}\rightarrow t$ and $\varepsilon
>0$ (here $A^{c}=[0,1]\setminus A$).

{\noindent \textbf{Proof. }} The fact that $I_{F}(t)$ is nested is trivial.
Now recall that the different interval components of $I_{F}(t)$ encode the
tree components of $\{v\in \mathcal{T}_{F}:\mathrm{ht}(v)>t\}$, call them $%
\mathcal{T}_{1}(t),\mathcal{T}_{2}(t),\ldots $. We already know that these
trees are rescaled independent copies of $\mathcal{T}_{F}$, that is, they
have the same law as $\mu _{F}(\mathcal{T}_{i}(t))^{-\alpha }\otimes
\mathcal{T}^{(i)},i\geq 1$, where $\mathcal{T}^{(i)},i\geq 1$ are
independent copies of $\mathcal{T}_{F}$. So let $\mathcal{T}^{(i)}=\mu _{F}(%
\mathcal{T}_{i}(t))^{\alpha }\otimes \mathcal{T}_{i}(t)$. Now, the orders
induced by $\preceq $ on the different $\mathcal{T}^{(i)}$'s have the same
law as $\preceq $ and are independent, because they only depend on the $%
L_{j} $'s that fall in each of them. Therefore, the trees $(\mathcal{T}%
^{(i)},\mu ^{(i)},\preceq ^{(i)})$ are independent copies of $(\mathcal{T}%
_{F},\mu _{F},\preceq )$, where $\mu ^{(i)}(\cdot )=\mu _{F}((\mu _{F}(%
\mathcal{T}_{i}(t))^{-\alpha }\otimes \cdot )\cap \mathcal{T}_{i}(t))/\mu
_{F}(\mathcal{T}_{i}(t))$ and $\preceq ^{(i)}$ is the order on $\mathcal{T}%
^{(i)}$ induced by the restriction of $\preceq $ to $\mathcal{T}_{i}(t)$. It
follows by our previous considerations that their respective height
processes $H^{(i)}$ are independent copies of $H_{F}$, and it is easy to
check that given $I_{F}(t)=\bigcup_{i\geq 1}I_{i}$ (where $I_{i}$ is the
interval corresponding to $\mathcal{T}_{i}(t)$), the excursions of $H_{F}$
above $t$ are precisely the processes $\mu (\mathcal{T}_{i}(t))^{-\alpha
}H^{(i)}=|I_{i}|^{-\alpha }H^{(i)}$. The self-similar fragmentation property
follows at once, as the fact that $I_{F}$ is Markov. Thanks to these
properties, we may just check the continuity in probability at time $0$, and
it is trivial because $H_{F}$ is a.s.\ continuous and positive on $(0,1)$. $%
\hfill \square $

It appears that besides these elementary properties, the process $H_{F}$ is
quite hard to study. In order to move one step further, we will try to give
a ``Poissonian construction'' of $H_{F}$, in the same way as we used
properties of the Poisson process construction of $\Pi ^{0}$ to study $%
\mathcal{T }_{F}$. To begin with, we move ``back to the homogeneous case''
by time-changing. For every $x\in \left( 0,1\right) $, let $I_{x}(t)$ be the
interval component of $I_{F}(t)$ containing $x$, and $|I_{x}(t)|$ be its
length ($=0$ if $I_{x}(t)=\varnothing $). Then set
\begin{equation*}
T_{t}^{-1}(x)=\inf \left\{ u\geq 0:\int_{0}^{u}|I_{x}(r)|^{\alpha }\mathrm{d}
r>t\right\} ,
\end{equation*}
and let $I_{F}^{0}(t)$ be the open set constituted of the union of the
intervals $I_{x}(T_{t}^{-1}(x)),x\in (0,1)$ (it suffices in fact to take the
union of the $I_{U_{i}}(T_{t}^{-1}(U_{i})),i\geq 1$). From \cite{bertsfrag02}
and Lemma \ref{isfrag}, $(I_F^0(t),t\geq 0)$ is a self-similar homogeneous
interval fragmentation.

%

\subsection{A Poissonian construction}


Recall that the process $(\Pi (t),t\geq 0)$ is constructed out of a
homogeneous fragmentation $(\Pi ^{0}(t),t\geq 0)$, which has been
appropriately time-changed, and where $(\Pi ^{0}(t),t\geq 0)$ has itself
been constructed out of a Poisson point process $(\Delta _{t},k_{t},t\geq 0)$
with intensity $\kappa _{\nu }\otimes \#$. Further, we mark this Poisson
process by considering, for each jump time $t$ of this Poisson process, a
sequence $(U_{i}(t),i\geq 1)$ of i.i.d.\ random variables that are uniform
on $\left( 0,1\right) $, so that these sequences are independent over
different such $t$'s. We are going to use the marks to build an order on the
non-void blocks of $\Pi ^{0}$. It is convenient first to formalize what we
precisely call an \emph{order} on a set $A$: it is a subset $\mathcal{O}$ of
$A\times A$ satisfying:

\begin{enumerate}
\item  $(i,i)\in \mathcal{O}$ for every $i\in A$

\item  $(i,j)\in \mathcal{O}$ and $(j,i)\in \mathcal{O}$ imply $i=j$

\item  $(i,j)\in \mathcal{O}$ and $(j,k)\in \mathcal{O}$ imply $(i,k)\in
\mathcal{O}$.
\end{enumerate}

If $B\subseteq A$, the restriction to $B$ of the order $\mathcal{O}$ is $
\mathcal{O}_{|B}=\mathcal{O}\cap (B\times B)$.
We now construct a process $(\mathcal{O}(t),t\geq 0)$, with values in the
set of orders of $\mathbb{N}$, as follows. Let $\mathcal{O}(0)=\{(i,i),i\in
\mathbb{N}\}$ be the trivial order, and let $n\in \mathbb{N}$. Let 
$0<t_{1}<t_{2}<\ldots <t_{K}$ be the times of occurrence of jumps of the
Poisson process $(\Delta _{t},k_{t},t\geq 0)$ such that both $k_{t}\leq n$
and $(\Delta _{t})_{|[n]}$ (the restriction of $\Delta_t$ to $[n]$)
is non-trivial. Let $\mathcal{O}^{n}(0)=\mathcal{O
}_{|[n]}(0)$, and define a process $\mathcal{O}^{n}(t)$ to be constant on
the time-intervals $[t_{i-1},t_{i})$ (where $t_{0}=0$), where inductively,
given $\mathcal{O}^{n}(t_{i-1})=\mathcal{O}^{n}(t_{i}-)$, $\mathcal{O}
^{n}(t_{i})$ is defined as follows. Let $J_{n}(t_{i})=\{j\in \Pi
_{k_{t_{i}}}^{0}(t_{i}-):j\leq n\mbox{ and }\Pi _{j}^{0}(t_{i})\neq
\varnothing \}$ so that $k_{t_{i}}\in J_{n}(t_{i})$ as soon as $\Pi
_{k_{t_{i}}}^{0}(t_{i}-)\neq \varnothing $. Let then
\begin{equation*}
\mathcal{O}^{n}(t_{i})=\bigcup_{\substack{ j,k\in J_{n}(t_{i}):  \\ 
U_{j}(t_{i})<U_{k}(t_{i})}}\{(j,k)\}\cup \bigcup_{\substack{ 
j:(j,k_{t_{i}})\in \mathcal{O}_{|[n]}(t_{i}-)  \\ k\in J_{n}(t_{i})}}
\{(j,k)\}\cup \bigcup_{\substack{ j:(k_{t_{i}},j)\in \mathcal{O}
_{|[n]}(t_{i}-)  \\ k\in J_{n}(t_{i})}}\{(k,j)\}.
\end{equation*}
In words, we order each set of new blocks in random order in accordance with
the variables $U_{m}(t_{i}),1\leq m\leq n$, and these new blocks have the
same relative position with other blocks as had their father, namely the
block $\Pi _{k_{t_{i}}}^{0}(t_{i}-)$.

It is not difficult that the orders thus defined are consistent as $n$
varies, i.e.\ $(\mathcal{O}^{n+1}(t))_{|[n]}=\mathcal{O}^{n}(t)$ for every $%
n,t$, and it easily follows that there exists a unique process $(\mathcal{O}%
(t),t\geq 0)$ such that $\mathcal{O}_{|[n]}(t)=\mathcal{O}^{n}(t)$ for every
$n,t$ (for existence, take the union over $n\in \mathbb{N}$, and unicity is
trivial). The process $\mathcal{O}$ thus obtained allows to build an
interval-valued version of the fragmentation $\Pi ^{0}(t)$, namely, for
every $t\geq 0$ and $j\geq 0$ let
\begin{equation*}
I_{j}^{0}(t)=\left( \sum_{k\neq j:(k,j)\in \mathcal{O}(t)}|\Pi
_{k}^{0}(t)|,\sum_{k:(k,j)\in \mathcal{O}(t)}|\Pi _{k}^{0}(t)|\right)
\end{equation*}
(notice that $I_{j}^{0}(t)=\varnothing $ if $\Pi _{j}^{0}(t)=\varnothing $).
Write $I^{0}(t)=\bigcup_{j\geq 1}I_{j}^{0}(t)$, and notice that the length $%
|I_{j}^{0}(t)|$ of $I_{j}^{0}(t)$ equals the asymptotic frequency of $\Pi
_{j}^{0}(t)$ for every $j\geq 1,t\geq 0$.


\begin{prp}
\label{tedious} The processes $(I_{F}^{0}(t),t\geq 0)$ and $(I^{0}(t),t\geq
0)$ have the same law.
\end{prp}

As a consequence, we have obtained a construction of an object with the same
law as $I^0_F$ with the help of a marked Poisson process in $\mathcal{P}%
_{\infty}$, and this is the one we are going to work with.

{\noindent \textbf{Proof. }} Let $I^0_F(i,t)$ be the interval component of $%
I^0_F(t)$ containing $U_i$ if $i$ is the least $j$ such that $U_j$ falls in
this component, and $I_F^0(i,t)=\varnothing$ else. Let $\mathcal{O}%
_F(t)=\{(i,i),i\in \mathbb{N}\}\cup\{(j,k):I^0_F(j,t)
\mbox{ is located to
the left of }I^0_F(k,t)\mbox
{ and both are nonempty}\}$. Since the lengths of the interval components of
$I^0_F$ and $I^0$ are the same, the only thing we need to check is that the
processes $\mathcal{O}$ and $\mathcal{O}_F$ have the same law. But then, for
$j\neq k$, $(j,k)\in\mathcal{O}_F(t)$ means that the branchpoint $b(L_j,L_k)$
of $L_j$ and $L_k$ has height less than $t$, and the subtree rooted at $%
b(L_j,L_k)$ containing $L_j$ has been placed before that containing $L_k$.
Using Lemma \ref{orderrk}, we see that given $\mathcal{T }_F,L_1,L_2,\ldots$%
, the subtrees rooted at any branchpoint $b$ of $\mathcal{T }_F$ are placed
in exchangeable random order independently over branchpoints. Precisely,
letting $\mathcal{T }^b_1$ be the subtree containing the leaf with least
label, $\mathcal{T }^b_2$ the subtree different from $\mathcal{T }^b_1$
containing the leaf with least label, and so on, the first subtrees $%
\mathcal{T }^b_1,\ldots, \mathcal{T }^b_k$ are placed in any of the $k!$
possible linear orders, consistently as $k$ varies. Therefore (see e.g.\
\cite[Lemma 10]{aldouscrt93}), there exist independent uniform$(0,1)$ random
variables $U^b_1,U^b_2,\ldots$ independent over $b$'s such that $\mathcal{T }%
^b_i$ is on the ``left'' of $\mathcal{T }^b_j$ (for the order $\mathcal{O}_F$%
) if and only if $U^b_i\leq U^b_j$. This is exactly how we defined the order
$\mathcal{O}(t)$. $\hfill \square$

\noindent \textbf{Remark. } As the reader may have noticed, this
construction of an interval-valued fragmentation has in fact little to do
with pure manipulation of intervals, and it is actually almost entirely
performed in the world of partitions. We stress that it is in fact quite
hard to construct directly such an interval fragmentation out of the plain
idea: ``start from the interval $(0,1)$, take a Poisson process $%
(s(t),k_{t},t\geq 0)$ with intensity $\nu (\mathrm{d}\mathbf{s})\otimes \#$,
and at a jump time of the Poisson process turn the $k_{t}$-th interval
component $I_{k_{t}}(t-)$ of $I(t-)$ (for some labeling convention) into the
open subset of $I_{k_{t}}(t-)$ whose components sizes are $%
|I_{k_{t}}(t-)|s_{i}(t),i\geq 1$, and placed in exchangeable order''. Using
partitions helps much more than plainly giving a natural ``labeling
convention'' for the intervals. In the same vein, we refer to the work of
Gnedin \cite{gnedin97}, which shows that exchangeable interval (composition)
structures are in fact equivalent to ``exchangeable partitions+order on
blocks''.

For every $x\in(0,1)$, write $I^0_x(t)$ for the interval component of $%
I^0_F(t) $ containing $x$, and notice that $I^0_x(t-)=\bigcap_{s\uparrow
t}I^0_x(s)$ is well-defined as a decreasing intersection. For $t\geq 0$ such
that $I^0_x(t)\neq I^0_x(t-)$, let $s^x(t)$ be the sequence $|I^0_F(t)\cap
I^0_x(t-)|/|I^0_x(t-)|$, where $|I^0_F(t)\cap I^0_x(t-)|$ is the decreasing
sequence of lengths of the interval components of $I^0_F(t)\cap I^0_x(t-)$.
The useful result on the Poissonian construction is given in the following

\begin{lmm}
\label{poisson} The process $(s^{x}(t),t\geq 0)$ is a Poisson point process
with intensity $\nu (\mathrm{d}\mathbf{s})$, and more precisely, the order
of the interval components of $I_{F}^{0}(t)\cap I_{x}^{0}(t-)$ is
exchangeable: there exists a sequence of i.i.d.\ uniform random variables $%
(U_{i}^{x}(t),i\geq 1)$, independent of $(\mathcal{G}^{0}(t-),s^{x}(t))$
such that the interval with length $s_{i}^{x}(t)|I_{x}^{0}(t-)|$ is located
on the left of the interval with length $s_{j}^{x}(t)|I_{x}^{0}(t-)|$ if and
only if $U_{i}^{x}(t)\leq U_{j}^{x}(t)$.
\end{lmm}

{\noindent \textbf{Proof. }} Let $i(t,x)=\inf \{i\in \mathbb{N}:U_{i}\in
I_x^{0}(t)\}$. Then $i(t,x)$ is an increasing jump-hold process in $\mathbb{N%
}$. If now $I_{x}^{0}(t)\neq I_{x}^{0}(t-)$, it means that there has been a
jump of the Poisson process $\Delta _{t},k_{t}$ at time $t$, so that $%
k_{t}=i(t,x)$, and then $s^{x}(t)$ is equal to the decreasing sequence $%
|\Delta _{t}|$ of asymptotic frequencies of $\Delta _{t}$, therefore $%
s^{x}(t)=|\Delta _{t}|$ when $k_{t}=i(t-,x)$, and since $i(t-,x)$ is
progressive, its jump times are stopping times so the process $%
(s^{x}(t),t\geq 0)$ is in turn a Poisson process with intensity $\nu (%
\mathrm{d}\mathbf{s})$. Moreover, by Proposition \ref{tedious} and the
construction of $I^0$, each time an interval splits, the corresponding
blocks are put in exchangeable order, which gives the second half of the
lemma. $\hfill \square$

\subsection{Proof of Theorem \ref{T4}}

\subsubsection{H\"{o}lder-continuity of $H_{F}$}

We prove here that the height process is a.s.\ H\"{o}lder-continuous of
order $\gamma $ for every $\gamma <\vartheta _{\mathrm{low}}\wedge |\alpha
|. $ The proof will proceed in three steps.

\noindent \textbf{First step: Reduction to the behavior of }$H_{F}$\textbf{\
near }$0.$ By a theorem of Garsia Rodemich and Rumsey (see e.g. \cite
{Del-Meyer}), the finiteness of $\int_{0}^{1}\int_{0}^{1}\frac{\left|
H_{F}(x)-H_{F}(y)\right| ^{n+n_{0}}}{\left| x-y\right| ^{\gamma n}}$d$x$d$y$
leads to the $\left( \frac{\gamma n-2}{n+n_{0}}\right) $-H%
\"{o}lder-continuity of $H_{F},$ so that when the previous integral is
finite for every $n,$ the height process $H_{F}$ is H\"{o}lder-continuous of
order $\delta $ for every $\delta <\gamma ,$ whatever is $n_{0}.$ To prove
Theorem \ref{T4} it is thus sufficient to show that for every $\gamma
<\vartheta _{\mathrm{low}}\wedge \left| \alpha \right| $ there exists a $%
n_{0}(\gamma )$ such that
\begin{equation*}
E\left[ \int_{0}^{1}\int_{0}^{1}\frac{\left| H_{F}(x)-H_{F}(y)\right|
^{n+n_{0}\left( \gamma \right) }}{\left| x-y\right| ^{\gamma n}}\text{d}x%
\text{d}y\right] <\infty \text{ for every integer }n.
\end{equation*}
Now take $V_{1},V_{2}$ uniform independent on $\left( 0,1\right) $,
independently of $H_{F}$.
The expectation above then rewrites $E\left[ \frac{\left|
H_{F}(V_{1})-H_{F}(V_{2})\right| ^{n+n_{0}\left( \gamma \right) }}{\left|
V_{1}-V_{2}\right| ^{\gamma n}}\right] .$

Consider next $I_{F}$ the interval fragmentation constructed from $H_{F}$
(see Section 4.1). By Lemma \ref{isfrag}, $H_{F}(V_{1})$ and $H_{F}(V_{2})$
may be rewritten as

\begin{equation*}
H_{F}(V_{i})=D_{\left\{ 1,2\right\} }+\lambda _{i}^{\left| \alpha \right|
}(D_{\left\{ 1,2\right\} })\widetilde{D}_{i},\text{ }i=1,2,
\end{equation*}
where $D_{\left\{ 1,2\right\} }$ is the first time at which $V_{1}$ and $%
V_{2}$ belong to different intervals of $I_{F}$ and $\widetilde{D}_{1},%
\widetilde{D}_{2}$ have the same law as $H_{F}(V_{1})$ and are independent
of $\mathcal{H}(D_{\left\{ 1,2\right\} })$, where $\mathcal{H}(t),t\geq 0$
is the natural completed filtration associated to $I_{F}$. The r.v. $%
\widetilde{D}_{1}$ and $\widetilde{D}_{2}$ can actually be described more
precisely. Say that at time $D_{\left\{ 1,2\right\} },$ $V_{1}$ belongs to
an interval $\left( a_{1},a_{1}+\lambda _{1}(D_{\left\{ 1,2\right\}
})\right) $ and $V_{2}$ to $\left( a_{2},a_{2}+\lambda _{2}(D_{\left\{
1,2\right\} })\right) .$ Then there exist two iid processes independent of $%
\mathcal{H}(D_{\left\{ 1,2\right\} })$ and with the same law as $H_{F},$ let
us denote them $H_{F}^{(1)}$ and $H_{F}^{(2)},$ such that $\widetilde{D}%
_{i}=H_{F}^{(i)}\left( \frac{V_{i}-a_{i}}{\lambda _{i}(D_{\left\{
1,2\right\} })}\right) ,$ $i=1,2$. Since $V_{i}\in \left(
a_{i},a_{i}+\lambda _{i}(D_{\left\{ 1,2\right\} })\right) $, the random
variables $\widetilde{V}_{i}=\left( V_{i}-a_{1}\right) \lambda
_{i}^{-1}(D_{\left\{ 1,2\right\} })$ are iid, with the uniform law on $%
\left( 0,1\right) $ and independent of $H_{F}^{(1)},H_{F}^{(2)}$ and $%
\mathcal{H}(D_{\left\{ 1,2\right\} }).$ And when $V_{1}>V_{2},$%
\begin{equation*}
V_{1}-V_{2}\geq \lambda _{1}(D_{\left\{ 1,2\right\} })\widetilde{V}%
_{1}+\lambda _{2}(D_{\left\{ 1,2\right\} })\left( 1-\widetilde{V}_{2}\right)
\end{equation*}
since $a_{1}$ is then largest than $a_{2}+\lambda _{2}(D_{\left\{
1,2\right\} }).$ This gives
\begin{eqnarray*}
E\left[ \frac{\left| D_{1}-D_{2}\right| ^{n+n_{0}\left( \gamma \right) }}{%
\left| V_{1}-V_{2}\right| ^{\gamma n}}\right] &=&2E\left[ \frac{\left|
D_{1}-D_{2}\right| ^{n+n_{0}\left( \gamma \right) }}{\left(
V_{1}-V_{2}\right) ^{\gamma n}}1_{\left\{ V_{1}>V_{2}\right\} }\right] \\
&\leq &2E\left[ \frac{\left( \lambda _{1}^{\left| \alpha \right|
}(D_{\left\{ 1,2\right\} })\widetilde{D}_{1}+\lambda _{2}^{\left| \alpha
\right| }(D_{\left\{ 1,2\right\} })\widetilde{D}_{2}\right) ^{n+n_{0}\left(
\gamma \right) }}{\left( \lambda _{1}(D_{\left\{ 1,2\right\} })\widetilde{V}%
_{1}+\lambda _{2}(D_{\left\{ 1,2\right\} })\left( 1-\widetilde{V}_{2}\right)
\right) ^{\gamma n}}\right]
\end{eqnarray*}
and this last expectation is bounded from above by
\begin{equation*}
2^{n+n_{0}\left( \gamma \right) }E\left[ \left( \lambda _{1}(D_{\left\{
1,2\right\} })\right) ^{\left( n+n_{0}\left( \gamma \right) \right) \left|
\alpha \right| -\gamma n}\right] \left( E\left[ \frac{H_{F}^{n+n_{0}\left(
\gamma \right) }(V_{1})}{V_{1}^{\gamma n}}\right] +E\left[ \frac{%
H_{F}^{n+n_{0}\left( \gamma \right) }(V_{1})}{\left( 1-V_{1}\right) ^{\gamma
n}}\right] \right) .
\end{equation*}
The expectation involving $\lambda _{1}$ is bounded by $1$ since $\gamma
<\left| \alpha \right| .$ And since $V_{1}$ is independent of $H_{F},$ the
two expectations in the parenthesis are equal (reversing the order $\preceq $
and performing the construction of $H_{F}$ gives a process with the same law
and shows that $H_{F}(x)\overset{law}{=}H_{F}(1-x)$ for every $x\in \left(
0,1\right) $) and finite as soon as
\begin{equation}
\sup_{x\in \left( 0,1\right) }E\left[ H_{F}(x)^{n+n_{0}(\gamma )}\right]
x^{-\gamma n}<\infty .  \label{202}
\end{equation}
The rest of the proof thus consists in finding an integer $n_{0}(\gamma )$
such that $\left( \ref{202}\right) $ holds for every $n.$ To do so, we will
have to observe the interval fragmentation $I_{F}$ at nice stopping times
depending on $x,$ say $\mathbb{T}_{x}^{(\gamma )},$ and then use the strong
fragmentation property at time $\mathbb{T}_{x}^{(\gamma )}.$ This gives
\begin{equation}
H_{F}(x)=\mathbb{T}_{x}^{(\gamma )}+\left( S_{x}(\mathbb{T}_{x}^{(\gamma
)})\right) ^{\left| \alpha \right| }\overline{H}_{F}(P_{x}(\mathbb{T}%
_{x}^{(\gamma )}))  \label{200}
\end{equation}
where $S_{x}(\mathbb{T}_{x}^{(\gamma )})$ is the length of the interval
containing $x$ at time $\mathbb{T}_{x}^{(\gamma )},$ $P_{x}(\mathbb{T}%
_{x}^{(\gamma )})$ the relative position of $x$ in that interval and $%
\overline{H}_{F}$ a process with the same law as $H_{F}$ and independent of $%
\mathcal{H}(\mathbb{T}_{x}^{(\gamma )})$. 

\noindent \textbf{Second step: Choice and properties of }$\mathbb{T}%
_{x}^{(\gamma )}$\textbf{. }Let us first introduce some notation in order to
prove the forthcoming Lemma $\ref{lemmeTxSx}.$ Recall that we have called $%
I_{F}^{0}$ the homogeneous interval fragmentation related to $I_{F}$ by the
time changes $T_{t}^{-1}(x)$ introduced in Section 4.1. In this homogeneous
fragmentation, let

$I_{x}^{0}(t)=\left( a_{x}(t),b_{x}(t)\right) $ be the interval containing $%
x $ at time $t$

$S^0_{x}(t)$ the length of this interval

$P^0_{x}(t)=(x-a_x(t))/S^0_x(t)$ the relative position of $x$ in $I_{x}(t).$

\noindent Similarly, we define $P_{x}^{0}(t-)$ to be the relative position
of $x$ in the interval $I_{x}^{0}(t-),$ which is well-defined as an
intersection of nested intervals. $S_{x}^{0}(t-)$ is the size of this
interval. We will need the following inequalities in the sequel:
\begin{equation*}
\begin{array}{cc}
P_{x}^{0}(t)\leq x/S_{x}^{0}(t) & \text{ \ \ \ }P_{x}^{0}(t-)\leq
x/S_{x}^{0}(t-).
\end{array}
\end{equation*}

Next recall the Poisson point process construction of the interval
fragmentation $I_{F}^{0}$, and the Poisson point process $\left(
s^{x}(t)\right) _{t\geq 0}$ of Lemma \ref{poisson}. Set
\begin{equation*}
\sigma (t):=-\ln \left( \prod_{s\leq t}s_{1}^{x}(t)\right) \,\ \ \ \ \ \
t\geq 0,
\end{equation*}
with the convention $s_{1}^{x}(t)=1$ when $t$ is not a time of occurrence of
the point process. By Lemma \ref{poisson}, the process $\sigma $ is a
subordinator with intensity measure $\nu (-\ln s_{1}\in x),$ which is
infinite. Consider then $T_{x}^{\mathrm{exit}},$ the first time at which $x$
is not in the largest sub-interval of $I_{x}^{0}$ when $I_{x}^{0}$ splits,
that is
\begin{equation*}
T_{x}^{\mathrm{exit}}:=\inf \left\{ t:S_{x}^{0}(t)<\exp (-\sigma
(t))\right\} .
\end{equation*}
By definition, the size of the interval containing $x$ at time $t<T_{x}^{%
\mathrm{exit}}$ is given by $S_{x}^{0}(t)=\exp (-\sigma (t)).$ We will need
to consider the first time at which this size is smaller than $a,$ for $a$
in $\left( 0,1\right) ,$ and so we introduce
\begin{equation*}
T_{a}^{\sigma }:=\inf \left\{ t:\exp (-\sigma (t))<a\right\} .
\end{equation*}
Note that $P_{x}^{0}(t)\leq x\exp (\sigma (t))$ when $t<T_{x}^{\mathrm{exit}%
} $ and that $P_{x}^{0}(T_{x}^{\mathrm{exit}}-)\leq x\exp (\sigma (T_{x}^{%
\mathrm{exit}}-)).$

Finally, to obtain a nice $\mathbb{T}_{x}^{\left( \gamma \right) }$ as
required in the preceding step, we stop the homogeneous fragmentation at
time
\begin{equation*}
T_{x}^{\mathrm{exit}}\wedge T_{x^{\varepsilon }}^{\sigma }
\end{equation*}
for some $\varepsilon $ to be determined (and depending on $\gamma $) and
then take for $\mathbb{T}_{x}^{\left( \gamma \right) }$ the self-similar
counterpart of this stopping time, that is $\mathbb{T}_{x}^{\left( \gamma
\right) }=T_{T_{x}^{\mathrm{exit}}\wedge T_{x^{\varepsilon }}^{\sigma
}}^{-1}(x).$ More precisely, we have

\begin{lmm}
\label{lemmeTxSx} For every $\gamma <\vartheta _{\mathrm{low}}\wedge \left|
\alpha \right| ,$ there exist a family of random stopping times $\mathbb{T}%
_{x}^{\left( \gamma \right) },x\in \left( 0,1\right) ,$ and an integer $%
N(\gamma )$ such that

(i) for every $n\geq 0,$ $\exists C_{1}(n):E\left[ \left( \mathbb{T}%
_{x}^{\left( \gamma \right) }\right) ^{n}\right] \leq C_{1}(n)x^{\gamma n}$ $%
\ \ \forall $ $x\in \left( 0,1\right) ,$

(ii) $\exists C_{2}$ such that $E\left[ \left( S_{x}(\mathbb{T}_{x}^{\left(
\gamma \right) })\right) ^{n}\right] \leq C_{2}x^{\gamma }$ for every $x$ in
$\left( 0,1\right) $ and $n\geq N(\gamma ).$
\end{lmm}

{\noindent \textbf{Proof. }}Fix $\gamma <\vartheta _{\mathrm{low}}\wedge
\left| \alpha \right| $ and then $\varepsilon <1$ such that $\gamma
/(1-\varepsilon )<\vartheta _{\mathrm{low}}$. The times $\mathbb{T}%
_{x}^{\left( \gamma \right) },x\in \left( 0,1\right) ,$ are constructed from
this $\varepsilon $ by
\begin{equation*}
\mathbb{T}_{x}^{\left( \gamma \right) }=T_{T_{x}^{\mathrm{exit}}\wedge
T_{x^{\varepsilon }}^{\sigma }}^{-1}(x),
\end{equation*}
and it may be clear that these times are stopping times with respect to%
\textbf{\ }$\mathcal{H}.$ A first remark is that the function $x\in \left(
0,1\right) \mapsto S_{x}(\mathbb{T}_{x}^{\left( \gamma \right) })$ is
bounded from above by $1$ and that $x\in \left( 0,1\right) \mapsto \mathbb{T}%
_{x}^{\left( \gamma \right) }$ is bounded from above by $\tau ,$ the first
time at which the fragmentation is entirely reduced to dust, that is, in
others words, the supremum of $H_{F}$ on $\left[ 0,1\right] .$ Since $\tau $
has moments of all orders, it is thus sufficient to prove statements (i) and
(ii) for $x\in \left( 0,x_{0}\right) $ for some well chosen $x_{0}>0.$
Another remark, using the definition of $T_{t}^{-1}(x),$ is that $\mathbb{T}%
_{x}^{\left( \gamma \right) }\leq T_{x}^{\mathrm{exit}}\wedge
T_{x^{\varepsilon }}^{\sigma }$ and $S_{x}(\mathbb{T}_{x}^{\left( \gamma
\right) })=S_{x}^{0}\left( T_{x}^{\mathrm{exit}}\wedge T_{x^{\varepsilon
}}^{\sigma }\right) ,$ so that we just have to prove (i) and (ii) by
replacing in the statement $\mathbb{T}_{x}^{\left( \gamma \right) }$ by $%
T_{x}^{\mathrm{exit}}\wedge T_{x^{\varepsilon }}^{\sigma }$ and $S_{x}(%
\mathbb{T}_{x}^{\left( \gamma \right) })$ by $S_{x}^{0}\left( T_{x}^{\mathrm{%
exit}}\wedge T_{x^{\varepsilon }}^{\sigma }\right) $.

We shall thus work with the homogeneous fragmentation. When $I_{x}^{0}$
splits to give smaller intervals, we divide these sub-intervals into three
groups: the largest sub-interval, the group of sub-intervals on its left and
the the group of sub-intervals on its right. With the notations of Lemma \ref
{poisson}, the lengths of the intervals belonging to the group on the left
are the $s_{i}^{x}(t)S_{x}^{0}(t-)$ with $i$ such that $%
U_{i}^{x}(t)<U_{1}^{x}(t)$ and similarly, the lengths of the intervals on
the right are the $s_{i}^{x}(t)S_{x}^{0}(t-)$ with $i$ such that $%
U_{i}^{x}(t)>U_{1}^{x}(t).$ An important point is that when $T_{x}^{\mathrm{%
exit}}<T_{x^{\varepsilon }}^{\sigma },$ then at time $T_{x}^{\mathrm{exit}},$
the point $x$ belongs to the group of sub-intervals on the left resulting
from the fragmentation of $I_{x}^{0}(T_{x}^{\mathrm{exit}}-).$ Indeed, when $%
T_{x}^{\mathrm{exit}}<T_{x^{\varepsilon }}^{\sigma },$ $\exp (-\sigma
(T_{x}^{\mathrm{exit}}))\geq x,$ which rewrites $s_{1}^{x}(T_{x}^{\mathrm{%
exit}})\exp (-\sigma (T_{x}^{\mathrm{exit}}-))\geq x.$ Using then that $%
P_{x}^{0}(T_{x}^{\mathrm{exit}}-)\leq x\exp (\sigma (T_{x}^{\mathrm{exit}%
}-)),$ we obtain $s_{1}^{x}(T_{x}^{\mathrm{exit}})\geq P_{x}^{0}(T_{x}^{%
\mathrm{exit}}-)$ and thus that $x$ does not belong to the group on the
right at time $T_{x}^{\mathrm{exit}}$ ($x$ belongs to the group on the right
at a time $t$ if and only if $P_{x}^{0}(t-)>\sum_{i:U_{i}^{x}(t)\leq
U_{1}^{x}(t)}s_{i}^{x}(t)$). Hence $x$ belongs to the union of intervals on
the left at time $T_{x}^{\mathrm{exit}}$ when $T_{x}^{\mathrm{exit}%
}<T_{x^{\varepsilon }}^{\sigma }$. In others words,
\begin{equation*}
T_{x}^{\mathrm{exit}}=\inf \left\{
t:\sum\nolimits_{i:U_{i}^{x}(t)<U_{1}^{x}(t)}s_{i}^{x}(t)>P_{x}^{0}(t-)%
\right\} \text{ when }T_{x}^{\mathrm{exit}}<T_{x^{\varepsilon }}^{\sigma }.
\end{equation*}

The key-point, consequence of Lemma \ref{poisson}, is that the process $%
\left( \sum_{i:U_{i}^{x}(t)<U_{1}^{x}(t)}s_{i}^{x}(t)\right) _{t\geq 0}$is a
marked Poisson point process with an intensity measure on $\left[ 0,1\right]
$ given by
\begin{equation*}
\mu (\mathrm{d}u):=\int_{S}p(\mathbf{s},\mathrm{d}u)\nu (\mathrm{d}\mathbf{s}%
),\text{ }u\in \left[ 0,1\right] ,
\end{equation*}
where for a fixed $\mathbf{s}$ in $S,$ $p(\mathbf{s},\mathrm{d}u)$ is the
law of $\sum_{i:U_{i}<U_{1}}s_{i},$ the $U_{i}$'s being uniform and
independent random variables. We refer to Kingman \cite{kingman93} for
details on marked Poisson point processes. Observing then that for any $a$
in $\left( 0,1/2\right) $ and for a fixed $\mathbf{s}$ in $S$
\begin{equation*}
1_{\left\{ 1-s_{1}>2a\right\} }\leq 1_{\left\{
\sum_{i:U_{i}<U_{1}}s_{i}>a\right\} }+1_{\left\{
\sum_{i:U_{i}>U_{1}}s_{i}>a\right\} },
\end{equation*}
we obtain that $1_{\left( 1-s_{1}>2a\right) }\leq 2P\left(
\sum_{i:U_{i}<U_{1}}s_{i}>a\right) $ and then the following inequality
\begin{equation*}
\mu \left( \left( a,1\right] \right) \geq \frac{1}{2}\nu \left(
s_{1}<1-2a\right) .
\end{equation*}
This, recalling the definition of $\vartheta _{\mathrm{low}}$ and that $%
\gamma /(1-\varepsilon )<\vartheta _{\mathrm{low}}$, leads to the existence
of a positive $x_{0}$ and a positive constant $C$ such that
\begin{equation}
\mu \left( \left( x^{1-\varepsilon },1\right] \right) \geq C\left(
x^{-\left( 1-\varepsilon \right) }\right) ^{\gamma /(1-\varepsilon
)}=Cx^{-\gamma }\text{ \ for all }x\text{ in }\left( 0,x_{0}\right) \text{.}
\label{201}
\end{equation}

\noindent \emph{Proof of (i)}. We again have to introduce a hitting time,
that is the first time at which the Poisson point process $\left(
\sum_{i:U_{i}^{x}(t)<U_{1}^{x}(t)}s_{i}^{x}(t)\right) _{t\geq 0}$ belongs to
$\left( x^{1-\varepsilon },1\right) :$
\begin{equation*}
H_{x^{1-\varepsilon }}:=\inf \left\{
t:\sum\nolimits_{i:U_{i}^{x}(t)<U_{1}^{x}(t)}s_{i}^{x}(t)>x^{1-\varepsilon
}\right\} .
\end{equation*}
By the theory of Poisson point processes, this time has an exponential law
with parameter $\mu \left( \left( x^{1-\varepsilon },1\right] \right) .$
Hence, given inequality $\left( \ref{201}\right) ,$ it is sufficient to show
that $T_{x}^{\mathrm{exit}}\wedge T_{x^{\varepsilon }}^{\sigma }\leq
H_{x^{1-\varepsilon }}$ to obtain (i) for $x$ in $\left( 0,x_{0}\right) $
and then (i) (we recall that it is already known that $\sup_{x\in \left[
x_{0},1\right) }x^{-\gamma n}E\left[ \left( \mathbb{T}_{x}^{\left( \gamma
\right) }\right) ^{n}\right] $ is finite). On the one hand, since $%
P_{x}^{0}(t)\leq x\exp (\sigma (t))$ when $t<T_{x}^{\mathrm{exit}},$%
\begin{equation*}
P_{x}^{0}(H_{x^{1-\varepsilon }}-)\leq x\exp (\sigma (H_{x^{1-\varepsilon
}}-))<x\exp (\sigma (H_{x^{1-\varepsilon }}))\text{ when }%
H_{x^{1-\varepsilon }}<T_{x}^{\mathrm{exit}}.
\end{equation*}
On the other hand, $H_{x^{1-\varepsilon }}<T_{x^{\varepsilon }}^{\sigma }$
yields
\begin{equation*}
x\exp (\sigma (H_{x^{1-\varepsilon }}))\leq x^{1-\varepsilon
}<\sum\nolimits_{i:U_{i}^{x}(H_{x^{1-\varepsilon
}})<U_{1}^{x}(H_{x^{1-\varepsilon }})}s_{i}^{x}(H_{x^{1-\varepsilon }}),
\end{equation*}
and combining these two remarks, we get that $H_{x^{1-\varepsilon }}<T_{x}^{%
\mathrm{exit}}\wedge T_{x^{\varepsilon }}^{\sigma }$ implies
\begin{equation*}
P_{x}^{0}(H_{x^{1-\varepsilon
}}-)<\sum\nolimits_{i:U_{i}^{x}(H_{x^{1-\varepsilon
}})<U_{1}^{x}(H_{x^{1-\varepsilon }})}s_{i}^{x}(H_{x^{1-\varepsilon }}).
\end{equation*}
Yet this is not possible, because this last relation on $H_{x^{1-\varepsilon
}}$ means that, at time $H_{x^{1-\varepsilon }}$, $x$ is not in the largest
sub-interval resulting from the splitting of $I_{x}^{0}(H_{x^{1-\varepsilon
}}-),$ which implies $H_{x^{1-\varepsilon }}\geq T_{x}^{\mathrm{exit}}$ and
this does not match with $H_{x^{1-\varepsilon }}<T_{x}^{\mathrm{exit}}\wedge
T_{x^{\varepsilon }}^{\sigma }.$ Hence $T_{x}^{\mathrm{exit}}\wedge
T_{x^{\varepsilon }}^{\sigma }\leq H_{x^{1-\varepsilon }}$ and (i) is proved$%
.$

\noindent \emph{Proof of (ii).} Take $N(\gamma )\geq \gamma /\varepsilon
\vee 1.$ When $T_{x^{\varepsilon }}^{\sigma }\leq T_{x}^{\mathrm{exit}},$
using the definition of $T_{x^{\varepsilon }}^{\sigma }$ and the right
continuity of $\sigma ,$ we have
\begin{equation*}
S_{x}^{0}(T_{x}^{\mathrm{exit}}\wedge T_{x^{\varepsilon }}^{\sigma })\leq
\exp (-\sigma (T_{x^{\varepsilon }}^{\sigma }))\leq x^{\varepsilon }
\end{equation*}
and consequently $\left( S_{x}^{0}(T_{x}^{\mathrm{exit}}\wedge
T_{x^{\varepsilon }}^{\sigma })\right) ^{N(\gamma )}\leq x^{\gamma }.$ Thus
it just remains to show that
\begin{equation*}
E\left[ \left( S_{x}^{0}(T_{x}^{\mathrm{exit}}\wedge T_{x^{\varepsilon
}}^{\sigma })\right) ^{N(\gamma )}1_{\left\{ T_{x}^{\mathrm{exit}%
}<T_{x^{\varepsilon }}^{\sigma }\right\} }\right] \leq x^{\gamma }\text{ for
}x<x_{0}.
\end{equation*}
When $T_{x}^{\mathrm{exit}}<T_{x^{\varepsilon }}^{\sigma },$ we know - as
explained at the beginning of the proof - that $x$ belongs at time $T_{x}^{%
\mathrm{exit}}$ to the group of sub-intervals on the left resulting from the
fragmentation of $I_{x}^{0}(T_{x}^{\mathrm{exit}}-)$ and hence that $%
S_{x}^{0}(T_{x}^{\mathrm{exit}}\wedge T_{x^{\varepsilon }}^{\sigma
})^{N(\gamma )}\leq s_{i}^{x}(T_{x}^{\mathrm{exit}})$ for some $i$ such that
$U_{i}^{x}(T_{x}^{\mathrm{exit}})<U_{1}^{x}(T_{x}^{\mathrm{exit}}).$ More
roughly,
\begin{equation*}
S_{x}^{0}(T_{x}^{\mathrm{exit}}\wedge T_{x^{\varepsilon }}^{\sigma
})^{N(\gamma )}1_{\left\{ T_{x}^{\mathrm{exit}}<T_{x^{\varepsilon }}^{\sigma
}\right\} }\leq \sum\nolimits_{i:U_{i}^{x}(T_{x}^{\mathrm{exit}%
})<U_{1}^{x}(T_{x}^{\mathrm{exit}})}s_{i}^{x}(T_{x}^{\mathrm{exit}%
})1_{\left\{ T_{x}^{\mathrm{exit}}<T_{x^{\varepsilon }}^{\sigma }\right\} }.
\end{equation*}
To evaluate the expectation of this random sum, recall from the proof of (i)
that $T_{x}^{\mathrm{exit}}\leq H_{x^{1-\varepsilon }}$ when $T_{x}^{\mathrm{%
exit}}<T_{x^{\varepsilon }}^{\sigma }$ and remark that either $T_{x}^{%
\mathrm{exit}}<H_{x^{1-\varepsilon }}$ and then
\begin{equation*}
\sum\nolimits_{i:U_{i}^{x}(T_{x}^{\mathrm{exit}})<U_{1}^{x}(T_{x}^{\mathrm{%
exit}})}s_{i}^{x}(T_{x}^{\mathrm{exit}})\leq x^{1-\varepsilon }\leq
x^{\gamma }\text{ \ \ (}\gamma <\vartheta _{\mathrm{low}}(1-\varepsilon
)\leq 1-\varepsilon \text{)}
\end{equation*}
or $T_{x}^{\mathrm{exit}}=H_{x^{1-\varepsilon }}$ and then
\begin{equation*}
\sum\nolimits_{i:U_{i}^{x}(T_{x}^{\mathrm{exit}})<U_{1}^{x}(T_{x}^{\mathrm{%
exit}})}s_{i}^{x}(T_{x}^{\mathrm{exit}})=\sum\nolimits_{i:U_{i}^{x}(H_{x^{1-%
\varepsilon }})<U_{1}^{x}(H_{x^{1-\varepsilon
}})}s_{i}^{x}(H_{x^{1-\varepsilon }}).
\end{equation*}
There we conclude with the following inequality
\begin{eqnarray*}
E\left[ \sum\nolimits_{i:U_{i}^{x}(H_{x^{1-\varepsilon
}})<U_{1}^{x}(H_{x^{1-\varepsilon }})}s_{i}^{x}(H_{x^{1-\varepsilon }})%
\right] &=&\frac{\int_{S}E\left[ \sum_{i:U_{i}<U_{1}}s_{i}1_{\left\{
\sum_{i:U_{i}<U_{1}}s_{i}>x^{1-\varepsilon }\right\} }\right] \nu (\text{d}s)%
}{\mu \left( \left( x^{1-\varepsilon },1\right] \right) } \\
&\leq &C^{-1}x^{\gamma }\int_{S}\left( 1-s_{1}\right) \nu (\text{d}s),\text{
}x\in \left( 0,x_{0}\right) .
\end{eqnarray*}
\cq

\noindent \textbf{Third step: Proof of }$\left( \ref{202}\right) .$ Fix $%
\gamma <\vartheta _{\mathrm{low}}\wedge \left| \alpha \right| $ and take $%
\mathbb{T}_{x}^{\left( \gamma \right) }$ and $N(\gamma )$ as introduced in
Lemma \ref{lemmeTxSx}. Let then $n_{0}(\gamma )$ be an integer larger than $%
N(\gamma )/\left| \alpha \right| $. According to the first step, Theorem \ref
{T4} is proved if $\left( \ref{202}\right) $ holds for this $n_{0}(\gamma )$
and every integer $n\geq 1.$ To show this, it is obviously sufficient to
prove that for every integers $n\geq 1$ and $m\geq 0,$ there exists a finite
constant $C(n,m)$ such that
\begin{equation*}
E\left[ H_{F}(x)^{m+n+n_{0}(\gamma )}\right] \leq C(n,m)x^{\gamma n}\text{ \
}\forall x\in \left( 0,1\right) .
\end{equation*}
This can be proved by induction: for $n=1$ and every $m\geq 0,$ using $%
\left( \ref{200}\right) ,$ we have
\begin{equation*}
\begin{array}{ll}
E\left[ H_{F}(x)^{m+1+n_{0}(\gamma )}\right] & \leq 2^{m+1+n_{0}(\gamma )}
\\
& \times \left( E\left[ \left( \mathbb{T}_{x}^{\left( \gamma \right)
}\right) ^{m+1+n_{0}(\gamma )}\right] +E\left[ \left( S_{x}(\mathbb{T}%
_{x}^{\left( \gamma \right) })\right) ^{\left| \alpha \right| \left(
m+1+n_{0}(\gamma )\right) }\widetilde{\tau }^{m+1+n_{0}(\gamma )}\right]
\right)
\end{array}
\end{equation*}
where $\widetilde{\tau }$ is the maximum of $\overline{H}_{F}$ on $\left(
0,1\right) .$ Recall that this maximum is independent of $S_{x}(\mathbb{T}%
_{x}^{\left( \gamma \right) })$ and has moments of all orders. Since
moreover $\left| \alpha \right| \left( m+1+n_{0}(\gamma )\right) \geq
N(\gamma ),$ we can apply Lemma \ref{lemmeTxSx} to deduce the existence of a
constant $C(1,m)$ such that
\begin{equation*}
E\left[ H_{F}(x)^{m+1+n_{0}(\gamma )}\right] \leq C(1,m)x^{\gamma }\text{ \
for }x\text{ in }\left( 0,1\right) .
\end{equation*}
Now suppose that for some fixed $n$ and every $m\geq 0,$
\begin{equation*}
E\left[ H_{F}(x)^{m+n+n_{0}(\gamma )}\right] \leq C(n,m)x^{\gamma n}\text{ \
}\forall x\in \left( 0,1\right) .
\end{equation*}
Then,
\begin{eqnarray*}
E\left[ \left( \overline{H}_{F}\left( P_{x}(\mathbb{T}_{x}^{\left( \gamma
\right) })\right) \right) ^{m+n+1+n_{0}(\gamma )}\mid \mathcal{H}\left(
\mathbb{T}_{x}^{\left( \gamma \right) }\right) \right] &\leq &C(n,m+1)\left(
P_{x}(\mathbb{T}_{x}^{\left( \gamma \right) })\right) ^{\gamma n} \\
&\leq &C(n,m+1)\left( S_{x}(\mathbb{T}_{x}^{\left( \gamma \right) })\right)
^{-\gamma n}x^{\gamma n}
\end{eqnarray*}
since $P_{x}(\mathbb{T}_{x}^{\left( \gamma \right) })\leq x/S_{x}(\mathbb{T}%
_{x}^{\left( \gamma \right) }).$ Next, by $\left( \ref{200}\right) ,$%
\begin{equation*}
\begin{array}{ll}
E\left[ H_{F}(x)^{m+n+1+n_{0}(\gamma )}\right] & \leq 2^{m+n+1+n_{0}(\gamma
)}E\left[ \left( \mathbb{T}_{x}^{\left( \gamma \right) }\right)
^{m+n+1+n_{0}(\gamma )}\right] \\
& +2^{m+n+1+n_{0}(\gamma )}C(n,m+1)E\left[ \left( S_{x}(\mathbb{T}%
_{x}^{\left( \gamma \right) })\right) ^{\left| \alpha \right| \left(
m+n+1+n_{0}(\gamma )\right) -\gamma n}\right] x^{\gamma n}.
\end{array}
\end{equation*}
Since $\gamma <\left| \alpha \right| ,$ the exponent $\left| \alpha \right|
\left( m+n+1+n_{0}(\gamma )\right) -\gamma n\geq N(\gamma ),$ and hence
Lemma \ref{lemmeTxSx} applies to give, together with the previous
inequality, the existence of a finite constant $C(n+1,m)$ such that
\begin{equation*}
E\left[ H_{F}(x)^{m+n+1+n_{0}(\gamma )}\right] \leq C(n+1,m)x^{\gamma \left(
n+1\right) }
\end{equation*}
for every $x$ in $\left( 0,1\right) .$ This holds for every $m$ and hence
the induction, formula $\left( \ref{202}\right) $ and Theorem \ref{T4} are
proved.

\subsubsection{Maximal H\"{o}lder exponent of the height process}

The aim of this subsection is to prove that a.s. $H_{F}$ cannot be
H\"{o}lder-continuous of order $\gamma $ for any $\gamma >\vartheta _{%
\mathrm{up}}\wedge |\alpha |/\varrho .$

We first\textbf{\ }prove that\textbf{\ }$H_{F}$ cannot be
H\"{o}lder-continuous with an exponent $\gamma $ larger than $\vartheta _{%
\text{up}}.$ To see this, consider the interval fragmentation $I_{F}$ and
let $U$ be a r.v. independent of $I_{F}$ and with the uniform law on $\left(
0,1\right) .$ By Corollary 2 in \cite{bertsfrag02}, there is a subordinator $%
\left( \theta (t),t\geq 0\right) $ with no drift and a L\'{e}vy measure
given by
\begin{equation*}
\pi _{\theta }(dx)=e^{-x}\sum_{i=1}^{\infty }\nu (-\log s_{i}\in dx),x\in
\left( 0,\infty \right) ,
\end{equation*}
such that the length of the interval component of $I_{F}$ containing $U$ at
time $t$ is equal to $\exp (-\theta (\rho _{\theta }(t))),$ $t\geq 0,$ $\rho
_{\theta }$ being the time-change
\begin{equation*}
\rho _{\theta }(t)=\inf \left\{ u\geq 0:\int_{0}^{u}\exp \left( \alpha
\theta (r)\right) dr>t\right\} ,t\geq 0.
\end{equation*}
Denoting by Leb the Lebesgue measure on $\left( 0,1\right) ,$ we then have
that
\begin{equation}
\text{Leb}\left\{ x\in \left( 0,1\right) :H_{F}(x)>t\right\} \geq \exp
(-\theta (\rho _{\theta }(t))).  \label{223}
\end{equation}
On the other hand, recall that $H_{F}$ is anyway a.s. continuous and
introduce for every $t>0$%
\begin{equation*}
x_{t}:=\inf \left\{ x:H_{F}(x)=t\right\} ,
\end{equation*}
so that $x<x_{t}\Rightarrow H_{F}(x)<t.$ Hence $x_{t}\leq $Leb$\left\{ x\in
\left( 0,1\right) :H_{F}(x)<t\right\} $ and this yields, together with $%
\left( \ref{223}\right) ,$ to
\begin{equation*}
x_{t}\leq 1-\exp (-\theta (\rho _{\theta }(t))\text{ a.s. for every }t\geq 0.
\end{equation*}
Now suppose that $H_{F}$ is a.s. H\"{o}lder-continuous of order $\gamma .$
The previous inequality then gives
\begin{equation}
t=H_{F}(x_{t})\leq Cx_{t}^{\gamma }\leq C\left( \theta (\rho _{\theta
}(t))\right) ^{\gamma }  \label{215}
\end{equation}
so that it is sufficient to study the behavior of $\theta (\rho _{\theta
}(t))$ as $t\rightarrow 0$ to obtain an upper bound for $\gamma .$ It is
easy that $\rho _{\theta }(t)\sim t$ as $t\downarrow 0$, so we just have to
focus on the behavior of $\theta (t)$ as $t\rightarrow 0.$ By \cite[Theorem
III.4.9]{bertLevy96}, for every $\delta >1$, $\lim_{t\rightarrow 0}\left(
\theta (t)/t^{\delta }\right) =0$ as soon as $\int_{0}^{1}\overline{\pi
_{\theta }}(t^{\delta })dt<\infty ,$ where $\overline{\pi _{\theta }}%
(t^{\delta })=\int_{t^{\delta }}^{\infty }\pi _{\theta }(dx).$ To see when
this quantity is integrable near $0$, remark first that
\begin{equation*}
\overline{\pi _{\theta }}(u)=\overline{\pi _{\theta }}(1)+\int_{u}^{1}e^{-x}%
\nu (-\log s_{1}\in dx)\text{ when }u<1,
\end{equation*}
(since $s_{i}\leq 1/2$ for $i\geq 2$) and second that
\begin{equation*}
\int_{u}^{1}e^{-x}\nu (-\log s_{1}\in dx)\leq \nu (s_{1}<e^{-u}).
\end{equation*}
Hence,
\begin{equation*}
\int_{0}^{1}\overline{\pi _{\theta }}(t^{\delta })dt\leq \overline{\pi
_{\theta }}(1)+\int_{0}^{1}\nu (s_{1}<e^{-t^{\delta }})dt
\end{equation*}
and by definition of $\vartheta _{\text{up}}$ this last integral is finite
as soon as $1/\delta >\vartheta _{\text{up}}.$ Thus $\lim_{t\rightarrow
0}\left( \theta (t)/t^{\delta }\right) =0$ for every $\delta <1/\vartheta _{%
\text{up}}$ and this implies, recalling $\left( \ref{215}\right) ,$ that $%
\gamma \delta <1$ for every $\delta <1/\vartheta _{\text{up}}.$ Which gives $%
\gamma \leq \vartheta _{\text{up}}.$

It remains to prove that $H_{F}$ cannot be H\"{o}lder-continuous with an
exponent $\gamma $ larger than $\left| \alpha \right| /\varrho .$ This is
actually a consequence of the results we have on the minoration of \textrm{d}%
$\mathrm{im\,}_{\mathcal{H}}(\mathcal{T}_{F}).$ Indeed, recall the
definition of the function $\tilde{H}_{F}:\left( 0,1\right) \rightarrow
\mathcal{T}_{F}$ introduced Section 4.1 and in particular that for $0<x<y<1$
\begin{equation*}
d\left( \tilde{H}_{F}(x),\tilde{H}_{F}(y)\right)
=H_{F}(x)+H_{F}(y)-2\inf_{z\in \left[ x,y\right] }H_{F}(z),
\end{equation*}
which shows that the $\gamma $-H\"{o}lder continuity of $H_{F}$ implies that
of $\tilde{H}_{F}$. It is easy and well known that since $\tilde{H}%
_{F}:\left( 0,1\right) \rightarrow \mathcal{T}_{F},$ the $\gamma $%
-H\"{o}lder continuity of $\tilde{H}_{F}$ leads to $\dim_{\cal H}
(\mathcal{T}_{F})\leq \dim_{\cal H}((0,1))/\gamma
=1/\gamma $. 
Hence $H_{F}$ cannot be H\"{o}lder-continuous with an order $\gamma >1/%
\dim_{\cal H}(\mathcal{T}_{F}).$ Recall then that $\mathrm{%
dim\,}_{\mathcal{H}}(\mathcal{T}_{F})\geq \varrho /\left| \alpha \right| .$
Hence $H_{F}$ cannot be H\"{o}lder-continuous with an order $\gamma >\left|
\alpha \right| /\varrho .$

\subsection{Height process of the stable tree}

\label{heightstable}

To prove Corollary \ref{C2}, we will check that $\nu _{-}(1-s_{1}>x)\sim
Cx^{1/\beta -1}$ for some $C>0$ as $x\downarrow 0$, where $\nu _{-}$ is the
dislocation measure of the fragmentation $F_{-}$ associated to the stable ($%
\beta $) tree. In view of Theorem \ref{T4} this is sufficient, since the
index of self-similarity is $1/\beta -1$ and $\int_{S}\left(
s_{1}^{-1}-1\right) \nu ($d$\mathbf{s})<\infty ,$ as proved in Sect.\ \ref
{checkst}. Recalling the definition of $\nu _{-}$ in Sect.\ \ref{checkst}
and the notations therein, we want to prove
\begin{equation*}
E\left[ T_{1}\mathbf{1}_{\{1-\Delta _{1}/T_{1}>x\}}\right] \sim Cx^{1/\beta
-1}\qquad \mbox{as }x\downarrow 0
\end{equation*}
Using the above notations, the quantity on the left can be rewritten as
\begin{equation*}
E\left[ (\Delta _{1}+T_{1}^{(\Delta _{1})})\mathbf{1}_{\left\{
T_{1}^{(\Delta _{1})}/(\Delta _{1}+T_{1}^{(\Delta _{1})})>x\right\} }\right]
=E\left[ \Delta _{1}(1+\Delta _{1}^{-1}T_{1}^{(\Delta _{1})})\mathbf{1}%
_{\left\{ \Delta _{1}^{-1}T_{1}^{(\Delta _{1})}>x(1-x)^{-1}\right\} }\right]
.
\end{equation*}
Recalling the law of $\Delta _{1}$ and the fact that $v^{-1}T_{1}^{(v)}$ has
same law as $T_{v^{-1/\beta }}^{(1)}$, this is
\begin{equation*}
c\int_{0}^{\infty }\mathrm{d}v\,v^{-1/\beta }e^{-c\beta v^{-1/\beta }}E\left[
(1+T_{v^{-1/\beta }}^{(1)})\mathbf{1}_{\left\{ T_{v^{-1/\beta
}}^{(1)}>x(1-x)^{-1}\right\} }\right] .
\end{equation*}
By \cite[Proposition 28.3]{Sato}, since $T^{(1)}$ and $T$ share the same
L\'{e}vy measure on a neighborhood of $0$, $T_{v}^{(1)}$ admits a continuous
density $q_{v}^{(1)}(x),x\geq 0$ for every $v>0$. We thus can rewrite the
preceding quantity as
\begin{equation*}
c\int_{0}^{\infty }\frac{\mathrm{d}v}{v^{1/\beta }}e^{-c\beta v^{-1/\beta
}}\int_{x/(1-x)}^{\infty }\!\!\!(1+u)q_{v^{-1/\beta }}^{(1)}(u)\mathrm{d}%
u=c\beta \int_{x/(1-x)}^{\infty }\!\!\!\mathrm{d}u(1+u)\int_{0}^{\infty }%
\frac{\mathrm{d}w}{w^{\beta }}e^{-c\beta w}q_{w}^{(1)}(u)
\end{equation*}
by Fubini's theorem and the change of variables $w=v^{-1/\beta }$. The
behavior of this as $x\downarrow 0$ is the same as that of $c\beta J(x)$
where $J(x)=\int_{x}^{\infty }\mathrm{d}uj(u)$, and where $%
j(u)=\int_{0}^{\infty }\mathrm{d}w\,w^{-\beta }e^{-c\beta w}q_{w}^{(1)}(u)$.
Write $\mathcal{J}(x)=\int_{0}^{x}J(u)\mathrm{d}u$ for $x\geq 0$, and
consider the Stieltjes-Laplace transform $\hat{\mathcal{J}}$ of $\mathcal{J}$
evaluated at $\lambda \geq 0$:
\begin{eqnarray*}
\hat{\mathcal{J}}(\lambda )=\int_{0}^{\infty }e^{-\lambda u}J(u)\mathrm{d}u
&=&\lambda ^{-1}\int_{0}^{\infty }(1-e^{-\lambda u})j(u)\mathrm{d}u \\
&=&\lambda ^{-1}\int_{0}^{\infty }\frac{\mathrm{d}w}{w^{\beta }}e^{-c\beta
w}\int_{0}^{\infty }\mathrm{d}uq_{w}^{(1)}(u)(1-e^{-\lambda u}) \\
&=&\lambda ^{-1}\int_{0}^{\infty }\frac{\mathrm{d}w}{w^{\beta }}e^{-c\beta
w}(1-e^{-w\Phi ^{(1)}(\lambda )})
\end{eqnarray*}
where as above $\Phi ^{(1)}(\lambda )=c\int_{0}^{1}u^{-1-1/\beta
}(1-e^{-\lambda u})\mathrm{d}u$. Integrating by parts yields
\begin{eqnarray*}
\hat{\mathcal{J}}(\lambda ) &=&\frac{\lambda ^{-1}}{\beta -1}%
\int_{0}^{\infty }\frac{\mathrm{d}w}{w^{\beta -1}}e^{-c\beta w}((c\beta
+\Phi ^{(1)}(\lambda ))e^{-w\Phi ^{(1)}(\lambda )}-c\beta ) \\
&=&\lambda ^{-1}\frac{\Gamma (2-\beta )}{\beta -1}((c\beta +\Phi
^{(1)}(\lambda ))^{\beta -1}-(c\beta )^{\beta -1})
\end{eqnarray*}
Is is easy by changing variables in the definition of $\Phi ^{(1)}$ that $%
\Phi ^{(1)}(\lambda )\sim C\lambda ^{1/\beta }$ as $\lambda \rightarrow
\infty $ for some $C>0$, so finally we obtain that $\hat{\mathcal{J}}%
(\lambda )\sim C\lambda ^{-1/\beta }$ as $\lambda \rightarrow \infty $ for
some other $C>0$. Since $\mathcal{J}$ is non-decreasing, Feller's version of
Karamata's Tauberian theorem \cite[Theorem 1.7.1']{bgt} gives $\mathcal{J}%
(x)\sim Cx^{1/\beta }$ as $x\downarrow 0$, and since $J$ is monotone, the
monotone convergence theorem \cite[Theorem 1.7.2b]{bgt} gives $J(x)\sim
\beta ^{-1}Cx^{1/\beta -1}$ as $x\downarrow 0$, as wanted.


\end{document}